\documentclass[12pt,twoside]{amsart}
\usepackage{amsmath}
\usepackage{amssymb}
\usepackage{mathrsfs}         

\usepackage{tikz-cd}
\tikzcdset{
  arrow style=tikz,
  diagrams={>={stealth}}
}
\usepackage{fontawesome}
\usepackage{stix}

\usepackage[colorlinks=true, urlcolor=blue,bookmarks=true,bookmarksopen=true, citecolor=blue]{hyperref}
\usepackage{slashed} 

\usepackage{graphicx}

\numberwithin{equation}{section}

\newtheorem{theorem}{Theorem}[section]

\newtheorem{corollary}[theorem]{Corollary}

\newtheorem{lemma}[theorem]{Lemma}
\newtheorem{prop}[theorem]{Proposition}

\theoremstyle{definition}
\newtheorem{defn}[theorem]{Definition}
\newtheorem{example}[theorem]{Example}
\newtheorem{remark}[theorem]{Remark}

\def \begineq{\begin{equation}}
\def \endeq{\end{equation}}

\def \bb{\mathbb}

\def \mc{\mathcal}
\def \mf{\mathfrak}
\def \ms{\mathscr}

\renewcommand{\tilde}{\widetilde}

\newcommand{\aair}[1]{\left\langle#1\right\rangle}

\newcommand{\half}[1]{\frac{#1}{2}}
\newcommand{\norm}[1]{{\left|\!\left| #1 \right|\!\right|}}

\newcommand{\bair}[1]{\left\{#1\right\}}
\newcommand{\pair}[1]{\left(#1\right)}
\newcommand{\sair}[1]{\left[#1\right]}
\newcommand{\se}[1]{C^\infty\left(#1\right)}

\def \dd{{\bb d}}
\def \uu{{\bb u}}
\def \vv{{\bb v}}
\def \ww{{\bb w}}

\def \CC{{\bb{C}}}

\def \GG{{\bb G}}

\def \JJ{{\bb{J}}}
\def \KK{{\bb{K}}}
\def \LL{{\bb{L}}}

\def \PP{{\bb{P}}}

\def \RR{{\bb{R}}}

\def \TT{{\bb{T}}}

\def \ZZ{{\bb{Z}}}

\def \TTM{{\bb{T}M}}

\def \BBC{{\mc B}}

\def \LLC{{\mc L}}

\def \WWC{{\mc W}}

\def \DDS{{\ms D}}

\def \JJS{{\ms J}}

\def \RRS{{\ms R}}

\def \({\left(}
\def \){\right)}
\def \<{\left\langle}
\def \>{\right\rangle}

\def \bar{\overline}

\def \cdiff{{\psi}}

\def \deg{\mathrm{deg}}
\def \dsum{\mathop{\oplus}}

\def \gcon{{\mc D}}
\def \gcurv{{\mc F}}

\def \grie{{\mc R}}
\def \gsca{{\mc S}}
\def \half{\frac{1}{2}}
\def \hat{\widehat}
\def \inter{\mathop{\cap}}
\def \into{\hookrightarrow}
\def \inverse{^{-1}}

\def \laplacian{{\Delta}}

\def \tensor{\otimes}

\def \vargeq{\geqslant}
\def \varleq{\leqslant}

\def \xto{\xrightarrow}

\def \Aut{{\rm Aut}}
\def \bismut{{\phi, \BBC}}
\def \bismutt{{\phi_t, \BBC_t}}

\def \Diff{{\rm Diff}}
\def \End{{\rm End}}

\def \id{{\rm id}}
\def \Id{{\rm Id}}
\def \img{\mathop{\rm img }}

\def \phib{{\phi, b}}
\def \phibt{{\phi_t, b_t}}

\def \rank{\mathop{\rm rank}}
\def \Ric{{\mc Ric}}
\def \Rc{{\mc Rc}}

\def \tr{{\rm tr}}
\def \TTM{{\TT M}}
\def \varphib{{\pair{\phib}}}
\def \Vol{{\rm Vol}}
\def \vol{{\rm vol}}
\def \YM{{\mc {Y\!M}}}

\newcommand{\1}{1\!\!1}

\def \qed{\hfill $\square$ \vspace{0.03in}}

\begin{document}

\title{Differential calculus for generalized geometry and geometric Lax flows}

\author{Shengda Hu}

\email{shu@wlu.ca}
\address{Department of Mathematics, Wilfrid Laurier University, 75 University Ave. West, Waterloo, Canada}

\abstract
 Employing a class of generalized connections, we describe certain differential complices $\pair{\tilde \Omega^*_\TT\pair{M}, \tilde \dd^\TT}$ constructed from $\wedge^* \TT M$ and study some of their basic properties, where $\TT M = T M \dsum T^*M$ is the generalized tangent bundle on $M$. A number of classical geometric notions are extended to $\TT M$, such as the curvature tensor for a generalized connection. In particular, we describe an analogue to the Levi-Civita connection when $\TT M$ is endowed with a generalized metric and a structure of exact Courant algebroid. We further describe in generalized geometry the analogues to the Chern-Weil homomorphism, a Weitzenb\"ock identity, the Ricci flow and Ricci soliton, the Hermitian-Einstein equation and the degree of a holomorphic vector bundle. Furthermore, the Ricci flows are put into the context of geometric Lax flows, which may be of independent interest.
\endabstract



\maketitle


\section{Introduction}\label{sec:intro}

In generalized geometry \`a la Hitchin \cite{Hitchin:genCalabiYau:03}, over a smooth manifold $M$ of real dimension $n$, the bundle $\TT M := TM \dsum T^*M$ is considered the analogue of the classical tangent bundle $TM$. It fits into the natural exact sequence
\[0 \to T^*M \into \TT M \xto{\pi} TM \to 0\]
and is endowed with the natural pairing
\[\aair{x, y} = \aair{X+\xi, Y+\eta} := \half \pair{\iota_X \eta + \iota_Y \xi}\]
where $x, y \in \se{\TT M}$ and $X, Y \in \se{TM}$, $\xi, \eta \in \se{T^*M}$ are their respective components. The dual of $\TT M$ can be identified with itself under the pairing $2\aair{,}$. Well-known geometric structures on $\TT M$ such as generalized complex, Riemannian, Hermitian, K\"ahler structures and generalized connections are natural extensions of the corresponding classical notions on $TM$. There are by now many references in the literature, including the pioneering works by Gualtieri \cite{Gualtieri:gencomplex:04, Gualtieri:Poisson:10, Gualtieri:gencomplex:11, Gualtieri:genKahler:14} on the subjects.

We show that the analogue can be pushed further, leading to coherent extensions of known geometric notions, as well as Lax forms for some geometric flows. For simplicity, we will restrict our considerations to compact connected orientable smooth manifolds without boundary, while aside from cohomology computations, most descriptions are of a local nature. It starts with a differential complex constructed with $\TT M$ in place of $TM$, as a quotient of $\Omega_\TT^*\pair{M} := \se{\wedge^*\TT M}$, whose differential is defined from a generalized connection on $\TT M$.
Let $V$ be a vector bundle, as in \cite{Gualtieri:Poisson:10}, a generalized connection on $V$ is a derivation
\[\gcon: \se{V} \to \se{\TT M \tensor V}\]
such that for $f \in \se{M}, v \in \se{V}$, $x \in \se{\TT M}$ and $X = \pi\pair{X}$
\[\gcon_x\pair{f v} = X\pair{f} v + f \gcon_x v\]
Given a generalized connection $\gcon^\TT$ on $\TT M$, it is \emph{$TM$-torsion free} if for all $x, y \in \se{\TT M}$
\begin{equation}\label{eq:TMtorsionfree}
 \pi\pair{x \diamond_\TT y} = \sair{\pi\pair{x}, \pi\pair{y}}
\end{equation}
where the \emph{$\gcon^\TT$-diamond bracket $\diamond_\TT$} is the skew symmetrization of $\gcon^\TT$:
\begin{equation}\label{eq:diamondbracketintro}
 x \diamond_\TT y := \gcon^\TT_x y - \gcon^\TT_y x
\end{equation}
The first observation is a differential complex obtained from such a generalized connection.
\begin{theorem}[\S \ref{subsec:TTMforms}]\label{thm:differentialcomplexintro}
 Consider the derivation $\dd^\TT : \Omega_\TT^*\pair{M} \to \Omega_\TT^{*+1}\pair{M}$ defined by
 \[\pair{\dd^\TT\theta}(x_0, x_1, \ldots, x_k) := \sum_{i}(-1)^i\pair{\gcon^\TT_{x_i} \theta}\pair{x_0, x_1, \ldots, \hat{x_i}, \ldots, x_k}\]
 where $\theta \in \Omega_\TT^k\pair{M}$ and $x_j \in \se{\TT M}$. Then $\dd^\TT \circ \dd^\TT$ is tensorial iff $\gcon^\TT$ is $TM$-torsion free, in which case, the quotient $\tilde \Omega_\TT^*\pair{M}$ of $\Omega_\TT^*\pair{M}$ by the image of $\dd^\TT \circ \dd^\TT$ is a differential complex with the induced derivation $\tilde \dd^\TT$.
\end{theorem}
\noindent
The cohomology of the resulting complex $\pair{\tilde \Omega_\TT^*\pair{M}, \tilde \dd^\TT}$ is denoted $\tilde H_\TT^*\pair{M}$, and the projection $\pi$ induces a natural map $\tilde \pi^*: H^*\pair{M} \to \tilde H_\TT^*\pair{M}$, where $H^*\pair{M}$ is the classical de Rham cohomology.

Suppose that $\TT M$ is endowed with a generalized Riemannian structure $\GG^b$ \cite{Gualtieri:gencomplex:04, Gualtieri:gencomplex:11}, defined by a Riemannian metric $g$ and a $2$-form $b \in \Omega^2\pair{M}$
\begin{equation}\label{eq:genmetric}
 \GG^b\pair{x, y} := \half \sair{g\pair{X, Y} + g\inverse\pair{\xi - \iota_X b, \eta - \iota_Y b}}
\end{equation}
The \emph{$\GG^b$-metric connections} are generalized connections on $\TT M$ that preserve $\aair{,}$ and $\GG^b$, among which the $TM$-torsion free ones are classified by $\phi \in \Omega^3\pair{M}$ and denoted $\gcon^\phib$.

The \emph{$\pair{\phib}$-connections} $\gcon^\phib$ can be seen as analogues of the Levi-Civita connections in this context, since they are also compatible (Definition \ref{defn:diamondmetriccompatible}) with the corresponding almost Dorfman brackets $*_\gamma$ defined by $\gamma = \phi - db$ :
\begin{equation}\label{eq:dorfmanbracket}
 (X + \xi) *_{\gamma} (Y+\eta) = [X, Y] + \LLC_X \eta - d\iota_Y \xi + \iota_Y\iota_X \gamma
\end{equation}
When $d\gamma = 0$, the bracket $*_\gamma$ is a \emph{Dorfman bracket}.
\begin{theorem}[Theorem \ref{thm:levicivita}, \S \ref{subsec:cohomology}]\label{thm:levicivitaintro}
 Let $\gamma \in \Omega^3\pair{M}$ be a closed $3$-form and $\phi = \gamma + db$. The $\varphib$-connection $\gcon^\phib$ associated with $\GG^b$ is the unique $TM$-torsion free $\GG^b$-metric connection on $\TT M$ that is metric compatible with a Dorfman bracket $*_\gamma$. The natural map $\tilde \pi^*: H^*\pair{M} \to \tilde H_\phib^*\pair{M}$ is injective. Moreover, $\tilde H_\phib^{2n}\pair{M} \cong \RR$.
\end{theorem}

The differential calculus thus established on $\TT M$ leads to a natural definition of curvature for any generalized connection $\gcon$ on any vector bundle $V$:
\[\gcurv^\TT\pair{\gcon} := \dd^\TT_\gcon \circ \gcon \in \Omega_\TT^2\pair{\End\pair{V}}\]
where $\dd^\TT_\gcon$ is the extension of $\dd^\TT$ to $\se{\TT M \tensor V}$ as usual. In terms of covariant derivatives, $\gcurv^\TT\pair{\gcon}$ is given by
\begin{equation}\label{eq:gencurvatureintro}
 \gcurv^\TT_{x, y}\pair{\gcon} v := \pair{\gcon_x \gcon_y - \gcon_y \gcon_x - \gcon_{x \diamond_\TT y}} v
\end{equation}
where $x, y \in \se{\TT M}$ and $v \in \se{V}$. A side-effect of this is that the resulting curvature tensor $\gcurv^\TT$ now depends on the generalized connection $\gcon^\TT$ on $\TT M$. Nonetheless, passing to the quotient $\tilde \Omega_\TT^*\pair{\End\pair{V}}$, the Chern-Weil homomorphism naturally extends.
\begin{theorem}[\S \ref{subsec:chernweil}]\label{thm:chernweilintro}
 Any invariant polynomial of the $\gcon^\TT$-curvature $\gcurv^\TT\pair{\gcon}$ defines a class in $\tilde H_\TT^*\pair{M}$, which coincides with the image under $\tilde \pi^*$ of the corresponding classical characteristic class in $H^*\pair{M}$.
\end{theorem}

The $\varphib$-curvature for $\gcon^\phib$ itself, denoted $\grie^\phib$, can be regarded as the analogue to the Riemannian curvature in this context. The $\pair{\phib}$-Ricci curvature $\Ric^\phib$, as well as the corresponding scalar curvature are defined via the usual contractions of $\grie^\phib$. Associated to the connection $\gcon^\phib$, further operators, including the analogues to the Bochner and Hodge Laplacians, can be defined on the tensor bundles. An analogue to the Weitzenb\"ock identity holds, which states that the Laplacians differ by the $\pair{\phib}$-Ricci curvature.
\begin{theorem}[Theorem \ref{thm:hodgelaplacianbochner}]\label{thm:Weitzenbockintro}
 On $\Omega^*_\TT\pair{M}$, let $\laplacian^\phib$ be the $\gcon^\phib$-Hodge Laplacian and $\laplacian_{\gcon^\phib}$ be the $\gcon^\phib$-Bochner Laplacian, then
 \[\laplacian^\phib = \laplacian_{\gcon^\phib} + \GG^b\Ric^\phib\GG^b\]
\end{theorem}

Consider next a generalized complex manifold $\pair{M, \gamma; \JJ}$. A generalized connection $\gcon^\TT$ is \emph{$\JJ$-compatible with $*_\gamma$} if the diamond bracket $\diamond_\TT$ coincides with $*_\gamma$ when restricted to sections of the same eigenbundle of $\JJ$. For such $\gcon^\TT$, $\dd^\TT$ admits the natural decomposition
\[\dd^\TT = \partial_\JJ^\TT + \bar \partial_\JJ^\TT\]
according to the types with respect to $\JJ$. Together with a generalized metric $\GG^b$ commuting with $\JJ$, we have a generalized Hermitian manifold $\pair{M, \gamma; \GG^b, \JJ}$, which corresponds classically to an almost bi-Hermitian structure $\pair{M, \gamma; g, I_\pm; b}$, where $g$ is Hermitian with respect to both almost complex structures $I_\pm$. We show that the $\JJ$-compatibility of $\gcon^\phib$ with $*_\gamma$ is a generalized K\"ahler condition, equivalent to the condition given in \cite{Gualtieri:Poisson:10}.
\begin{theorem}[Theorem \ref{thm:genkahlercondition}]\label{thm:genkahlerconditionintro}
 On a generalized Hermitian manifold $\pair{M, \gamma; \GG^b, \JJ}$, let $\phi = \gamma + db$. Then it is a generalized K\"ahler manifold iff $\gcon^\phib$ is $\JJ$-compatible with $*_\gamma$.
\end{theorem}
\noindent
In terms of $I_\pm$, the $\JJ$-compatibility is equivalent to $\nabla^{\pm \phi} I_\pm = 0$.
 On a generalized K\"ahler manifold, $I_\pm$ are integrable. Working with $\gcon^\phib$, we recover a well-known result obtained via holomorphic reduction \cite{Gualtieri:genKahler:14}, that the $I_\pm$-(anti)holomorphic tangent bundles on a generalized K\"ahler manifold carry natural $I_\mp$-holomorphic structures respectively.

For a $\JJ$-holomorphic Hermitian vector bundle $\pair{V, \bar\partial_\JJ, h}$, the notion of Chern connection extends naturally. Over a generalized Hermitian manifold $\pair{M, \gamma; \GG^b, \JJ}$, there is a natural contraction $\Lambda_{\JJ_-}$ on $\Omega_\TT^2\pair{M}$, where for $x, y \in \se{\TT M}$
\[\Lambda_{\JJ_-}\pair{x \wedge y} := \GG^b\pair{\JJ x, y}\]
Given a $\gamma$-$\JJ$-connection $\gcon^\TT$ (Definition \ref{defn:diamondgcplxcompatible}) on $\TT M$, analogous to the classical case (L\"ubke-Teleman \cite{LubkeTeleman:KobayashiHitchin:95}), we propose the following \emph{$\gcon^\TT$-$\JJ$-Hermitian-Einstein equation} for the Hermitian metric $h$:
\[\sqrt{-1}\Lambda_{\JJ_-}\pair{\gcurv^{\TT, C}(V)} = 2 c \Id_V\]
where $c \in \RR$ and $\gcurv^{\TT, C}$ denotes the $\gcon^\TT$-curvature of the generalized Chern connection. The contraction also leads to the notion of degree for a $\JJ$-holomorphic Hermitian vector bundle. The degree is independent to the choice of Hermitian metric on $V$ if the generalized Hermitian manifold is \emph{$\gcon^\TT$-$\JJ$-Gauduchon}, i.e. for all $f \in \se{M}$
\[\int_M \Lambda_{\JJ_-}\pair{\partial_\JJ^\TT \bar \partial_\JJ^\TT f} d\vol_g = 0\]
For such manifolds, the notions of slope and stability naturally extend and one should expect a version of Kobayashi-Hitchin correspondence to hold. On a generalized K\"ahler manifold, these notions relate to their classical counterparts, in particular, the $\JJ$-Hermitian-Einstein equation is equivalent to an equation proposed by Hitchin \cite{Hitchin:GenHoloBundles:11} (Remark in \S 3.3 \emph{ibid.}).
\begin{theorem}[\S \ref{subsec:JHermitianEinstein}]\label{thm:genkahlerheintro}
 Let $\pair{M, \gamma; \GG^b, \JJ}$ be a $\JJ$-Gauduchon generalized K\"ahler manifold and $\omega_\pm \in \Omega^2\pair{M}$ be the K\"ahler forms for $I_\pm$ respectively. Let $\pair{V, \partial_\JJ}$ be a $\JJ$-holomorphic vector bundle, then the $\JJ$-Hermitian-Einstein equation is equivalent to
 \[\dfrac{\sqrt{-1}}{2}\pair{F_+^C\pair{V} \wedge \omega_+^{m-1} + \pair{-1}^\varepsilon F_-^C\pair{V} \wedge \omega_-^{m-1}} = c \pair{m-1}! \Id_V d\vol_g\]
 where $F^C_\pm$ are the classical Chern curvatures with respect to $I_\pm$, $\varepsilon = 0$ if $I_\pm$ induce the same orientation on $TM$ and $\varepsilon = 1$ otherwise. The $\GG^b$-degree of $\pair{V, h}$ is computed by
 \[\deg_{\GG^b}\pair{V, h} = \half \pair{\deg_+\pair{V, h} + \deg_-\pair{V, h}}\]
 where $\deg_\pm$ are the degrees with respect to the $I_\pm$-holomorphic structure.
\end{theorem}

Geometric flows such as the mean curvature flow (Brakke \cite{Brakke:meancurvature:78}) and the Ricci flow (Hamilton \cite{Hamilton:ricciflow:82}) are very important in understanding smooth manifolds and structures associated to them. In generalized geometry, it is natural to consider flows involving the structures on $\TT M$ such as the generalized metrics or generalized complex structures. A generalized metric $\GG^b$ defines a bundle automorphism of $\TT M$ via the pairing $\aair{,}$:
\begin{equation}\label{eq:metricmap}
 \aair{\GG^b x, y} := \GG^b\pair{x, y}
\end{equation}
which enjoys the following properties
\begin{equation}\label{eq:metricmapproperties}
 \pair{\GG^b}^2 = \1 \text{ and } \aair{\GG^b x, y} = \aair{x, \GG^b y}
\end{equation}
For a smooth family of generalized metrics $\GG^{b_t}_t$, differentiating the above identites shows that $\KK_t = \dfrac{d}{dt}\GG^{b_t}_t$ is skew with respect to $\GG^{b_t}_t$ and symmetric with respect to $\aair{,}$, i.e. it satisfies
\begin{equation}\label{eq:metricmappropertiesdiff}
 \GG^{b_t}_t \KK_t + \KK_t \GG^{b_t}_t = 0 \text{ and } \aair{\KK_t x, y} = \aair{x, \KK_t y}
\end{equation}
Hence, a flow of generalized metrics can be defined from geometric quantities satisfying the above identities. Similar identities can be derived for a smooth family of generalized complex structures. For any smooth family $\LL_t : \TT M \to \TT M$, the commutator $\sair{\LL_t, \GG^{b_t}_t}$ satisfies the first identity in \eqref{eq:metricmappropertiesdiff}, which lends to the consideration of Lax pairs (Lax \cite{Lax:laxpair:68}).

In this context, we generally assume that there is a fixed closed $3$-form $\gamma$ such that $\phi = \gamma + db$ holds for all $\varphib$ under consideration. Hence the flow generated can be regarded as over a fixed Courant algebroid, defined by the Dorfman bracket $*_\gamma$. For instance, any $\TT M$-form $\theta \in \Omega_\TT^2\pair{M}$ defines a \emph{$2$-form Lax flow} via the induced map $\theta_t : \TT M \to \TT M$. The action of generalized symmetries of $\TT M$ on a generalized metric can be represented as Lax flows generated by the $\varphib$-curvatures of a Hermitian line bundle.
\begin{theorem}[Theorem \ref{thm:linebundlemorphism}]\label{thm:genmorphismaslaxflowintro}
 Consider a family of unitary generalized connections $\bair{\gcon_t}$ on a Hermitian line bundle $V$. Let $\phi_t = \gamma + db_t$, then the Lax flow
 \[\dfrac{d}{dt} \GG^{b_t}_t = \sair{\gcurv^{\phibt}\pair{\gcon_t}, \GG^{b_t}_t}\]
 corresponds to the push-forward of an initial generalized metric $\GG^{b}$ by a family of generalized diffeomorphisms of the exact Courant algebroid defined by $*_\gamma$.
\end{theorem}
\noindent
Here we recall (\cite{Gualtieri:gencomplex:04}, Hu-Uribe \cite{HuUribe:gensymmetry:07}) that a generalized diffeomorphism for the Courant algebroid structure on $\TT M$ defined by $*_\gamma$ is a pair $\pair{\lambda, B}$, where $\lambda \in \Diff\pair{M}$ is a diffeomorphism and $B \in \Omega^2\pair{M}$, such that
\[\lambda^*\gamma = \gamma - dB\]
whose action on $\TT M$ by push-forward is given by
\[\pair{\lambda, B}_*\pair{X + \xi} := \lambda_*\pair{X + \iota_X B + \xi}\]
The action of $B \in \Omega^2\pair{M}$ is also denoted $e^B$:
\[e^B\pair{X + \xi} := X + \iota_XB + \xi\]
We comment that the diffeomorphism components of the generalized diffeomorphisms in Theorem \ref{thm:genmorphismaslaxflowintro} arise from the vector components of the generalized connections on $V$, which are $0$ when $\gcon_t$ are liftings of classical connections on $V$.

Even though the Bianchi identites do not hold for $\grie^\phib$ in general, it turns out that the $\varphib$-Ricci tensor $\Rc^\phib$ is symmetric. The corresponding Lax flow is the \emph{Ricci Lax flow}
\begin{equation}\label{eq:laxricciintro}
 \dfrac{d}{dt} \GG^{b_t}_t = \sair{\Ric^{\phibt}, \GG^{b_t}_t}
\end{equation}
which defines the Ricci flow in this context. We note that in \eqref{eq:laxricciintro}, only the components of $\Rc^\phib$ that involve sections from different eigenbundles of $\GG^b$ contribute, similar to the Ricci flows for Courant algebroids considered in \v{S}evera-Valach \cite{SeveraValach:ricciflow:16, SeveraValach:typeIIsupergravity:18} and Streets \cite{Streets:genTdualityRenorm:13}.
\begin{theorem}[Theorem \ref{thm:genricciflow}]\label{thm:genricciflowintro}
 The Ricci Lax flow equation above is equivalent to the generalized Ricci flow equation
 \begin{equation*}
  \begin{cases}
   \dfrac{d}{dt} g_t & = -2 Rc_{t} + \dfrac{1}{2}\pair{\gamma + db_t}^2\\
   \dfrac{d}{dt} b_t & = - d^*\pair{\gamma + db_t}
  \end{cases}
 \end{equation*}
\end{theorem}
\noindent
The above system has been studied in mathematics and physics literatures, for instance, Garcia-Fernandez-Streets \cite{GarciaFernandezStreets:generalizedRicciflow:20}, and references therein. We note that an equation similar to \eqref{eq:laxricciintro} appeared in \cite{Streets:genTdualityRenorm:13} citing communication with Gualtieri (reference [11] \emph{ibid.}) and the equivalence with the generalized Ricci flow appeared in \cite{GarciaFernandezStreets:generalizedRicciflow:20} (Remark 4.8 \emph{ibid.}), via a somewhat different approach to the notion of generalized Ricci curvature.

A few more examples of geometric Lax flows are given. Conformal deformations of the Riemannian metric $g$ can be represented as a Lax flow, of the form
\[\dfrac{d}{dt}\GG^{b_t}_t = \sair{\theta_t, \GG^{b_t}_t}\]
where $\theta_t \in \Omega_\TT^2\pair{M}$ are $\GG^{b_t}_t$-conformal $\TT M$-forms (Definition \ref{defn:conformalform}). The \emph{Ricci soliton equation} (Streets \cite{Streets:connRicciFlow:08} and references therein) takes the following form (Definition \ref{defn:laxriccisoliton})
\[ \sair{\Ric^{\phibt} - \dd^\phibt u_t - \theta_t, \GG^{b_t}_t} = 0\]
where the $\GG^{b_t}_t$-conformal $\TT M$-forms $\theta_t$ have constant conformal weights. We also see that the classical K\"ahler-Ricci flow can be recast as a geometric Lax flow (\S \ref{subsec:bismutlaxricci}). 

We expect that many classical constructions should admit natural extensions to $\TT M$ via the differential calculus developed here. For instance, the notions of second fundamental form and mean curvature for submanifolds extend naturally. Spinors, which are behind the notion of $\JJ_-$-contraction in Definition \ref{defn:Jminuscontraction}, relate the geometry on $\TT M$ back to $\Omega^*\pair{M}$, which, in particular, lead to the canonical line of a generalized (almost) complex structure \cite{Gualtieri:gencomplex:04, Gualtieri:gencomplex:11} as well as the notion of scalar curvature in generalized K\"ahler geometry (Goto \cite{Goto:scalarcurvature:16, Goto:scalarcurvature:21}, Wang \cite{Wang:ToricGenKahlerIII:19}). Functionals involving curvatures, such as the Yang-Mills functional, can be extended (\S \ref{subsec:yangmills}) and lead to natural questions on extremal / critical (generalized) connections / metrics with respect to them. Explicit examples such as compact Lie groups (\cite{Goto:scalarcurvature:21}, Hu \cite{Hu:LiegenKahler:15}) could provide further insights into understanding these extensions. It should be worth exploring the interaction of the Riemannian, the complex and the Poisson geometric methods in generalized Hermitian geometry.
Equations in Lax form admit geometric interpretations (Griffiths \cite{Griffiths:LaxEquation:85}), it would be interesting to understand if it provides new perspective for the related geometric flows.

We briefly give some further descriptions of the structure of the paper.
In \S \ref{sec:TMtorsionless}, we set up the differential calculus on $\TT M$ and compute in \S \ref{subsec:compactLiegroups} the group $\tilde H_{\gamma, 0}^*\pair{G}$ for a semi-simple real Lie group $G$, with the bi-invariant metric and the Cartan $3$-form $\gamma$. The generalized curvature tensors are introduced in \S \ref{sec:phibcurvature}. The rest of the article applies the constructions in various contexts. The analogue to the Riemann curvature is discussed in \S \ref{sec:phibriemannian}, together with the associated Ricci and scalar curvatures, as well as the generalized Bismut connections \cite{Gualtieri:Poisson:10}. In \S \ref{sec:gencomplex}, we apply the differential calculus to generalized complex and Hermitian manifolds. The degree, stability and Hermitian-Einstein equation for a generalized holomorphic bundle over a generalized Hermitian manifold are discussed in \S \ref{sec:Jholobundles}. In the last section \S \ref{sec:laxflows}, we discuss the notion of geometric Lax flows.

\subsection*{Acknowledgement} I would like to thank Spiro Karigiannis for helpful discussions and enlightening questions, and Yucong Jiang for pointing out the references \cite{SeveraValach:ricciflow:16, SeveraValach:typeIIsupergravity:18}. This research was partially supported by the NSERC Discovery Grant RGPIN-2019-05899.

\section{Differential calculus on $\TT M$}\label{sec:TMtorsionless}

Let $V \to M$ be a vector bundle. Recall that a generalized connection on $V$ is a derivation:
\[\gcon: \se{V} \to \se{\TT M \tensor V} : \gcon\pair{f v} = df \tensor v + f \gcon v\]
where $\TT M = TM \dsum T^*M$, $f \in \se{M}$ and $v \in \se{V}$. It is the \emph{lift} of a classical connection $\nabla_0$ on $V$ if
\begin{equation}\label{eq:liftclassical}
 \gcon_x v = \nabla_{0,\pi\pair{x}} v
\end{equation}
for all $x \in \se{\TT M}$ and $v \in \se{V}$. The generalized connections naturally extend to tensor bundles in the standard fashion.

\subsection{$\TT M$-forms}\label{subsec:TTMforms}
Under the pairing $2\aair{,}$, sections of $\wedge^* \TT M$ can be seen as \emph{$\TT M$-forms} and the space of such forms will be suggestively denoted by
\[\Omega_\TT^*\pair{M} := \se{\wedge^* \TTM}\]
Let $\gcon^\TT$ be a generalized connection on $\TT M$. 
The skew-symmetrization of the covariant derivative by $\gcon^\TT$ induces the \emph{$\gcon^\TT$-derivation $\dd^\TT$}. Namely, for $\theta \in \Omega_\TT^k\pair{M}$:
\begin{equation}\label{eq:TTMderivation}
 \begin{split}
  & \pair{\dd^\TT\theta}(x_0, x_1, \ldots, x_k) := \sum_{i}(-1)^i\pair{\gcon^\TT_{x_i} \theta}\pair{x_0, x_1, \ldots, \hat{x_i}, \ldots, x_k}
 \end{split}
\end{equation}
where $x_i \in \se{\TTM}$. For $f \in \se{M}$, $\dd^\TT$ coincides with the usual differential:
\begin{equation}\label{eq:functiondiff}
 \dd^\TT f := df \in \Omega^1_\TT\pair{M}
\end{equation}
Moverover, $\dd^\TT$ is a graded derivation on $\Omega^*_\TT \pair{M}$, i.e. for $\theta_1 \in \Omega^k\pair{M}$ and $\theta_2 \in \Omega^*\pair{M}$:
\begin{equation}\label{eq:TTMgradedderivation}
 \dd^\TT\pair{\theta_1 \wedge \theta_2} = \pair{\dd^\TT\theta_1} \wedge \theta_2 + \pair{-1}^k \theta_1 \wedge \pair{\dd^\TT \theta_2}
\end{equation}

Standard computations lead to
\begin{equation}\label{eq:dsquared}
 \begin{split}
  & \pair{\dd^\TT \circ \dd^\TT \theta} \pair{x_0, x_1, \ldots, x_{k + 1}}\\
  = & \sum_{i<j} \pair{-1}^{i+j} \sair{\tau_T\pair{x_i, x_j}} \theta\pair{x_0, \ldots, \hat x_i, \ldots, \hat x_j, \ldots, x_{k + 1}}\\
   & - \sum_{i < j < \ell} \pair{-1}^{i+j+\ell} \theta\pair{\sair{x_i \diamond_\TT x_j \diamond_\TT x_\ell}, x_0, \ldots, \hat x_i, \ldots, \hat x_j, \ldots, \hat x_\ell, \ldots, x_{k+1}}
 \end{split}
\end{equation}
where $\tau_T$ is the \emph{$TM$-torsion} of $\gcon^\TT$ given in \eqref{eq:TMtorsion} and $\sair{x \diamond_\TT y \diamond_\TT z}$ is the \emph{Jacobiator} of the diamond bracket $\diamond_\TT$ \eqref{eq:diamondbracketintro}:
\begin{equation}\label{eq:jacobiator}
 \sair{x \diamond_\TT y \diamond_\TT z} := \pair{x \diamond_\TT y} \diamond_\TT z + c.p.
\end{equation}
\begin{defn}\label{defn:Ttorsionless}
 Let $\gcon^\TT$ be a generalized connection on $\TT M$. Its \emph{$TM$-torsion} is
 \begin{equation}\label{eq:TMtorsion}
  \tau_T\pair{x, y} := \pi\pair{x \diamond_\TT y} - \sair{\pi\pair{x}, \pi\pair{y}}
 \end{equation}
 where $x, y \in \se{\TT M}$. Then $\gcon^\TT$ is \emph{$TM$-torsion free} if its $TM$-torsion vanishes.
\end{defn}
Let $\grie^\TT$ be the (not necessarily tensorial) curvature operator for $\gcon^\TT$:
\[\grie^\TT_{x,y} z := \gcon^\TT_x \gcon^\TT_y z - \gcon^\TT_y \gcon^\TT_x z - \gcon^\TT_{\gcon^\TT_x y} z + \gcon^\TT_{\gcon^\TT_y x} z\]
where $x, y, z \in \se{\TT M}$, then straightforward rearrangements lead to
\begin{equation}\label{eq:jacobiatorcurvature}
 \sair{x \diamond_\TT y \diamond_\TT z} = - \grie^\TT_{x,y} z - c.p.
\end{equation}
It follows that $\dd^\TT \circ \dd^\TT$ is tensorial when $\gcon^\TT$ is $TM$-torsion free, since in this case $\grie^\TT$ is tensorial by \eqref{eq:curvtensorialproof}. Furthermore, $\gcon^\TT$ being $TM$-torsion free also implies its Jacobiator has values in $T^*M$:
\begin{equation}\label{eq:projectedjacobi}
 \pi\sair{x \diamond_\TT y \diamond_\TT z} = \sair{\sair{\pi\pair{x}, \pi\pair{y}}, \pi\pair{z}} + c.p. = 0
\end{equation}

Any torsion free affine connection $\nabla^T$ on $TM$ lifts to a $TM$-torsion free generalized connection on $\TT M$:
\[\gcon^\TT_x y := \nabla^T_X y = \nabla^T_X Y + \nabla^T_X \eta\]
for $x, y \in \se{\TT M}$ with $x = X + \xi$ and $y = Y + \eta$. The affine space $\ms D\pair{\TT M}$ of generalized connections on $\TT M$ is modeled on the space of bundle homomorphisms:
\[\ms D\pair{\TT M} \cong \bair{A : \TT M \tensor \TT M \to \TT M}\]
in which, the subspace $\ms D_\tau\pair{\TT M}$ of $TM$-torsion free ones is modeled on the subspace of the right hand side consisting of the symmetric ones valued in $T^*M$:
\begin{equation}\label{eq:Ttorionfreespace}
 \ms D\pair{\TT M} \supseteq \ms D_\tau\pair{\TT M} \cong \bair{A : \TT M \odot \TT M \to T^*M}
\end{equation}

The \emph{contraction} by $x \in \se{\TT M}$ is a graded derivation on $\Omega^*_\TT\pair{M}$ defined by
\begin{equation}\label{eq:contractionconvention}
 \iota_x y := 2\<x, y\>
\end{equation}
where $y \in \Omega^1_\TT \pair{M}$. The \emph{Lie derivative} along $x \in \se{\TT M}$ is given by
\begin{equation}\label{eq:genLiederivative}
 \LLC^\TT_x \theta := \iota_x \dd^\TT \theta + \dd^\TT\iota_x \theta
\end{equation}
where $\theta \in \Omega^*_\TT\pair{M}$. In particular, for $f \in \se{M}$ and $X = \pi\pair{x} \in \se{TM}$:
\[\LLC^\TT_x f = Xf\]
Suppose that $\gcon^\TT$ is $TM$-torsion free, then the familiar relations among the operators $\dd^\TT$, $\iota_x$ and $\LLC^\TT_x$ almost hold, up to possible terms involving the Jacobiator, similar to \eqref{eq:dsquared}.
\begin{prop}\label{prop:operatorrelations}
 Let $x, y, z, w \in \se{\TT M}$, $\theta \in \Omega_\TT^k\pair{M}$, $\alpha \in \Omega^1\pair{M}$ and $X = \pi\pair{x}$. Suppose that $\gcon^\TT$ is $TM$-torsion free, then
 \begin{enumerate}
  \item
  $d\alpha = 0 \Longrightarrow \dd^\TT \alpha = 0$;
  \item
  $\sair{x \diamond_\TT y \diamond_\TT z} \in \se{T^*M}$;
  \item
  $X\<y, z\> = \<x \diamond_\TT y, z\> + \<y, \LLC^\TT_x z\>$;
  \item
  $\LLC^\TT_x\iota_y \theta - \iota_y \LLC^\TT_x \theta = \iota_{x\diamond_\TT y} \theta$;
  \item
  $\<\LLC^\TT_x y - x\diamond_\TT y, z\> = \<\gcon^\TT_y x, z\> + \<\gcon^\TT_z x, y\>$;
  \item
  $\<\sair{\LLC^\TT_x, \LLC^\TT_y} z, w\> = \<\LLC^\TT_{x \diamond_\TT y} z, w\> + \iota_{\sair{x \diamond_\TT y \diamond_\TT w}} z$.
  \item
  For $x_1, \ldots, x_{k+1} \in \se{\TT M}$:
  \begin{equation*}
   \begin{split}
    & \pair{\sair{\dd^\TT, \LLC^\TT_x}\theta}\pair{x_1, \ldots, x_{k+1}} \\
    = & \sum_{i<j} \pair{-1}^{i+j+1}\theta\pair{\sair{x \diamond_\TT x_i \diamond_\TT x_j}, x_1, \ldots, \hat x_i, \ldots, \hat x_j, \ldots, x_{k + 1}}
   \end{split}
  \end{equation*}
 \end{enumerate}
\end{prop}
\begin{proof}
 The verification follows from standard computations and is left for the reader.
\end{proof}

If $\gcon^\TT$ is $TM$-torsion free, by item $(1)$ of Proposition \ref{prop:operatorrelations}, the cohomology of $\pair{\Omega^*_\TT\pair{M}, \dd^\TT}$ is well-defined for degrees $k < 2$. Furthermore, the Jacobiator \eqref{eq:jacobiator} defines a degree $2$ map $\ms J^\TT$ on $\Omega^*_\TT \pair{M}$, which commutes with $\dd^\TT$
\begin{equation}\label{eq:curvaturemap}
 \begin{split}
  & \pair{\ms J^\TT\theta}\pair{x_0, \ldots, x_{k+1}} := \pair{\dd^\TT \circ \dd^\TT \theta}\pair{x_0, \ldots, x_{k+1}}\\
  = & \sum_{i < j < \ell} \pair{-1}^{i+j+\ell+1} \theta\pair{\sair{x_i \diamond_\TT x_j \diamond_\TT x_\ell}, x_0, \ldots, \hat x_i, \ldots, \hat x_j, \ldots, \hat x_\ell, \ldots, x_{k+1}}
 \end{split}
\end{equation}
where $\theta \in \Omega_\TT^k\pair{M}$ and $x_j \in \se{\TT M}$. It follows that $\dd^\TT$ induces a differential $\tilde \dd^\TT$ on the quotient space of the \emph{$\gcon^\TT$-reduced $\TT M$-forms}:
\begin{equation}\label{eq:reducedTTMform}
 0 \to \img \ms J^\TT \to \Omega^*_\TT\pair{M} \xto{\ms R^\TT} \tilde\Omega^*_\TT\pair{M} \to 0
\end{equation}
where $\RRS^\TT$ denotes the quotient map. In particular
\begin{equation}\label{eq:atmost2reduced}
 \tilde\Omega_\TT^k\pair{M} = \Omega_\TT^k\pair{M} \text{ for } k \varleq 2
\end{equation}

\begin{defn}\label{defn:TTcohomology}
 Let $\gcon^\TT$ be a $TM$-torsion free generalized connection on $\TT M$. The complex $\pair{\tilde \Omega_\TT^*\pair{M}, \tilde \dd^\TT}$ is the \emph{$\gcon^\TT$-de Rham complex} and its $k$-th cohomology is the \emph{$k$-th $\gcon^\TT$-de Rham cohomology of $M$}:
 \begin{equation}\label{eq:TTcohomology}
  \tilde H^k_\TT\pair{M} := \dfrac{\ker\pair{\tilde \dd^\TT : \tilde \Omega_\TT^k\pair{M} \to \tilde \Omega_\TT^{k+1}\pair{M}}}{\img \pair{\tilde \dd^\TT : \tilde \Omega_\TT^{k-1}\pair{M} \to \tilde \Omega_\TT^k\pair{M}}}
 \end{equation}
\end{defn}
Because $\ms J^\TT$ is of degree $2$, when $k < 2$, $\tilde H^k_\TT\pair{M}$ compute the corresponding cohomology groups of $\Omega^*_\TT\pair{M}$. By \eqref{eq:functiondiff}, it is evident that $\tilde H^0_\TT\pair{M} = \RR = H^0\pair{M}$, which consists of the constant functions.
The item $(1)$ in Proposition \ref{prop:operatorrelations} then gives a natural inclusion
\begin{equation}\label{eq:firstcohomologyinclusion}
 H^1\pair{M} \subseteq \tilde H^1_\TT\pair{M}
\end{equation}
in which the equality may not hold in general (c.f. Proposition \ref{prop:firstphibcohomology}).

In general, the map $\pi$ induces a natural injection $\pi^*: \Omega^k\pair{M} \into \Omega_\TT^k\pair{M}$ for all $k$:
\begin{equation}\label{eq:pipullback}
 \pair{\pi^*\alpha}\pair{x_1, \ldots, x_k} := \alpha\pair{\pi\pair{x_1}, \ldots, \pi\pair{x_k}}
\end{equation}
where $\alpha \in \Omega^k\pair{M}$ and $x_j \in \se{\TT M}$. Alternatively, $\pi^*$ is induced from the inclusion $T^*M \into \TT M$. When $\gcon^\TT$ is $TM$-torsion free, $\pi^*$ commutes with the derivations
\begin{equation}\label{eq:picommute}
 \pi^*\pair{d\alpha} = \dd^\TT\pair{\pi^*\alpha}
\end{equation}
Thus $\pi^*$ defines a morphism of cochain complices after passing to the quotient
\[\tilde \pi^*: \Omega^k\pair{M} \to \tilde \Omega_\TT^k\pair{M}\]
which induces the corresponding maps on the cohomology groups:
\begin{equation}\label{eq:picohomology}
 \tilde \pi^*: H^k\pair{M} \to \tilde H_\TT^k\pair{M}
\end{equation}

\subsection{$\varphib$-connections}\label{subsec:phibconnections}
To restrict $\gcon^\TT$ further, consider a generalized metric $\GG^b$ \eqref{eq:genmetric}, corresponding to the pair $(g, b)$ of Riemannian metric $g$ on $M$ and $2$-form $b \in \Omega^2(M)$.
\begin{defn}\label{defn:metricconnection}
 A \emph{generalized metric connection with respect to $\GG^b$}, or a \emph{$\GG^b$-metric connection}, is a generalized connection $\gcon^\TT$ on $\TTM$ that preserves both $\GG^b$ and $\<,\>$:
 \begin{equation}\label{eq:gmetricconnection}
  X\<y, z\> = \<\gcon^\TT_x y, z\> + \<y, \gcon^\TT_x z\> \text{ and } X \GG^b\pair{y, z} = \GG^b\pair{\gcon^\TT_x y, z} + \GG^b\pair{y, \gcon^\TT_x z}
 \end{equation}
 where $x, y, z \in \se{\TT M}$ and $X = \pi\pair{x}$.
\end{defn}

A $\GG^b$-metric connection $\gcon^\TT$ preserves the $\pm 1$-eigenbundles $C_\pm^b$ of $\GG^b$:
\begin{equation}\label{eq:geigenbundles}
 C_\pm^b := \bair{x^{b\pm} := X + \pair{b\pm g}X : X \in TM}
\end{equation}
which induces a \emph{$\GG^b$-eigendecomposition} into a quadruple of metric connections $\nabla_\star^\bullet$ on $TM$:
\begin{equation}\label{eq:metricconndecomposition}
 \nabla_{\star, X}^\bullet Y := \pi\pair{\gcon^\TT_{x^{b\star}} y^{b\bullet}}
\end{equation}
where $\star$ and $\bullet$ respectively stand for $+$ or $-$.
Furthermore, the corresponding $\dd^\TT$ admits the induced $\GG^b$-eigendecomposition.
\begin{lemma}\label{lemma:plusminusdecomposition}
 Let $\gcon^\TT$ be a $\GG^b$-metric connection, then the operator $\dd^\TT$ decomposes into components as follows:
 \begin{equation}\label{eq:plusminusdecomposition}
  \dd^\TT = \dd^\TT_+ + \dd^\TT_-: \Omega^{p,q}_\GG \pair{M} \to \Omega^{p+1, q}_\GG \pair{M} \dsum \Omega^{p, q+1}_\GG \pair{M}
 \end{equation}
 where
 \begin{equation}\label{eq:gtypedecomposition}
  \Omega^{p,q}_\GG \pair{M} := \se{\wedge^p C_+^b \tensor \wedge^q C_-^b} \cong \Omega^p\pair{M} \tensor \Omega^q\pair{M}
 \end{equation}
\end{lemma}
\begin{proof}
 Consider $\theta \in \Omega^{p, 0}_\GG\pair{M}$. For $k > 1$, since $\gcon^\TT$ preserves $C^b_\pm$, it is straightforward to verify that for $x_j^{b\pm} \in \se{C^b_\pm}$
 \[\pair{\dd^\TT \theta}\pair{x_0^{b-}, x_1^{b-}, \ldots, x_{k-1}^{b-}, x_{k}^{b+}, \ldots, x_q^{b+}} = 0\]
 Thus $\dd^\TT \theta$ cannot contain any components in $\Omega^{p-k+1, k}_\GG\pair{M}$ for $k > 1$.
 The general situation follows from the analogue for $\theta' \in \Omega^{0,q}_\GG\pair{M}$ then noticing that $\dd^\TT$ is a derivation \eqref{eq:TTMgradedderivation}.
\end{proof}

In its $\GG^b$-eigendecomposition, a $\GG^b$-metric connection $\gcon^\TT$ is $TM$-torsion free iff $\nabla_+^+ = \nabla_-^- = \nabla$ is the Levi-Civita connection for $g$ and for all $X, Y \in \se{TM}$:
\[\nabla_{+, X}^- Y - \nabla_{-, Y}^+ X - [X, Y] = 0\]
It follows that a $3$-form $\phi \in \Omega^3\pair{M}$, which may not be closed, can be defined by
\begin{equation}\label{eq:phidefn}
 \phi\pair{X,Y,Z} := 2g\pair{\pair{\nabla_-^+ - \nabla}_X Y, Z} = 2g\pair{\pair{\nabla_+^- - \nabla}_Y X, Z}
\end{equation}
The mixed components in the $\GG^b$-eigendecomposition can be expressed in terms of $\phi$, e.g.
\begin{equation}\label{eq:mixcomponents}
 \nabla_{-,X}^+ Y = \nabla^{+\phi}_X Y := \nabla_X Y + \dfrac{1}{2} g\inverse \iota_Y\iota_X \phi
\end{equation}
which is a metric connection on $TM$ with totally skew torsion $\phi$. The computations are summarized in the following Theorem / Definition.
\begin{theorem}[$\varphib$-connections]\label{thm:torsionlessmetric}
 The $TM$-torsion free $\GG^b$-metric connections on $\TT M$ are classified by $3$-forms $\phi \in \Omega^3\pair{M}$, which is denoted $\gcon^\phib$, and referred to as the \emph{$\varphib$-connection}.
The $\GG^b$-eigendecomposition of $\gcon^\phib$ is explicitly presented below:
\begin{equation}\label{eq:phibconnectionexplicit}
 \begin{split}
  & \gcon^\phib_{x^{b+}} y^{b+} = \nabla_XY + (b+g)\pair{\nabla_XY} ~~~ \text{ and } ~~~ \gcon^\phib_{x^{b+}} y^{b-} = \nabla^{-\phi}_XY + (b-g)\pair{\nabla^{-\phi}_XY} \\
  & \gcon^\phib_{x^{b-}} y^{b+} = \nabla^{+\phi}_XY + (b+g)\pair{\nabla^{+\phi}_XY} ~~~ \text{ and } ~~~ \gcon^\phib_{x^{b-}} y^{b-} = \nabla_XY + (b-g)\pair{\nabla_XY}
 \end{split}
\end{equation}
where $X, Y \in \se{TM}$, $x^{b\pm} = X + \pair{b\pm g}X \in \se{C^b_\pm}$ and so on. 
In the splitting $\TT M = TM \dsum T^*M$, the $\varphib$-connection takes the form:
\begin{equation}\label{eq:phibconnectionstandard}
 \gcon^\phib_{e^b x} {e^b y} = e^b\sair{\nabla_X y + \dfrac{1}{4} g\inverse \pair{\iota_{g\inverse \eta} \iota_X \phi - \iota_Y \iota_{g\inverse \xi} \phi} + \dfrac{1}{4}\pair{\iota_Y\iota_X \phi - \iota_{g\inverse \eta}\iota_{g\inverse \xi} \phi}}
\end{equation}
where $x = X + \xi, y = Y + \eta \in \se{\TT M}$.
\end{theorem}
\begin{proof}
 Under the pairing $2\aair{,}$, the space $\ms D\pair{\TT M}$ of generalized connections on $\TT M$ is modeled on
 \[\ms D\pair{\TT M} \cong \bair{A : \TT M \to \TT M \tensor \TT M}\]
 The subspace $\ms D_{\GG^b}\pair{\TT M}$ of $\GG^b$-metric connections is then modeled on
 \begin{equation}\label{eq:gbmetricconnspace}
  \ms D_{\GG^b}\pair{\TT M} \cong \bair{A : \TT M \to \pair{\wedge^2 C_+^b} \dsum \pair{\wedge^2 C_-^b} \cong \End\pair{TM, g}^{\dsum 2}}
 \end{equation}
 
 From \eqref{eq:Ttorionfreespace} and \eqref{eq:gbmetricconnspace}, it follows that space of $TM$-torsion free $\GG^b$-metric connections on $\TT M$ is modeled on the intersection. The identity \eqref{eq:phidefn} can be seen also as
 \[\phi\pair{X, Y, Z} = 2\aair{A_{x^{b-}} y^{b+}, z^{b+}} = -2\aair{A_{y^{b+}} x^{b-}, z^{b-}}\]
 which shows that the intersection is isomorphic to $\Omega^3\pair{M}$. The rest follows from straightforward computations.
\end{proof}

Using the explicit formulae in \eqref{eq:phibconnectionexplicit}, it is straightforward to verify that $\gcon^\phib$ is metric compatible with the (almost) Dorfman bracket $*_{\phi-db}$, in the sense of Definition \ref{defn:diamondmetriccompatible} below.
\begin{defn}\label{defn:diamondmetriccompatible}
 A $\GG^b$-metric connection $\gcon^\TT$ is \emph{metric compatible with $*_\gamma$} \eqref{eq:dorfmanbracket} if the diamond bracket $\diamond_\TT$ coincides with $*_\gamma$ on mixed $\GG^b$-eigensections, i.e.
 \begin{equation}\label{eq:diamondmetriccompatible}
  x^{b+} \diamond_\TT y^{b-} = x^{b+} *_{\gamma} y^{b-}
 \end{equation}
 where $x^{b+} \in \se{C^b_+}$ and $y^{b-} \in \se{C^b_-}$.
\end{defn}

The space of $\GG^b$-metric connections that are metric compatible with $*_\gamma$ but not necessarily $TM$-torsion free is modeled on the following subspace of \eqref{eq:gbmetricconnspace}:
\[\ms D_{\GG^b, \gamma}\pair{\TT M} \cong \bair{A : C_+^b \to \wedge^2 C_+^b } \dsum \bair{A : C_-^b \to \wedge^2 C_-^b} \cong \Omega^1\pair{\End\pair{TM, g}}^{\dsum 2}\]

\begin{theorem}\label{thm:levicivita}
 For a closed form $\gamma \in \Omega^3\pair{M}$, let $\phi = \gamma + db$. Then $\gcon^\phib$ is the unique $TM$-torsion free $\GG^b$-metric connection on $\TT M$ that is metric compatible with $*_\gamma$.
 \qed
\end{theorem}

\subsection{$\gcon^\phib$-de Rham cohomology}\label{subsec:cohomology}

The $\GG^b$-eigendecomposition of $\TT M$ induces two left inverses of the natural injection $\pi^*$ in \eqref{eq:pipullback}. There are two obvious projections $p^\pm: \Omega_\TT^k\pair{M} \to \se{\wedge^k C^b_\pm}$ for each $k$:
\[\pair{p^\pm \theta}\pair{x_1^{b\pm}, \ldots, x_k^{b\pm}} := \theta\pair{x_1^{b\pm}, \ldots, x_k^{b\pm}}\]
where $\theta \in \Omega_\TT^k\pair{M}$, $X_j \in \se{TM}$ and $x_j^{b\pm} = X_j + \pair{b\pm g} X_j$ for all $j$. The projection $\pi$ induces the natural isomorphisms $\pi_*: \se{\wedge^* C^b_\pm} \cong \Omega^*\pair{M}$:
\[\pair{\pi_* \theta_\pm}\pair{X_1, \ldots, X_k} := \theta_\pm\pair{x_1^{b\pm}, \ldots, x_k^{b\pm}}\]
where $\theta_\pm \in \se{\wedge^*C^b_\pm}$. Then $\pi^\pm_* = \pi_* \circ p^\pm : \Omega_\TT^k\pair{M} \to \Omega^k\pair{M}$ is given by
\[\pair{\pi^\pm_* \theta}\pair{X_1, \ldots, X_k} := \theta\pair{x_1^{b\pm}, \ldots, x_k^{b\pm}}\]
It is now straightforward to see that
\[\pi^\pm_* \circ \pi^* = \id : \Omega^*\pair{M} \to \Omega^*\pair{M}\]

Let $\dd^\phib$ denote the $\gcon^\phib$-derivation \eqref{eq:TTMderivation} induced by $\gcon^\phib$. 
The decomposition \eqref{eq:plusminusdecomposition} can be explicitly described by the classical de Rham differential and covariant derivatives. For instance, for $\alpha \in \Omega^k\pair{M}$, let $\alpha_+ \in \Omega^{p,0}_\GG\pair{M}$ such that $\alpha = \pi^+_*\pair{\alpha_+}$. Then $\dd^\phib_+ \alpha_+$ is essentially the de Rham differential:
\begin{equation}\label{eq:dphibcomponentplus}
 \begin{split}
  & \pair{\dd^\phib \alpha_+}\pair{x^{b+}_0, x^{b+}_1, \ldots, x^{b+}_k} = \sum_{i}(-1)^i\pair{\gcon^\phib_{x^{b+}_i} \alpha_+}\pair{x^{b+}_0, x^{b+}_1, \ldots, \hat{x^{b+}_i}, \ldots, x^{b+}_k}\\
  = & \sum_{i}(-1)^i\pair{\nabla_{X_i} \alpha}_+\pair{x^{b+}_0, x^{b+}_1, \ldots, \hat{x^{b+}_i}, \ldots, x^{b+}_k} = \pair{d\alpha}_+ \pair{x^{b+}_0, x^{b+}_1, \ldots, x^{b+}_k}
 \end{split}
\end{equation}
while the component $\dd^\phib_- \alpha_+ \in \Omega^{p,1}_\GG\pair{M}$ is essentially given by $\nabla^{+\phi}$:
\begin{equation}\label{eq:dphibcomponent}
 \pair{\dd^\phib_- \alpha_+}\pair{x_0^{b-}, x_1^{b+}, \ldots, x_q^{b+}} 
 = \pair{\nabla^{+\phi}_{X_0} \alpha}_+ \pair{x_1^{b+}, \ldots, x_q^{b+}}
\end{equation}
where $x_j^{b\pm} \in \se{C^b_\pm}$, and $\pi\pair{x_j^{b\pm}} = X_j$.

\begin{lemma}\label{lemma:jacobiatorphib}
 The Jacobiator for $\diamond_\phib$ is given by
 \begin{equation}\label{eq:jacobiatorphib}
  \sair{x^{b\pm} \diamond_\phib y^{b\pm} \diamond_\phib z^{b\mp}} = \pm 2 g\pair{R_{X, Y}^{\mp \phi} Z}
 \end{equation}
 where $x^{b\pm}, y^{b\pm}, z^{b\pm} \in \se{C_\pm^b}$ with $\pi\pair{x^{b\pm}} = X$ and so on, and $R^{\pm\phi}$ are respectively the classical curvature for $\nabla^{\pm \phi}$. All other components in the $\GG^b$-eigendecomposition vanish.
\end{lemma}
\begin{proof}
 Follows by relatively lengthy but standard computations from the definitions.
\end{proof}

\begin{prop}\label{prop:injectivecohomology}
 The natural map $\tilde \pi^*: H^*\pair{M} \to \tilde H_\phib^*\pair{M}$ in \eqref{eq:picohomology} is injective.
\end{prop}
\begin{proof}
 For $\alpha, \beta \in \Omega^*\pair{M}$, suppose that $\pi^*\alpha = \pi^*\pair{d\beta}$, then the injectivity of $\pi^*$ implies that $\alpha = d\beta$. The statement then follows from $\img \ms J^\phib \inter \img \pi^* = \bair{0}$. In fact, suppose $\ms J^\phib \theta = \pi^*\alpha \in \img \ms J^\phib \inter \img \pi^*$, then
 \[\alpha = \pi^+_* \pair{\pi^*\alpha} = \pi^+_* \pair{\ms J^\phib \theta} = 0\]
 The last equality is due to $\sair{x^{b+} \diamond_\phib y^{b+} \diamond_\phib z^{b+}} = 0$, which follows from Lemma \ref{lemma:jacobiatorphib}.
\end{proof}

Since $\tilde \Omega_\TT^k\pair{M} \cong \Omega_\TT^k\pair{M}$ for $k \varleq 2$ and the derivations coincide for $k < 2$, the groups $\tilde H^k_\phib\pair{M}$ for $k < 2$ can be determined. When $k = 0$, it is straightforward to see that $\tilde H^0_\phib\pair{M} \cong \RR$ consists of the constant functions.
\begin{prop}\label{prop:firstphibcohomology}
 Let $P^1_\phi\pair{M}$ be the space of $\nabla$-parallel $1$-forms on $M$, which also annilates $\phi$, i.e.
 \[\xi \in P^1_\phi\pair{M} \iff \nabla \xi = 0 \text{ and } \iota_{g\inverse \xi} \phi = 0\]
 then
 \begin{equation}\label{eq:firstphibcohomology}
  \tilde H^1_\phib\pair{M} \cong H^1\pair{M} \dsum P^1_\phi\pair{M}
 \end{equation}
\end{prop}
\begin{proof}
 Let $\theta \in \Omega^1_\TT\pair{M}$. For $X, Y \in \se{TM}$, let $x^{b\pm} = X + \pair{b\pm g} X$ and so on, define
 \begin{equation}\label{eq:oneformdecomposition}
  \alpha\pair{X} := \theta\pair{x^{b+}} \text{ and } \beta\pair{X} := \theta\pair{x^{b-}}
 \end{equation}
 It then gives the identification 
 \begin{equation}\label{eq:oneformidentificationpair}
  Q:\Omega^1_\TT\pair{M} \cong \Omega^1\pair{M} \dsum \Omega^1\pair{M} : \theta \mapsto  \half \pair{\alpha + \beta, \alpha - \beta}
 \end{equation}
 
 Suppose that $\dd^\phib\theta = 0$, then \eqref{eq:dphibcomponentplus} implies that
 \[d\alpha = d\beta = 0\]
 Hence
 \[0 = \pair{\dd^\phib_-\theta}\pair{x^{b+}, y^{b-}} 
  = -\sair{\nabla^{+\phi}_Y \pair{\alpha - \beta}}\pair{X}\]
 i.e. $\nabla^{+\phi} \pair{\alpha - \beta} = 0$. Set $\xi = \alpha - \beta$, then $d\xi = 0$ implies that
 \[0 = \pair{\nabla^{+\phi}_X \xi} \pair{Y} - \pair{\nabla^{+\phi}_Y \xi} \pair{X} = -\phi\pair{X, Y, g\inverse \xi}\]
 for all $X, Y \in \se{TM}$, from which $\xi \in P^1_\phi\pair{M}$ follows. Thus on $\Omega^1_\TT\pair{M}$
 \[Q\pair{\ker \dd^\phib} = \ker d \dsum P^1_\phi\pair{M}\]
 Then \eqref{eq:firstphibcohomology} follows from $Q\pair{df} = \pair{df, 0}$ for $f \in \se{M}$.
\end{proof}

\begin{corollary}\label{coro:firstphibcohomology}
 $\tilde H^1_\phib\pair{M} \cong H^1\pair{M}$ if one of the following holds
 \begin{enumerate}
  \item
  $\phi$ is non-degenerate, i.e. the following map is injective
  \[\iota_\bullet \phi: \se{TM} \to \Omega^2\pair{M}: X \mapsto \iota_X \phi\]
  \item
  $M$ admits no non-trivial $\nabla$-parallel vector fields.
  \qed
 \end{enumerate}
\end{corollary}

\begin{example}\label{ex:firstphibcohomology}
 Let $M = \RR^n / \ZZ^n$, with the induced flat metric. Suppose that $\phi = 0$, then the $\nabla$-parallel forms on $M$ are the constant forms. Thus
 \[\tilde H^1_{0,b}\pair{M} \cong H^1\pair{M} \dsum T_{0}M \cong \RR^{2n}\]
 where $T_0 M$ is the tangent space at $0 \in M$. On the other hand, for $n = 3$, let $\phi = \vol_g$, the volume form of the flat metric on $M$, then it is non-degenerate. By Corollary \ref{coro:firstphibcohomology}
 \[\tilde H^1_{\vol_g, b}\pair{M} = H^1\pair{M} \cong \RR^3\]
\end{example}

Let $n = \dim_\RR M$. The group $H_\TT^{2n}\pair{M}$ is always well-defined for a $TM$-torsion free $\gcon^\TT$
\[H_\TT^{2n}\pair{M} := \dfrac{\Omega_\TT^{2n}\pair{M}}{\img \pair{\dd^\TT : \Omega_\TT^{2n-1}\pair{M} \to \Omega_\TT^{2n}\pair{M}}}\]
When $\gcon^\TT = \gcon^\phib$, the groups $\tilde H_\phib^{2n}\pair{M}$ and $H_\phib^{2n}\pair{M}$ can be determined.
\begin{prop}\label{prop:topcohomology}
 For any $\varphib \in \Omega^3\pair{M} \dsum \Omega^2\pair{M}$, $\tilde H_\phib^{2n}\pair{M} \cong H_\phib^{2n}\pair{M} \cong \RR$.
\end{prop}
\begin{proof}
 The derivation $\dd^\phib$ on $\Omega_\TT^{2n-1}\pair{M}$ splits in the decomposition \eqref{eq:gtypedecomposition} as
 \[\dd^\phib = \dd^\phib_+ \dsum \dd^\phib_- : \Omega_\TT^{n-1,n}\pair{M} \dsum \Omega_\TT^{n,n-1}\pair{M} \to \Omega_\TT^{n,n}\pair{M} = \Omega_\TT^{2n}\pair{M}\]
 For instance, $\dd^\phib_+ : \Omega_\TT^{n-1,n}\pair{M} \to \Omega_\TT^{n,n}\pair{M}$ is given by the lifting of $d: \Omega^{n-1}\pair{M} \to \Omega^n\pair{M}$ as in \eqref{eq:dphibcomponentplus}. 
 It follows that
 \[\tilde H_\phib^{2n}\pair{M} \cong H_\phib^{2n}\pair{M} \cong \RR\]
 since $\img \JJS^\phib \inter \Omega^{2n}\pair{M} = \{0\}$ by Proposition \ref{prop:topdifferential}, which gives $\tilde \Omega_\TT^{2n}\pair{M} = \Omega_\TT^{2n}\pair{M}$.
\end{proof}

Let $M$ and $M'$ be two smooth manifolds. It is straightforward to see that if a diffeomorphism $\lambda : M \to M'$ relates the corresponding data on both manifolds, i.e.
\[\pair{g, b, \phi} = \pair{\lambda^*g', \lambda^*b', \lambda^*\phi'}\]
it induces isomorphisms throughout the constructions. In particular, it induces the natural isomorphism of cohomology groups $\lambda^*: \tilde H^*_{\phi', b'} \pair{M'} \xto{\cong} \tilde H^*_{\phi, b} \pair{M}$.

\subsection{Laplacians}\label{subsec:laplacians}

Analogous to the classical case, a generalized connection $\gcon$ on $V$ defines the corresponding \emph{$\gcon$-Bochner Laplacian} on $\se V$. Let $\{X_i\}$ be a local orthonormal frame on $TM$ for the Riemannian metric $g$ and let $\{e_i^{b\pm}\}$ be the corresponding local $\GG^b$-orthonormal frame on $\TT M$:
\[e_i^{b\pm} := e^b\sair{X_i \pm g\pair{X_i}}\]
\begin{defn}\label{defn:bochnerlaplacian}
 Let $\gcon$ be a generalized connection on $V$. The \emph{$\gcon$-Bochner Laplacian (with respect to $\gcon^\phib$)} on $V$ is given for $v \in \se{V}$ by:
 \begin{equation}\label{eq:bochnerlaplacian}
  \laplacian_\gcon v := - \sum_i\sair{ \pair{\gcon_{e_i^{b+}}\gcon_{e_i^{b+}} - \gcon_{\gcon^{\phib}_{e_i^{b+}} e_i^{b+}}} v +  \pair{\gcon_{e_i^{b-}}\gcon_{e_i^{b-}} - \gcon_{\gcon^{\phib}_{e_i^{b-}} e_i^{b-}}} v }
 \end{equation}
\end{defn}
It is straightforward to verify that \eqref{eq:bochnerlaplacian} is independent of the choice of $\bair{X_i}$. Because $\gcon^{\phib}_{e_i^{b\pm}} e_i^{b\pm}$ involve only the Levi-Civita connection for $g$, it is evident that
\[\laplacian_\gcon v = \pair{\laplacian^b_+ + \laplacian^b_-} v\]
where $\laplacian^b_{\pm}$ are the classical Bochner Laplacians for $\nabla^b_\pm$ respectively. Thus $\laplacian_\gcon$ is a second order elliptic operator, which in general depends on both $g$ and $b$, i.e. on $\GG^b$, but not on $\phi$. When $\gcon$ is a lift of a classical connection, $\laplacian_\gcon$ reduces to (a constant multiple of) the corresponding classical Bochner Laplacian.

The operator $\dd^\phib$ can alternatively be written as
\begin{equation}\label{eq:TTMderivationframe}
 \dd^\phib \theta = \sum_{j}\pair{e_j^{b+} \wedge \gcon^\phib_{e_j^{b+}} \theta - e_j^{b-} \wedge\gcon^\phib_{e_j^{b-}} \theta}
\end{equation}
where $\theta \in \Omega^*_\TT\pair{M}$. Thus the symbol of $\dd^\phib$ is
\[\sigma\pair{\dd^\phib} = 2 \sqrt{-1} \xi \wedge : \wedge^k \TT M \to \wedge^{k+1} \TT M\]
where $\xi \in T^*M$. It follows that whenever $\dd^\phib$ squares to $0$ it defines an elliptic complex. Analogous to the classical case, for $\theta\in \Omega^*_\TT\pair{M}$, define
\begin{equation}\label{eq:phibdformaldual}
 \dd^{\phib*}\theta := - \half \sum_j \pair{\iota_{e_j^{b+}}\gcon^\phib_{e_j^{b+}} \theta + \iota_{e_j^{b-}}\gcon^\phib_{e_j^{b-}} \theta }
\end{equation}
where the $\half$ is due to the convention \eqref{eq:contractionconvention}. The symbol of $\dd^{\phib*}$ is then
\[\sigma\pair{\dd^{\phib*}} = -\sqrt{-1}\iota_{g\inverse\xi + bg\inverse \xi} : \wedge^k \TT M \to \wedge^{k-1} \TT M\]
where $\xi \in T^*M$.
\begin{defn}\label{defn:hodgelaplacian}
 The \emph{$\gcon^\phib$-Hodge Laplacian} is the operator on $\Omega^*_\TT\pair{M}$ given by
 \begin{equation}\label{eq:hodgelaplacian}
  \laplacian^\phib := \dd^\phib \dd^{\phib*} + \dd^{\phib *} \dd^\phib
 \end{equation}
\end{defn}

\begin{prop}\label{prop:hodgelaplacian}
 The $\gcon^\phib$-Hodge Laplacian $\laplacian^\phib$ is a second order elliptic operator.
\end{prop}
\begin{proof}
 It is clear from the discussions above that the symbol of $\laplacian^\phib$ is
 \[\sigma\pair{\laplacian^\phib} = 2\norm{\xi}^2_g: \wedge^k \TT M \to \wedge^{k} \TT M\]
 where $\xi \in T^*M$, from which the statement follows. 
\end{proof}

\begin{theorem}\label{thm:phibhodgeforh1}
 Let $M$ be a closed manifold. Then the following holds:
 \begin{equation}\label{eq:phibhodgeforh1}
  \tilde H^1_\phib\pair{M} \cong \left.\ker \dd^\phib\right|_{\Omega^1_\TT\pair{M}} \inter \left.\ker \dd^{\phib*}\right|_{\Omega^1_\TT\pair{M}} \subseteq \left.\ker \laplacian^\phib\right|_{\Omega^1_\TT\pair{M}}
 \end{equation}
\end{theorem}
\begin{proof}
 Let $\theta \in \Omega^1_\TT\pair{M}$ and consider $\alpha, \beta \in \Omega^1\pair{M}$ as in \eqref{eq:oneformdecomposition}, then it can be shown that
 \[\dd^{\phib*}\theta = d^*\pair{\alpha + \beta}\]
 Combining the identity above with the proof of Proposition \ref{prop:firstphibcohomology}, the first isomorphism follows from the classical Hodge theorem. The last inclusion is obvious.
\end{proof}

Even though $\pair{\Omega_\TT^*\pair{M}, \dd^\phib}$ is generally not a chain complex, Theorem \ref{thm:phibhodgeforh1} nonetheless hints at the analogue of its ``cohomology groups''.
\begin{defn}\label{defn:phibcohomology}
 Let $\GG^b$ be a generalized metric on $M$. The \emph{$\varphib$-pseudo-cohomology groups $\check H^*_\phib\pair{M}$} of $M$ consist of the intersection of the kernels of $\dd^\phib$ and $\dd^{\phib*}$:
 \[\check H^k_\phib\pair{M} := \left.\ker \dd^\phib\right|_{\Omega^k_\TT\pair{M}} \inter \left.\ker \dd^{\phib*}\right|_{\Omega^k_\TT\pair{M}}\]
 The \emph{$\varphib$-Laplacian kernels $\hat H^k_\phib\pair{M}$} are the subspaces
 \[\hat H^k_\phib\pair{M} := \left.\ker \laplacian^\phib\right|_{\Omega^k_\TT\pair{M}}\]
\end{defn}

\begin{corollary}\label{coro:cohomologyfinitedim}
  For a closed manifold $M$, both $\hat H^*_\phib\pair{M}$ and $\check H^*_\phib\pair{M}$ are finite dimensional.
\end{corollary}
\begin{proof}
 It follows from $\check H^k_\phib\pair{M} \subseteq \hat H^k_\phib\pair{M}$ and the ellipticity of $\laplacian^\phib$.
\end{proof}

\begin{remark}\label{remark:inequalitytokerlaplacian}
 The inclusions $\check H^k_\phib\pair{M} \subseteq \hat H^k_\phib\pair{M}$ may not be equalities. For $k = 1$, direct computation using the $\GG^b$-eigendecomposition shows that
 \[\theta = \hat H^1_\phib \pair{M} \iff \begin{cases}\laplacian\pair{\alpha + \beta} + \laplacian_{+ \phi}\pair{\alpha - \beta} {- \half \sum_{j,k} \pair{d\beta}\pair{X_j, X_k} \iota_{X_k} \iota_{X_j} \phi}= 0 \\ \\ \laplacian\pair{\alpha + \beta} - \laplacian_{- \phi}\pair{\alpha - \beta} {+ \half \sum_{j,k}\pair{d\alpha}\pair{X_j, X_k} \iota_{X_k} \iota_{X_j}\phi}= 0\end{cases}\]
 where $\alpha, \beta$ are as in \eqref{eq:oneformdecomposition} and $\{X_j\}$ is a local $g$-orthonormal frame of $TM$. In the case that $d\alpha = d\beta = 0$, it can be shown that the right hand side is exactly equivalent to $\theta \in \check H^1_\phib\pair{M}$. Namely, under the identification \eqref{eq:oneformidentificationpair}
 \[\check H^1_\phib\pair{M} = \hat H^1_\phib\pair{M} \inter \pair{\ker d \dsum \ker d}\]
\end{remark}

\subsection{Compact Lie groups}\label{subsec:compactLiegroups}

Manifolds admitting flat metric connections with non-trivial completely skew torsions are known by Cartan and Schouten (\cite{CartanSchouten:26a}, \cite{CartanSchouten:26b}), which are essentially compact Lie groups and $S^7$ (see also Agricola-Friedrich \cite{AgricolaFriedrich:FlatMetricConnection:09}).

Suppose that $G$ is a real semi-simple Lie group, endowed with the bi-invariant Killing metric $g$ and the corresponding bi-invariant Cartan $3$-form $\gamma \in \Omega^3\pair{G}$. In this case, as will become clear below, the computations are very much parallel to those for the classical de Rham cohomology of the doubled group $G \times G$.

The metric connections $\nabla^{\pm \gamma}$, with torsions $\pm \gamma$ respectively, are flat. Let $\pair{\phib} = \pair{\gamma, 0}$ and consider the $\pair{\gamma, 0}$-connection $\gcon^{\gamma,0}$.
By Lemma \ref{lemma:jacobiatorphib}, for all $x, y, z \in \se{\TT G}$
\[\sair{x \diamond_{\gamma,0} y \diamond_{\gamma,0} z} = 0\]
which implies that $\JJS^{\gamma,0} = \dd^{\gamma,0} \circ \dd^{\gamma,0} = 0$. Hence $\pair{\Omega_\TT^*\pair{G}, \dd^{\gamma,0}}$ is the $\gcon^{\gamma,0}$-de Rham complex, whose cohomology is the $\gcon^{\gamma,0}$-de Rham cohomology $\tilde H_{\gamma,0}^*\pair{G}$.

Let $X_u^l$ denote the left invariant vector field on $G$ such that $X_u^l\pair{e} = u \in \mf g := T_e G$. The corresponding right invariant vector fields are denoted $X_u^r$. The Lie algebra structure on $\mf g$ is identified with the Lie algebra of left invariant vector fields:
\[\sair{X_u^l, X_v^l} = X_{\sair{u, v}}^l\]
For $u, v \in \mf g$, set $\theta_u^r := g\pair{X_u^r}$ and $\theta_v^l := g\pair{X_v^l}$, then
\[x_u^+ = X_u^r + \theta_u^r \in \se{C_+} \text{ and } x_v^- = X_v^l - \theta_v^l \in \se{C_-}\]
It is straightforward to see that
\[\gcon^{\gamma,0}_{x_u^+} x_v^+ = -\half x_{\sair{u, v}}^+, \gcon^{\gamma,0}_{x_u^-} x_v^- = \half x_{\sair{u, v}}^-, \text{ and } \gcon^{\gamma,0}_{x_u^+} x_v^- = \gcon^{\gamma,0}_{x_u^-} x_v^+ = 0\]

Let $\uu = \pair{u, u'}, \vv = \pair{v, v'} \in \mf g \dsum \mf g$ and
\[x_{\uu} = -x_u^+ + x_{u'}^-, x_\vv = -x_v^+ + x_{v'}^- \in \se{\TT G}\]
then direct computation leads to
\begin{equation}\label{eq:diamondtrivialization}
 \gcon^{\gamma,0}_{x_\uu} x_\vv = \half x_{\sair{\uu, \vv}} \Longrightarrow 
 x_\uu \diamond_{\gamma,0} x_\vv = x_\uu *_{\gamma} x_\vv = x_{\sair{\uu, \vv}}
\end{equation}
where the Lie bracket on $\mf g \dsum \mf g$ is the direct sum of those on each factor. The last equality in \eqref{eq:diamondtrivialization} can be seen also from the \emph{Courant trivialization} of Alekseev-Bursztyn-Meinrenken \cite{AlekseevBursztynMeinrenken:PureSpinor:09}. 
The $\pair{\gamma,0}$-curvature can then be computed as
\[\grie^{\gamma,0}_{x_\uu, x_\vv} x_\ww = -\dfrac{1}{4}\sair{\sair{x_\uu, x_\vv}, x_\ww} = -\dfrac{1}{4} x_{\sair{\sair{\uu, \vv}, \ww}}\]
which gives the $\pair{\gamma,0}$-Ricci tensor
\[\Rc^{\gamma, 0}\pair{x_\uu, x_\vv} = \dfrac{1}{4} \GG\pair{x_\uu, x_\vv}\]
In particular, $\pair{G, \gamma; \GG}$ may be seen as an example of a \emph{$\pair{\gamma, 0}$-Einstein manifold}, where the $\pair{\gamma,0}$-Ricci curvature is proportional to the generalized metric.

To compute $\tilde H_{\gamma, 0}\pair{G}$, recall that $\TT G$ is the dual of itself via $2\aair{, }$, which leads to
\[\pair{\dd^{\gamma,0} x_\uu}\pair{x_\vv, x_\ww} = - x_\uu \pair{x_{\sair{\vv, \ww}}}\]
Let $\theta \in \Omega_\TT^*\pair{G}$ be decomposable as the product of $k$ sections of the form $x_{\vv}$, then
\[\pair{\dd^{\gamma,0}\theta}\pair{x_{\uu_0}, \ldots, x_{\uu_k}} = \sum_{i<j} \pair{-1}^{i+j} \theta\pair{x_{\sair{\uu_i, \uu_j}}, x_{\uu_0}, \ldots, \hat x_{\uu_i}, \ldots, \hat x_{\uu_j}, \ldots, x_{\uu_k}}\]
Let $f_\uu \in \mf g^* \dsum \mf g^*$ be defined by
\[f_\uu \pair{\vv} := x_\uu\pair{x_\vv}\]
It induces an inclusion of the Chevalley-Eilenberg complex of $\mf g\dsum \mf g$ for the trivial module
\[\pair{\wedge^*\pair{\mf g^* \dsum \mf g^*}, \delta} \into \pair{\Omega^*_\TT\pair{G}, \dd^{\gamma,0}}: f_\uu \mapsto x_\uu\]
Similar to the classical case, this induces an isomorphism on the cohomology:
\[\tilde H_{\gamma, 0}^*\pair{G} \cong H^*\pair{\mf g \dsum \mf g} \cong H^*\pair{G \times G}\]
The isomorphism above in fact is an isomorphism of rings, where on $\tilde H^*_{\gamma, 0}\pair{G}$ the product is induced by the wedge product in $\Omega_\TT^*\pair{M}$.

\section{Curvature tensors}\label{sec:phibcurvature}

Consider a generalized connection $\gcon$ on $V$. Let $\gcon^\TT$ be any generalized connection on $\TT M$. The $\gcon^\TT$-derivation $\dd^\TT$ \eqref{eq:TTMderivation} extends to $\Omega_\TT^k\pair{V} := \se{\wedge^k \TTM \tensor V}$:
\begin{equation}\label{eq:derivationextension}
 \dd^\TT_\gcon(\theta \tensor v) := \pair{\dd^\TT\theta} \tensor v + (-1)^k \theta \wedge \gcon v
\end{equation}
where $v \in \se V$. The \emph{$\gcon^\TT$-curvature operator} $\gcurv^{\TT}\pair{\gcon}$ of $\gcon$ is then given by
\begin{equation}\label{eq:gencurvaturedefn}
 \gcurv^{\TT}\pair{\gcon} := \dd^\TT_\gcon \circ \gcon
\end{equation}
which generally is not tensorial in $v$ if $\gcon^\TT$ is not $TM$-torsion free. When $\gcon$ is understood, it will be dropped from the notations $\gcurv^\TT \pair{\gcon}$ and $\dd^\TT_\gcon$.

In terms of covariant derivatives, the $\gcon^\TT$-curvature operator $\gcurv^{\TT}$ is given by \eqref{eq:gencurvatureintro}. It is tensorial iff $\gcon^\TT$ is $TM$-torsion free, in which case, for any $f \in \se{M}$
\begin{equation}\label{eq:curvtensorialproof}
 \gcurv^{\TT}_{x, y} (fv) - f\gcurv^{\TT}_{x, y} v = \pair{\sair{\pi(x), \pi(y)} - \pi\pair{x \diamond_\TT y}} (f) v = 0
\end{equation}
The resulting tensor $\gcurv^{\TT} \in \Omega_\TT^2\pair{\End\pair{V}}$ is the \emph{$\gcon^\TT$-curvature} of $\gcon$.  Similar to \eqref{eq:derivationextension}, let $\tilde \dd^\TT_\gcon$ be the extension of $\tilde \dd^\TT$ to $\tilde \Omega^*_\TT\pair{V} := \tilde \Omega_\TT^*\pair{M} \tensor \se{V}$. By \eqref{eq:atmost2reduced}, $\gcurv^\TT \pair{\gcon}$ can be seen as an element in $\tilde \Omega_\TT^2 \pair{\End\pair{V}}$ and \eqref{eq:gencurvaturedefn} can also be rewritten as
\begin{equation}\label{eq:reducedcurvature}
 \gcurv^\TT \pair{\gcon} = \tilde \dd^\TT_\gcon \circ \gcon
\end{equation}

\begin{example}\label{ex:hermitianline}
 Let $\pair{V, h}$ be a Hermitian vector bundle on $M$. A generalized connection $\gcon$ on $V$ \emph{preserves $h$}, or is \emph{($h$-)unitary}, if for $v_j  \in \se V$ and $x \in \TT M$ with $X = \pi(x)$:
 \[X h(v_1, v_2) = h\pair{\gcon_x v_1, v_2} + h\pair{v_1, \gcon_x v_2}\]
 Suppose now that $V$ is a Hermitian line bundle and $s$ is a local section of $V$ such that $h(s, s) = 1$. Since $\gcon$ is unitary, it is determined by a local section $u \in \se \TTM$, such that for $x \in \TTM$
 \[\sqrt{-1}\gcon_x s = 2\<x, u\> s\]
 Analogous to the classical computation, the $\gcon^\TT$-curvature for the line bundle $V$ is then
 \begin{equation}\label{eq:linbundlecurvature}
  \sqrt{-1}\gcurv^{\TT}_{x,y}(\gcon) = 2\pair{\<y, \gcon^{\TT}_x u\> - \<x, \gcon^{\TT}_y u\>} = \pair{\dd^{\TT}u}\pair{x, y}
 \end{equation}
\end{example}

\subsection{Chern-Weil homomorphism}\label{subsec:chernweil}

Let $\gcon^\TT$ be a $TM$-torsion free generalized connection on $\TT M$. Since $\dd^\TT$ generally does not square to $0$ \eqref{eq:dsquared}, the Bianchi identity generally does not hold for $\gcurv^\TT\pair{\gcon}$. In terms of covariant derivatives, $\dd^\TT_\gcon \gcurv^\TT\pair{\gcon}$ can be expanded into
\begin{equation}\label{eq:differentialbianchi}
 \begin{split}
  & \pair{\dd^\TT_\gcon \gcurv^\TT}_{x, y, z} = \gcon_x \gcurv^\TT_{y,z} - \gcurv^\TT_{\gcon^\TT_x y, z} - \gcurv^\TT_{y, \gcon^\TT_x z} - \gcurv^\TT_{y,z}\gcon_x + c.p. \text{ in } x, y, z\\
  = & - \gcon_{\grie^\TT_{x, y} z} - c.p. \text{ in } x, y, z 
 \end{split}
\end{equation}
where $x, y, z \in \se{\TT M}$. By \eqref{eq:jacobiatorcurvature} and \eqref{eq:projectedjacobi}, it gives
\begin{equation}\label{eq:differentialbianchitensor}
 \pair{\dd^\TT_\gcon \gcurv^\TT}_{x, y, z} = \gcon_{\sair{x \diamond_\TT y \diamond_\TT z}} = \psi_{\sair{x \diamond_\TT y \diamond_\TT z}}
\end{equation}
which leads to the Bianchi identiy over $\tilde \Omega_\TT^*\pair{M}$.
\begin{lemma}\label{lemma:bianchireduced} 
 Let 
 $\gcurv^\TT \in \Omega_\TT^2\pair{M} \tensor \End\pair{V}$ be the $\gcon^\TT$-curvature of a generalized connection $\gcon$ on $V$, then
 \begin{equation}\label{eq:reducedbianchi}
  \tilde \dd^\TT_\gcon \gcurv^\TT = 0
 \end{equation}
\end{lemma}
\begin{proof}
 It follows from \eqref{eq:differentialbianchitensor} that
\[\dd^\TT_\gcon \gcurv^\TT \in \img \ms J^\TT \tensor \End\pair{V}\]
Thus $\tilde \dd^\TT_\gcon \gcurv^\TT = \ms R^\TT\pair{\dd^\TT_\gcon \gcurv^\TT} = 0$.
\end{proof}

The space $\DDS(V)$ of generalized connections on $V$ is an affine space modeled on $\Omega_\TT^1\pair{\End\pair{V}}$, which coincides with $\tilde \Omega_\TT^1\pair{\End\pair{V}}$. For $A \in \tilde \Omega_\TT^1\pair{ \End\pair{V}}$, standard computation gives
\[\gcurv^{\TT}\pair{\gcon + A} - \gcurv^{\TT}\pair{\gcon} = \dd^\TT_\gcon A + A\wedge A\]
It follows that
\begin{equation}\label{eq:chernindependentofconnection}
 \tr_V\pair{\gcurv^{\TT}\pair{\gcon + A}} - \tr_V\pair{\gcurv^{\TT}\pair{\gcon}} = \tr_V\pair{\tilde \dd^\TT_\gcon A} = \tilde \dd^\TT \tr_V\pair{A}
\end{equation}
As in the classical case, \eqref{eq:reducedbianchi} implies that
\[\tilde \dd^\TT \tr_V \gcurv^\TT\pair{\gcon} = \tr_V \tilde \dd^\TT_\gcon \gcurv^\TT\pair{\gcon} = 0\]

The gauge group $\Aut\pair{V}$ acts on $\DDS(V)$ by push-forward
\[\pair{\lambda \gcon}_x v := \lambda\inverse \sair{\gcon_x\pair{\lambda v}}\]
where $\lambda \in \Aut\pair{V}$, $x \in \TT M$ and $v \in \se{V}$. It induces the action on the curvature by conjugation
\[\gcurv^\TT\pair{\lambda\gcon} = \lambda\inverse \gcurv^\TT\pair{\gcon} \lambda\]
The Chern-Weil homomorphism extends to define characteristic classes for $V$ in $\tilde H^*_\TT\pair{M}$.
\begin{defn}\label{defn:chernclasses}
 For a Hermitian vector bundle $\pair{V, h}$ over $M$, its \emph{$k$-th $\gcon^\TT$-Chern class } is
 \begin{equation}\label{eq:chernclasses}
  c^\TT_k\pair{V} := \sair{\tr_V\pair{\sqrt{-1}\gcurv^\TT\pair{\gcon}}^k} \in \tilde H^{2k}_\TT \pair{M}
 \end{equation}
 where $\gcon$ is a generalized connection on $V$. For real vector bundles, their \emph{$\gcon^\TT$-Euler} and \emph{$\gcon^\TT$-Pontrjagin} classes can similarly be defined, as elements of $\tilde H^*_\TT\pair{M}$ of appropriate degrees.
\end{defn}

By \eqref{eq:chernindependentofconnection}, the $\gcon^\TT$-Chern classes do not depend on the choice of $\gcon$ on $V$. Let $\gcon$ be the lift of a classical connection $\nabla_0$ on $V$ and let $F_0$ be the classical curvature of $\nabla_0$, then
\[\gcurv^\TT_{x, y} = F_{0; \pi\pair{x}, \pi\pair{y}}\]
for $x, y \in \se{\TT M}$. This relates $c^\TT_*\pair{V}$ to the classical Chern classes $c_*\pair{V}$.
\begin{prop}\label{prop:phibchernclass}
 Let $\pair{V, h}$ be a Hermitian vector bundle, then 
 \[c^\TT_k \pair{V} = \tilde \pi^* c_k \pair{V}\]
 for all $k$. In particular, $c^\TT_k \pair{V} = 0$ for all $k$ such that $2k > n$. Similarly, the $\gcon^\TT$-Euler and Pontrjagin classes are the images of the respective classical classes under $\tilde \pi^*$.
 \qed
\end{prop}

\subsection{$\varphib$-curvatures}\label{subsec:phibcurvatures}

A generalized Riemannian metric $\GG^b$ induces an eigendecomposition of a generalized connection $\gcon$ on $V$. For $X \in \se{TM}$ and $v \in \se{V}$, let $x^{b\pm} = X + \pair{b\pm g} X$, then:
\begin{equation}\label{eq:gconeigendecomposition}
 \nabla^b_{\pm, X} v := \gcon_{x^{b\pm}} v
\end{equation}
The connections $\nabla^b_\pm$ depend on $b$, while their difference does not (c.f. \cite{Gualtieri:Poisson:10}):
\begin{equation}\label{eq:gconeigendiff}
 \cdiff_X := \dfrac{1}{2}\pair{\nabla^b_{+, X} - \nabla^b_{-, X}} = \dfrac{1}{2}\pair{\gcon_{x^{b+}} - \gcon_{x^{b-}}} = \gcon_{g(X)}
\end{equation}
The average of $\nabla^b_\pm$ gives the \emph{$b$-neutral connection} of $\gcon$:
\begin{equation}\label{eq:neutralconnection}
 \nabla^b_{0, X} := \half \pair{\nabla^b_{+,X} + \nabla^b_{-, X}} = \gcon_{X+\iota_X b}
\end{equation}
which leads to 
\begin{equation}\label{eq:gconneutraldecomp}
 \nabla^b_\pm = \nabla^b_0 \pm \psi
\end{equation}

When $\cdiff = 0$, the generalized connection $\gcon$ is the lift of a classical connection $\nabla_0$ on $V$, in which case $\nabla^b_\pm = \nabla_0$ are independent of $b$ as well. The dependence on $b$ of the $\GG^b$-eigendecomposition of $\gcon$ can be described in terms of $\cdiff$ as follows.
\begin{prop}\label{prop:eigendecompositionbdependence}
 For $b \in \Omega^2(M)$, let 
 \[j_b := g\inverse b: TM \to TM: X \mapsto g\inverse\pair{\iota_X b}\]
 and $\pair{\nabla^b_+, \nabla^b_-}$ be the $\GG^b$-eigendecomposition of $\gcon$ on $V$. For $a\in \Omega^2(M)$, then
 \[\nabla^{b+a}_\pm = \nabla^b_\pm + \cdiff_{j_a}\]
 gives the $\GG^{b+a}$-eigendecomposition of $\gcon$. \qed
\end{prop}

\begin{defn}\label{defn:phibcurvature}
 Let $\gcon$ be a generalized connection on a vector bundle $V$ over $M$. The \emph{$\varphib$-curvature} $\gcurv^\phib\pair{\gcon}$ of $\gcon$ is its $\gcon^\phib$-curvature \eqref{eq:gencurvaturedefn}. When $\gcon$ is understood it is also simply denoted $\gcurv^\phib$. 
\end{defn}

Given the pair of classical connections $\pair{\nabla^b_+, \nabla^b_-}$, besides the curvature $F_\pm^b$ of each of them, there is also a \emph{mixed curvature} $F_{+, -}^b$ (\cite{Wang:ToricGenKahlerIII:19}):
\begin{equation}\label{eq:mixedcurvature}
 F^\phib_{+, -; X, Y} v := \pair{\nabla^b_{+,X}\nabla^b_{-,Y} - \nabla^b_{-, \nabla^{-\phi}_X Y}} v - \pair{\nabla^b_{-,Y}\nabla^b_{+,X} - \nabla^b_{+, \nabla^{+\phi}_Y X}} v
\end{equation}
where $X, Y \in \se{TM}$ and $v \in \se{V}$. It can also be expressed using the tensor $\cdiff$
\begin{equation}\label{eq:mixedcurvdecompose}
 F^\phib_{+, -; X, Y} = F^b_{+, X, Y} - 2\pair{\nabla_{+,X}^b \cdiff}_Y - \pair{\iota_{g\inverse \cdiff} \phi}\pair{X, Y}
\end{equation}
Let $F^b_0$ be the classical curvature for $\nabla^b_0$ \eqref{eq:neutralconnection}, it gives another decomposition for the mixed curvature \eqref{eq:mixedcurvature}:
\begin{equation}\label{eq:mixedcurvdecomposealternative}
 F_{+,-;X,Y}^\phib = F^b_{0;X,Y} + \pair{\iota_{g\inverse\cdiff}\phi}\pair{X, Y} - \sair{\cdiff_X, \cdiff_Y} - \sair{\pair{\nabla^b_{0,X}\cdiff}_Y + \pair{\nabla^b_{0,Y}\cdiff}_X}
\end{equation}

\begin{theorem}\label{thm:phibcurvaturedecompose}
 The $\varphib$-curvature admits $\GG^b$-eigendecomposition in terms of the (mixed) curvatures of the pair $\pair{\nabla^b_+, \nabla^b_-}$ of classical connections as follows:
 \begin{equation}\label{eq:phibcurvaturedecompose}
  \gcurv^\phib_{x^{b\pm}, y^{b\pm}} v = F^b_{\pm, X, Y} v ~~~ \text{ and } ~~~ \gcurv^\phib_{x^{b+}, y^{b-}} v = F^\phib_{+, -; X, Y} v
 \end{equation}
\end{theorem}
\begin{proof}
 Straightforward from the definition, via the $\GG^b$-eigendecomposition.
\end{proof}

\begin{example}\label{ex:hermitianlinephib}
 Continue with Example \ref{ex:hermitianline} for $\gcon^\TT = \gcon^\phib$. In this case, the local section $u \in \se{\TT M}$ decomposes into
 \[u = \half  \sair{\pair{g\inverse \nu_+ + bg\inverse \nu_+ + \nu_+} - \pair{g\inverse \nu_-+ bg\inverse \nu_- - \nu_-}}\]
 where $\sqrt{-1}\nu_\pm \in \Omega^1\pair{M}$ are the local $1$-forms defining the connections $\nabla^b_\pm$ respectively. It follows that
 \[\sqrt{-1}F^b_\pm = d\nu_\pm\]
 and the mixed component in $\gcurv^\phib$ is given by
 \[\sqrt{-1}F^b_{+,-; X, Y} = \sair{\nabla^{-\phi}_X \nu_-}\pair{Y} - \sair{\nabla^{+\phi}_Y \nu_+}\pair{X}\]
 Note that $F^b_{+,-}$ is neither symmetric nor skew-symmetric in $X$ and $Y$, and decomposes into symmetric and skew-symmetric parts as
 \[\sqrt{-1}F^b_{+,-} = - \LLC_{g\inverse\psi} g + \pair{\sqrt{-1}F^b_0 - \iota_{g\inverse \psi} \phi}\]
 where $F^b_0$ is the curvature of the $b$-neutral connection $\nabla^b_0$ and $\psi = \half \pair{\nu_+ - \nu_-}$, in the decomposition \eqref{eq:gconneutraldecomp} of $\gcon$.
\end{example}

\begin{corollary}\label{coro:phibcurvaturestandard}
 In the decomposition $\TT M = e^b TM \dsum T^*M$ the $\varphib$-curvature is
 \begin{equation}\label{eq:phibcurvaturestandard}
  \begin{split}
   \gcurv^\phib_{x, y} = & F^b_{0; X, Y} + \pair{\nabla^b_{0, X} \cdiff}_{g\inverse \eta} - \pair{\nabla^b_{0, Y} \cdiff}_{g\inverse \xi} + \sair{\cdiff_{g\inverse \xi}, \cdiff_{g\inverse \eta}} \\
   & + \half \sair{\pair{\iota_{g\inverse \cdiff} \phi} \pair{X, Y} -  \pair{\iota_{g\inverse \cdiff} \phi} \pair{g\inverse \xi, g\inverse \eta}}
  \end{split}
 \end{equation}
 where $x = X + \iota_X b + \xi$ and $y = Y + \iota_Y b + \eta$.
 \qed
\end{corollary}

\begin{remark}\label{remark:dependence}
 The dependence of $\gcurv^\phib$ on $\phi$ is completely contained in the last term 
 of \eqref{eq:mixedcurvdecompose}, e.g. for $\varphi \in \Omega^3(M)$:
 \[F^{\phi+\varphi, b}_{+,-; X, Y} - F^{\phi, b}_{+,-; X, Y} = - \pair{\iota_{g\inverse \varphi}\phi} \pair{X, Y}
 \]
 The dependence of $\gcurv^\phib$ on $b$ is more complicated, which can be derived from \eqref{eq:phibcurvaturestandard} by relatively lengthy computations, noting that $\xi$ and $\eta$ in the expression depend on $b$ as well.
\end{remark}

\subsection{Yang-Mills functional}\label{subsec:yangmills}
It is straightforward to extend the Yang-Mills functional to this context. 
The generalized metric $\GG^b$ induces natural inner product on $\wedge^*\TT M$. For a Hermitian bundle $\pair{V, h}$, it induces a natural norm on $\wedge^*\TT M \tensor \End\pair{V}$, denoted $\norm{\bullet}_{\GG^b, h}$. The \emph{$\gcon^\TT$-Yang-Mills functional} on $\DDS(V)$ is given by
\begin{equation}\label{eq:yangmillsfunctional}
 \YM_\TT\pair{\gcon} := \int_M \norm{\gcurv^\TT\pair{\gcon}}^2_{\GG^b, h} d\vol_g
\end{equation}
It is evidently invariant under the gauge action on $\DDS(V)$. When restricted to the subspace of the lifts of classical connections on $V$, $\YM_\TT\pair{\gcon}$ reduces to (a constant multiple of) the classical Yang-Mills functional. It can also be regarded as a functional of the pair $\pair{\gcon, \gcon^\TT}$ of generalized connections on $V$ and $\TT M$ respectively.

When $\gcon^\TT = \gcon^\phib$, it can be represented as
\begin{equation}\label{eq:genYangMills}
 \YM_\phib\pair{\gcon} = YM\pair{\nabla^b_+} + YM\pair{\nabla^b_-} + 2\int_M \norm{F^\phib_{+,-}}^2_{g,h} d\vol_g
\end{equation}
where $YM\pair{\bullet}$ denotes the classical Yang-Mills functional. The $2$-form $b$ affects only the $\GG^b$-eigendecomposition of $\gcon$. The right hand side can be seen as a functional for a pair of classical connections $\pair{\nabla_+, \nabla_-}$, where the last term encodes the dependence on $\phi \in \Omega^3\pair{M}$ (Remark \ref{remark:dependence}), as well as the interaction within the pair.

\section{Curvatures on $\TT M$} \label{sec:phibriemannian}

For a $\GG^b$-metric connection $\gcon^\TT$ on $\TT M$, its $\varphib$-curvature is denoted $\grie^{\TT, \phib}$ and the associated \emph{curvature tensor} is given by
\begin{equation}\label{eq:genRiemann}
 \grie^{\TT, \phib}(x, y, z, w) := \GG^b\pair{\grie^{\TT, \phib}_{x,y} z, w}
\end{equation}
where $x, y, z, w \in \se \TTM$. Similar to the classical situation, it is skew in the first two and the last two entries respectively:
\begin{equation*}
 \grie^{\TT, \phi, b}(x, y, z, w) = - \grie^{\TT, \phi, b}(y, x, z, w) = \grie^{\TT, \phi, b}(y, x, w, z)
\end{equation*}
\begin{defn}\label{defn:phibRiemanniantensor}
 The \emph{$\varphib$-Riemannian curvature} $\grie^\phib$ for $(M, g)$ is the $\varphib$-curvature for $\gcon^{\phib}$, and the corresponding curvature tensor is the \emph{$\varphib$-Riemann tensor}, which is also denoted by $\grie^\phib$.
\end{defn}

\subsection{Bianchi identities}\label{subsec:bianchi}
Since $\gcon^\phib$ preserves $C_\pm^b$, $\grie^\phib(x, y, z, w)$ vanishes when the last two entries are sections of different $\GG^b$-eigenbundles. The non-vanishing components in the $\GG^b$-eigendecomposition of $\grie^\phib$ are given by the classical Riemann tensor $R$ of $g$ as well as the curvature tensors $R^{\pm\phi}$ of $\nabla^{\pm \phi}$.
\begin{prop}\label{prop:grieeigendecomposition}
 Let $X, Y, Z, W \in TM$ and $x^{b\pm} \in C_\pm^b$ with $X = \pi\pair{x^{b\pm}}$ etc, then
 \begin{enumerate}
  \item
  $\grie^\phib\pair{x^{b\pm}, y^{b\pm}, z^{b\pm}, w^{b\pm}} = R(X, Y, Z, W)$
  \item
  $\grie^\phib\pair{x^{b\mp}, y^{b\pm}, z^{b\pm}, w^{b\pm}} = R^{\pm\phi}(X, Y, Z, W) \mp \dfrac{1}{2}\pair{\nabla^{\pm\phi}_X\phi}\pair{Y, Z, W}$
  \item
  $\grie^\phib\pair{x^{b\mp}, y^{b\mp}, z^{b\pm}, w^{b\pm}} = R^{\pm\phi}(X, Y, Z, W)$
 \end{enumerate}
 All other components of $\grie^\phib$ vanish.
\end{prop}
\begin{proof}
 Standard computations from the definitions, which is left for the reader.
\end{proof}
By $(1)$ above, the algebraic Bianchi identity for $\grie^\phib$ holds when all entries involved are from the same $\GG^b$-eigenbundle. The analogues to the algebraic and differential Bianchi identities follow from previous discussions.
\begin{lemma}\label{lemma:bianchis}
 In their respective $\GG^b$-eigendecompositions:
 \begin{enumerate}
  \item
  For $\grie^{\phi, b}_{x, y} z + c.p.$:
  \begin{equation}\label{eq:firstbianchi}
   \grie^{\phi, b}_{x^{b\mp}, y^{b\mp}} z^{b\pm} + \grie^{\phi, b}_{y^{b\mp}, z^{b\pm}} x^{b\mp} + \grie^{\phi, b}_{z^{b\pm}, x^{b\mp}} y^{b\mp} = \pm2g\pair{R^{\pm\phi}_{X,Y} Z}
  \end{equation}
  \item
  For $\dd^\phib_\gcon \gcurv^\phib$:
  \begin{equation}\label{eq:secondbianchi}
   \pair{\dd^\phib_\gcon \gcurv^\phib}_{x^{b\mp}, y^{b\mp}, z^{b\pm}} = \mp 2\gcon_{g\pair{R^{\pm \phi}_{X, Y} Z}} = \mp 2\cdiff_{R^{\pm \phi}_{X, Y} Z}
  \end{equation}
 \end{enumerate}
 All other components vanish. 
\end{lemma}
\begin{proof}
 The identity \eqref{eq:firstbianchi} follows from \eqref{eq:jacobiatorcurvature} and \eqref{eq:jacobiatorphib}. If $d\phi = 0$, \eqref{eq:firstbianchi} can also be obtained from the explicit expressions in Proposition \ref{prop:grieeigendecomposition} and the identity \cite{Bismut:localindex:89} below
 \begin{equation}\label{eq:bismutdoubleswitch}
  R^{+\phi}(X, Y, Z, W) = R^{-\phi}(Z, W, X, Y)
 \end{equation}
 Then \eqref{eq:secondbianchi} follows from \eqref{eq:differentialbianchitensor} and \eqref{eq:firstbianchi}.
\end{proof}

In particular, the differential Bianchi identity holds for the $\varphib$-curvature when $\cdiff = 0$, i.e. $\gcon$ is the lifting of a classical connection on $V$. Another special case is when $\nabla^{\pm\phi}$ are flat (\cite{AgricolaFriedrich:FlatMetricConnection:09, CartanSchouten:26a, CartanSchouten:26b}), i.e. $R^{\pm \phi} = 0$, then both Bianchi identities hold. In this second special case, $\grie^\phib$ enjoys all the symmetries of a classical Riemann curvature.
\begin{theorem}\label{thm:classicalidentities}
Suppose that the connections $\nabla^{\pm\phi}$ are flat on $TM$. Then
 \begin{enumerate}
  \item
  $\grie^\phib(x, y, z, w) = - \grie^\phib(y, x, z, w) = \grie^\phib(y, x, w, z)$
  \item
  $\grie^\phib_{x, y}z + \grie^\phib_{y, z} x + \grie^\phib_{z, x} y = 0$
  \item
  $\dd^\phib_\gcon \grie^\phib = 0$
  \item
  $\grie^\phib(x, y, z, w) = \grie^\phib(z, w, x, y)$
 \end{enumerate}
 In this case, $\grie^\phib$ defines a symmetric pairing on $\wedge^2 \TT M$:
 \[\grie^\phib: \wedge^2\TT M \tensor \wedge^2\TT M \to \RR: \grie^\phib(x\wedge y, w \wedge z) := \grie^\phib(x, y, z, w)\]
 which defines the corresponding operator on $\wedge^2 \TT M$ via $\GG^b$.
\end{theorem}
\begin{proof}
 Items $(1)$ -- $(3)$ follow from preceeding discussions and the flatness assumption, while $(4)$ follows from $(1)$ -- $(3)$ as in the classical situation. The last statement is a consequence of $(1)$ and $(4)$.
\end{proof}

The following consequence of Lemma \ref{lemma:bianchis} was used in the proof of Proposition \ref{prop:topcohomology}.
\begin{prop}\label{prop:topdifferential}
 $\dd^\phib$ is a differential at $\Omega_\TT^{2n-1}\pair{M}$.
\end{prop}
\begin{proof}
Let $\bair{X_j}$ be a local $g$-orthonormal frame of $TM$ and $\bair{e_j^{b\pm}}$ the induced $\GG^b$-orthonormal frame of $\TT M$. Let $\theta \in \Omega_\TT^{2n-2}\pair{M}$, then by \eqref{eq:jacobiatorcurvature}, \eqref{eq:curvaturemap} and \eqref{eq:firstbianchi}:
\begin{equation*}
 \begin{split}
  & \pair{\dd^\phib \circ \dd^\phib \theta}\pair{e_1^{b+}, \ldots, e_n^{b+}, e_1^{b-}, \ldots, e_n^{b-}} \\
  = & - \sum_{i<j, k} \pair{-1}^{i+j+k+n} \theta\pair{- g\pair{R^{-\phi}_{X_i, X_j} X_k}, \ldots, \hat e_i^{b+}, \ldots, \hat e_j^{b+}, \ldots, \hat e_k^{b-}, \ldots} \\
  & - \sum_{i, j<k} \pair{-1}^{i+j+k} \theta\pair{g\pair{R^{+\phi}_{X_j, X_k} X_i}, \ldots,\hat e_i^{b+}, \ldots, \hat e_j^{b-}, \ldots, \hat e_k^{b-}, \ldots} \\
  = & \sum_{i,j,k} \pair{-1}^{j+k+n} \sair{ - R^{-\phi}\pair{X_i, X_j, X_k, X_i} + R^{+\phi}\pair{X_i, X_k, X_j, X_i}} \theta\pair{\ldots, \hat e_j^{b+}, \ldots, \hat e_k^{b-}, \ldots} \\
  = & 0
 \end{split}
\end{equation*}
The last step above follows from \eqref{eq:bismutdoubleswitch}.
\end{proof}

\subsection{Ricci curvature}\label{subsec:riccicurvature}
The trace of the $\varphib$-curvature on $\TT M$ defines the corresponding $\varphib$-Ricci curvature.
\begin{defn}\label{defn:standardRicci}
 For a $\GG^b$-metric connection $\gcon^\TT$, the \emph{$\varphib$-Ricci curvature} $\Ric^\TT: \TT M \to \TT M$ for $(M, \GG^b)$ is the trace of the $\varphib$-curvature $\grie^\TT$. For $x, y \in \se \TTM$, in the local orthonormal frame $\left\{e_i^{b+}, e_j^{b-}\right\}$ of $\TT M$ induced from a local $g$-orthonormal frame $\{X_i\}$
 \begin{equation}\label{eq:phibricci}
  \Ric^\TT(x) := \sum_{i}\sair{\grie^\TT_{x, e_i^{b+}} e_i^{b+} + \grie^\TT_{x, e_i^{b-}} e_i^{b-}}
 \end{equation}
 The \emph{$\varphib$-Ricci tensor} $\Rc^\TT \in \se{\TTM \tensor \TTM}$ is
 \begin{equation}\label{eq:phibriccitensor}
  \Rc^\TT(x, y) := \GG^b\pair{\Ric^\TT(x), y}
 \end{equation}
\end{defn}

Specialized to the $\varphib$-connection $\gcon^\phib$, the $\GG^b$-eigendecomposition of $\Rc^\phib$ can be determined from that of $\grie^\phib$ as follows:
\begin{equation}\label{eq:phibriccicomponents}
 \Rc^\phib\pair{x^{b\pm}, y^{b\pm}} = Rc\pair{X, Y} ~~~ \text{ and } ~~~ \Rc^\phib\pair{x^{b\pm}, y^{b\mp}} = Rc^{\mp\phi}\pair{X, Y}
\end{equation}
where $Rc$ is the Ricci tensor for $\nabla$, $Rc^{\pm\phi}$ are the Ricci tensors for $\nabla^{\pm\phi}$ respectively:
\begin{equation}\label{eq:torsionricci}
 Rc^{\pm\phi} = Rc \mp \dfrac{1}{2}d^*\phi - \dfrac{1}{4}\phi^2
\end{equation}
and
\[\phi^2(X, Y) := \sum_{i,j} \phi(X, X_i, X_j)\phi(Y, X_i, X_j)\]
It follows that $\Rc^\phib$ is symmetric, in other words
\begin{equation}\label{eq:phibriccisymmetry}
 \<\Ric^\phib(x), \GG^b y\> =\<\GG^b\Ric^\phib(x), y\> = \<\GG^b\Ric^\phib(y), x\>= \<\GG^b x, \Ric^\phib(y)\>
\end{equation}

Similar to the classical case, the $\varphib$-Ricci curvature appears in a Weitzenb\"ock identity relating two natural Laplacians on $\wedge^*\TT M$ described in \S \ref{subsec:laplacians}.
\begin{theorem}\label{thm:hodgelaplacianbochner}
 On $\Omega^*_\TT\pair{M}$, the following Weitzenb\"ock identity holds:
 \begin{equation}\label{eq:hodgelaplacianbochner}
  \laplacian^\phib = \laplacian_{\gcon^\phib} + \GG^b\Ric^\phib\GG^b
 \end{equation}
 where $\laplacian^\phib$ is the $\gcon^\phib$-Hodge Laplacian while $\laplacian_{\gcon^\phib}$ is the $\gcon^\phib$-Bochner Laplacian.
\end{theorem}
\begin{proof}
 It can be shown following standard computation that
 \[ \laplacian^\phib \theta = \laplacian_{\gcon^\phib} \theta - \half \sum_{\alpha, \beta}\GG^b\pair{e_\alpha^b} \wedge \iota_{e_\beta^b} \pair{\grie^\phib_{e_\alpha^b, e_\beta^b} \theta}\]
 where $e^b_\alpha, e^b_\beta$ run through $\{e^{b+}_i, e^{b-}_j\}$. Set $\theta = \GG^b\pair{x}$ for $x \in \TT M$, then
 \begin{equation*}
  \begin{split}
   & \sum_{\alpha, \beta}\<\GG^b\pair{e_\alpha^b} \wedge \iota_{e_\beta^b} \pair{\grie^\phib_{e_\alpha^b, e_\beta^b} \GG^b\pair{x}}, y\> = 2\sum_{\alpha, \beta}\< \GG^b\pair{\grie^\phib_{e_\alpha^b, e_\beta^b} x}, e_\beta^b\> \GG^b\pair{e_\alpha^b, y}\\
   = & 2\sum_{\beta}\grie^\phib\pair{y, e_\beta^b, x, e_\beta^b} = -2\Rc^\phib\pair{y, x} = -2\<\GG^b\Ric^\phib\pair{x}, y\>
  \end{split}
 \end{equation*}
 from which \eqref{eq:hodgelaplacianbochner} follows.
\end{proof}

The symmetry of $\Rc^\phib$ implies that $\dd^{\phib*}$ is a differential at $\Omega^1_\TT\pair{M}$. For $x_i \in \se{\TT M}$, $i = 1, \ldots, k$, set $y_i = \GG^b\pair{x_i}$, then straightforward computation gives
\begin{equation}\label{eq:dstarsquared}
 \begin{split}
  & \pair{\dd^{\phib*} \circ \dd^{\phib*}} \sair{{y_1} \wedge \ldots \wedge {y_k}}\\
  = & \sum_{i<j} \pair{-1}^{i+j} \sair{\Rc^\phib\pair{x_j, x_i} - \Rc^\phib\pair{x_i, x_j}} {y_1} \wedge \ldots \wedge \hat{{y_i}} \wedge \ldots \wedge \hat{{y_j}} \wedge \ldots  \wedge {y_k}\\
  & - \sum_{i<j<\ell} \pair{-1}^{i+j+\ell}\GG^b\pair{\grie^\phib_{x_i, x_j} x_\ell + c.p.}  {y_1} \wedge \ldots \wedge \hat{{y_i}} \wedge \ldots \wedge \hat{{y_j}} \wedge \ldots \wedge \hat{{y_\ell}} \wedge \ldots \wedge {y_k} 
 \end{split}
\end{equation}
Due to the symmetry of $\Rc^\phib$, the terms in the second line vanish. Furthermore, when $k = 2$, the terms in the last line vanish as well.
\begin{prop}\label{prop:chaincomplex}
 The operator $\dd^{\phib*}$ is a differential at $\Omega^1_\TT\pair{M}$, i.e.
 \begin{equation}\label{eq:dstarat1}
  \theta \in \Omega^2_\TT\pair{M} \Longrightarrow \pair{\dd^{\phib*} \circ \dd^{\phib*}} \theta = 0
 \end{equation}
 If the connections $\nabla^{\pm\phi}$ are flat on $TM$, then $\pair{\Omega^*_\TT \pair{M}, \dd^\phib}$ and $\pair{\Omega^*_\TT \pair{M}, \dd^{\phib*}}$ are both chain complices.
\end{prop}
\begin{proof}
 For the last statement, that $\pair{\Omega^*_\TT\pair{M}, \dd^\phib}$ is a chain complex follows from \eqref{eq:dsquared}, \eqref{eq:jacobiatorcurvature} and the algebraic Bianchi identity, which is item $(2)$ in Theorem \ref{thm:classicalidentities}. The statement for $\pair{\Omega^*_\TT \pair{M}, \dd^{\phib*}}$ follows from \eqref{eq:dstarsquared} and the algebraic Bianchi identity.
\end{proof}

\subsection{Bismut connection}\label{subsec:bismutconnection}
The \emph{generalized Bismut connection $\gcon^\bismut$} introduced by Gualtieri \cite{Gualtieri:Poisson:10} is a $\GG^b$-metric connection and is the lift of a classical connection $\nabla^\bismut$ on $\TT M$, which depends on both $\phi$ and $b$:
\begin{equation}\label{eq:genbismut}
 \nabla^\bismut_X y^{b+} := \nabla^{+\phi}_X Y + \pair{b+g}\pair{\nabla^{+\phi}_X Y} ~~~ \text{ and } ~~~ \nabla^\bismut_X y^{b-} := \nabla^{-\phi}_X Y + \pair{b-g}\pair{\nabla^{-\phi}_X Y}
\end{equation}
Note that $\gcon^\bismut$ is compatible with the (almost) Dorfman bracket $*_{\phi - db}$ (Definition \ref{defn:diamondmetriccompatible}).

Since $\gcon^\bismut$ is a lift of a classical connection, by \eqref{eq:secondbianchi}, the \emph{$\varphib$-Bismut curvature} $\grie^\bismut$ of $\gcon^\bismut$ satisfies the differential Bianchi identity:
\[\dd^\phib_\gcon \grie^\bismut = 0\]
More explicitly, $\grie^\bismut$ is determined by the classical curvature of $\nabla^\bismut$, which in turn is given by $R^{\pm \phi}$:
\begin{equation}\label{eq:bismutriemann}
 \grie^\bismut\pair{x, y, z^{b\pm}, w^{b\pm}} = R^{\pm\phi}\pair{X,Y,Z,W}
\end{equation}
The $\GG^b$-eigendecomposition of the corresponding \emph{$\varphib$-Bismut Ricci tensor} $\Rc^\bismut$ is thus
\begin{equation}\label{eq:bismutriccieigen}
 \Rc^\bismut\pair{x, y^{b+}} = Rc^{+\phi}\pair{X, Y} \text{ and } \Rc^\bismut\pair{x, y^{b-}} = Rc^{-\phi}\pair{X, Y}
\end{equation}

\subsection{Scalar curvatures} \label{subsec:scalarcurvature}
The traces of the $\varphib$-Ricci curvatures give the corresponding \emph{$\varphib$-scalar curvature}s, which depend on the $\GG^b$-metric connection $\gcon^\TT$. For instance, the 
\emph{$\varphib$-Riemann scalar curvature} $\gsca^\phib$ is the trace of $\Ric^\phib$:
\begin{equation}\label{eq:scalarcurvature}
 \gsca^\phib = \sum_j \sair{\Rc^\phib\pair{e^{b+}_i, e^{b+}_i} + \Rc^\phib\pair{e^{b-}_i, e^{b-}_i}} = 2S
\end{equation}
where $S$ is the classical scalar curvature of $g$. On the other hand, the \emph{$\varphib$-Bismut scalar curvature} $\gsca^\bismut$ is the trace of $\Ric^\bismut$:
\[\gsca^\bismut = \sum_j \sair{\Rc^\bismut\pair{e^{b+}_i, e^{b+}_i} + \Rc^\bismut\pair{e^{b-}_i, e^{b-}_i}} = 2S - 3\norm{\phi}_g^2\]
where $\norm{\phi}_g$ is the norm of $\phi$ with respect to $g$:
\[\norm{\phi}_g^2 = \sum_{i<j<k}\phi\pair{X_i,X_j,X_k}^2\]

\section{Generalized complex manifolds}\label{sec:gencomplex}

Let $\JJ$ be a generalized almost complex structure on $M$. It induces a polarization of $\TT_\CC M : =\TT M \tensor_\RR \CC$ as the direct sum of its $\pm \sqrt{-1}$-eigenbundles:
\begin{equation}\label{eq:complexpolarization}
 \TT_\CC M = \TT_\JJ^{1,0} M \dsum \TT_\JJ^{0,1} M
\end{equation}
Here $\TT_\JJ^{1,0} M$ denotes the $\sqrt{-1}$-eigenbundle of $\JJ$, and $\TT^{0,1}_\JJ M$ its complex conjugate. They are maximally isotropic and are dual to each other under the pairing $2\<,\>$ on $\TT_\CC M$. For instance, the space of \emph{$\pair{0,1}$-forms with respect to $\JJ$} is identified with the sections of $\TT_\JJ^{1,0} M$:
\[\Omega^{0,1}_\JJ\pair{M} := \se{\TT_\JJ^{1,0} M}\]
Similar to the classical case, the $\pair{0,k}$-forms with respect to $\JJ$ are sections of $\wedge^k\TT_\JJ^{1,0} M$:
\[\Omega^{0,k}_\JJ\pair{M} := \se{\wedge^k\TT_\JJ^{1,0} M}\]
In general, the type decomposition of $\Omega_\TT^*\pair{M}$ with respect to $\JJ$ is given by
\begin{equation}\label{eq:Jtypedecomposition}
 \Omega^*_\TT\pair{M} = \dsum_{p,q} \Omega^{p,q}_\JJ \pair{M} := \dsum_{p,q} \se{\wedge^p \TT^{0,1}_\JJ M \tensor \wedge^q \TT^{1,0}_\JJ M} 
\end{equation}
For notational convenience, sometimes $\TT_\JJ^{1,0} M$ is denoted $L$, while $\TT_\JJ^{0,1} M$ is denoted $\bar L$.

\begin{defn}\label{defn:diamondgcplxcompatible}
 Let $\pair{M, \JJ}$ be a generalized almost complex manifold and let $\gamma \in \Omega^3\pair{M}$ be a closed $3$-form. Then a generalized connection $\gcon^\TT$ on $\TT M$ is \emph{$\JJ$-compatible with $*_\gamma$} if $\diamond_\TT$ coincides with $*_\gamma$ on the sections from the same eigenbundle of $\JJ$, i.e.
 \begin{equation}\label{eq:diamondgcplxcompatible}
  x \diamond_\TT y = x *_\gamma y \text{ and } \bar x \diamond_\TT \bar y = \bar x *_\gamma \bar y
 \end{equation}
 where $x, y \in \se{\TT^{1,0}_\JJ M}$. Such generalized connection $\gcon^\TT$ is a \emph{$\gamma$-$\JJ$-connection} if it furthermore is $TM$-torsion free.
\end{defn}
\noindent
It is straightforward to see that \eqref{eq:diamondgcplxcompatible} is equivalent to the following, where $x, y \in \se{\TT M}$:
\[\pair{\JJ x} \diamond_\TT y + x \diamond_\TT \pair{\JJ y} = \pair{\JJ x} *_\gamma y + x *_\gamma \pair{\JJ y}\]
Thus, if non-empty, the space $\ms D_{\JJ, \gamma}\pair{\TT M}$ of generalized connections that are $\JJ$-compatible with $*_\gamma$ is modeled on
\[\ms D_{\JJ, \gamma}\pair{\TT M} \cong \bair{A : \TT M \tensor \TT M \to \TT M : A_{\JJ x} y - A_y\pair{\JJ x} = A_{\JJ y} x - A_x \pair{\JJ y}}\]
It's then evident by \eqref{eq:Ttorionfreespace} that the subspace $\ms D_{\JJ, \gamma, \tau}\pair{\TT M}$ of the $\gamma$-$\JJ$-connections, if non-empty, is modeled on
\begin{equation}\label{eq:gammaJconnectionspace}
 \ms D_{\JJ, \gamma, \tau}\pair{\TT M} \cong \bair{A : \TT M \odot \TT M \to T^*M : A_{\JJ x} y - A_y\pair{\JJ x} = A_{\JJ y} x - A_x \pair{\JJ y}}
\end{equation}

Let $\gamma \in \Omega^3\pair{M}$ be a closed $3$-form, then $\JJ$ is \emph{integrable (with respect to $\gamma$)} if $\TT^{1,0}_\JJ M$ is involutive under the Dorfman bracket $*_\gamma$, i.e.
\[x, y \in \se{\TT^{1,0}_\JJ M} \Longrightarrow x *_\gamma y \in \se{\TT^{1,0}_\JJ M}\]
In this case, $\pair{M, \gamma; \JJ}$ is a \emph{generalized complex manifold}.

\begin{lemma}\label{lemma:dphibdecompose}
 Let $\pair{M, \gamma; \JJ}$ be a generalized complex manifold with a generalized connection $\gcon^\TT$ on $\TT M$ that is $\JJ$-compatible with $*_\gamma$. Then the operator $\dd^\TT$ decomposes into components as follows:
 \begin{equation}\label{eq:Jdecomposition}
  \dd^\TT = \partial^\TT_\JJ + \bar\partial^\TT_\JJ: \Omega^{p,q}_\JJ \pair{M} \to \Omega^{p+1, q}_\JJ \pair{M} \dsum \Omega^{p, q+1}_\JJ \pair{M}
 \end{equation}
\end{lemma}
\begin{proof}
 The proof is similar to that of Lemma \ref{lemma:plusminusdecomposition}, employing \eqref{eq:diamondgcplxcompatible} and the integrability of $\JJ$. The details are left for the reader.
\end{proof}

The integrability of $\JJ$ implies that both of its eigenbundles are complex Lie algebroids, with their Lie brackets given by the restriction of $*_\gamma$. 
The corresponding Lie algebroid de Rham differential for $\bar L = \TT_\JJ^{0,1} M$ will be denoted by $d_{\bar L}$:
\begin{equation}\label{eq:liealgebroiddifferential}
 d_{\bar L}: \Omega^{0,k}_\JJ\pair{M} \to \Omega^{0,k+1}_\JJ\pair{M}
\end{equation}
If $\gcon^\TT$ is $\JJ$-compatible with $*_\gamma$, then $d_{\bar L}$ coincides with the restriction of $\dd^\TT$ on $\Omega^{0,*}_\JJ \pair{M}$. More precisely, for $\theta \in \Omega^{0,q}_\JJ \pair{M}$, direct computation shows that
\begin{equation}\label{eq:gcplxdrestriction}
 \pair{d_{\bar L} \theta}\pair{\bar x_0, \bar x_1, \ldots, \bar x_q} = \sum_j (-1)^j \pair{\gcon^{\TT}_{\bar x_j} \theta}\pair{\bar x_0, \ldots, \hat{\bar x_j}, \ldots, \bar x_q}
\end{equation}
where $x_j \in \se{\TT^{1,0}_\JJ M}$ for all $j$, which gives
\begin{equation}\label{eq:dTcomponentbar}
 \bar\partial^\TT_\JJ \theta = d_{\bar L} \theta \in \Omega^{0, q+1}_\JJ\pair{M}
\end{equation}
The other component $\partial^\TT_\JJ \theta \in \Omega^{1, q}_\JJ\pair{M}$ is given by
\begin{equation}\label{eq:dTJcomponent}
 \pair{\partial^\TT_\JJ \theta}\pair{x_0, \bar x_1, \ldots, \bar x_q} = \pair{\gcon^{\TT}_{x_0} \theta}\pair{\bar x_1, \ldots, \bar x_q} + \sum_{j} \theta\pair{\bar x_1, \ldots, \gcon^\TT_{\bar x_j} x_0, \ldots, \bar x_q}
\end{equation}
where $x_j \in \se{\TT^{1,0}_\JJ M}$. Via complex conjucation, the analogous versions of the identities \eqref{eq:gcplxdrestriction}, \eqref{eq:dTcomponentbar} and \eqref{eq:dTJcomponent} are valid for $\bar \theta \in \Omega^{q,0}_\JJ\pair{M}$ with $L = \TT^{1,0}_\JJ M$ in place of $\bar L$.

\begin{example}\label{ex:ddbarfgeneral}
 Suppose that $\gcon^\TT$ is $\JJ$-compatible with $*_\gamma$ and consider $f \in \se{M}$. It follows from Proposition \ref{prop:operatorrelations} and Lemma \ref{lemma:dphibdecompose} that
 \begin{equation}\label{eq:vanishingcomponents}
  \partial^\TT_\JJ \partial^\TT_\JJ f = 0, \bar \partial^\TT_\JJ \bar\partial^\TT_\JJ f = 0, \text{ and } \partial^\TT_\JJ \bar\partial^\TT_\JJ f + \bar \partial^\TT_\JJ \partial^\TT_\JJ f = 0
 \end{equation}
 The first two identities in the above can also be seen as the consequences of \eqref{eq:dTcomponentbar} and the corresponding version for $d_L$.
\end{example}

\subsection{Generalized (almost) Hermitian manifolds}\label{subsec:almosthermitian}
Let $\pair{M; \GG^b, \JJ}$ be a \emph{generalized almost Hermitian manifold}, i.e. $\JJ$ and $\JJ_- := \GG^b \JJ$ are commuting generalized almost complex structures. The eigenbundles of $\JJ$ (and $\JJ_-$) decompose into the common eigenbundles of $\GG^b$ and $\JJ$. Let
\begin{equation}\label{eq:commoneigen}
 \ell_\pm := \TT_\JJ^{1,0} M \inter \pair{C_\pm^b \tensor \CC}
\end{equation}
then for instance
\[\TT_\JJ^{1,0} M = \ell_+ \dsum \ell_- \text{ and } \TT_{\JJ_-}^{1,0} M = \ell_+ \dsum \bar \ell_-\]
The restriction of $\JJ$ to $C^b_\pm$ induces two almost complex structures $I_\pm$ on $TM$:
\begin{equation}\label{eq:inducedalmostcomplex}
 I_\pm X+ \pair{b\pm g}\pair{I_\pm X} := \JJ_+\sair{X+\pair{b\pm g}X}
\end{equation}
It follows that $\pair{M; \GG^b, \JJ}$ is equivalent to a pair of almost Hermitian structures $\pair{g, I_\pm}$ together with $b \in \Omega^2(M)$.

\begin{lemma}\label{lemma:hermitiancompatible}
 For a generalized almost Hermitian manifold $\pair{M; \GG^b, \JJ}$, the space of $\GG^b$-metric $\gamma$-$\JJ$-connections, if non-empty, is modeled on the following subspace of $\Omega^3\pair{M}$:
 \[\bair{\phi \in \Omega^3\pair{M}: \phi\pair{I_- X, Y, Z} + \phi\pair{X, I_+ Y, Z} = 0, \text{ for all } X, Y, Z \in \se{TM}}\]
\end{lemma}
\begin{proof}
 It follows from \eqref{eq:gammaJconnectionspace} and Theorem \ref{thm:torsionlessmetric}.
\end{proof}

Let $\phi = \gamma + db$, then by Theorem \ref{thm:levicivita}, $\gcon^\phib$ is metric compatible with $*_\gamma$. Hence, the $\JJ$-compatibility of $\gcon^\phib$ with $*_\gamma$ is equivalent to $\diamond_\phib$ and $*_\gamma$ coincide on sections of the same common eigenbundle of $\JJ$ and $\GG^b$, e.g.
\begin{equation}\label{eq:Jcompatibledorfman}
 x^{b\pm} *_\gamma y^{b\pm} = x^{b\pm} \diamond_\phib y^{b\pm}
\end{equation}
where $x^{b\pm}, y^{b\pm} \in \se{\ell_\pm}$ with $\pi\pair{x^{b\pm}} = X$ and so on. 
\begin{lemma}\label{lemma:phibJcompatiblealmosthermitian}
 On a generalized almost Hermitian manifold $\pair{M; \GG^b, \JJ}$, let $\gamma \in \Omega^3\pair{M}$ be a closed $3$-form and $\phi = \gamma + db$. Then $\gcon^\phib$ is $\JJ$-compatible with $*_\gamma$ iff $\nabla^{\pm \phi} I_\pm = 0$.
\end{lemma}
\begin{proof}
 For any $X, Y, Z \in \se{TM}$, let $x^{b\pm} \in \se{C_\pm^b}$ such that $\pi\pair{x^{b\pm}} = X$ and so on, then
 \[\aair{x^{b\pm} *_\gamma y^{b\pm} - x^{b\pm} \diamond_\phib y^{b\pm}, z} = \pm g\pair{\nabla^{\pm\phi}_Z X, Y}\]
 Thus \eqref{eq:Jcompatibledorfman} is equivalent to
 \[g\pair{\nabla^{\pm\phi}_Z X_\pm, Y_\pm} = 0\]
 for all $Z \in \se{TM}$ and $X_\pm, Y_\pm \in \se{T_{\pm; 1,0} M}$ respectively. Hence $\nabla^{\pm\phi}$ preserves $T_{\pm; 1,0} M$ respectively, from which the statement follows.
\end{proof}

\begin{remark}\label{remark:GrayHervella}
 The condition $\nabla^{\pm\phi} I_\pm = 0$ can be rewritten as
 \begin{equation}\label{eq:levicivitaIpm}
  \pair{\nabla_X I_\pm} Y = \pm \half \pair{I_\pm g\inverse \iota_Y\iota_X \phi - g\inverse \iota_{I_\pm Y} \iota_X \phi}
 \end{equation}
 for $X, Y \in \se{TM}$, which implies that $I_\pm$ are of class $\WWC_1 \dsum \WWC_3 \dsum \WWC_4$ in the Gray-Hervella classification of almost Hermitian structures \cite{GrayHervella:SixteenClasses:80}. The Nijenhuis tensors of $I_\pm$
 \[N_{I_\pm}(X,Y) := [X, Y] + I_\pm[I_\pm X, Y] + I_\pm[X, I_\pm Y] - [I_\pm X, I_\pm Y]\]
 can be expressed in terms of $\phi$:
 \[g\pair{N_{I_\pm}\pair{X, Y}, Z} = \phi\pair{I_\pm X, I_\pm Y, Z} + \phi\pair{I_\pm X, Y, I_\pm Z} + \phi\pair{X, I_\pm Y, I_\pm Z} - \phi\pair{X, Y, Z}\]
 which implies that the integrability of $I_\pm$ is equivalent to $\phi$ being of type $\pair{1,2} + \pair{2,1}$ with respect to $I_\pm$ respectively. Furthermore, by Friedrich-Ivanov \cite{FriedrichIvanov:connectionwithskewtorsion:02}, if almost Hermitian connections for $\pair{g, I_\pm}$ that admit completely skew torsions exist, they must be unique. 
\end{remark}

In $\pair{M; \GG^b, \JJ}$, when $\JJ$ is integrable with respect to $\gamma\in \Omega^3(M)$, the structure $\pair{M, \gamma; \GG^b, \JJ}$ defines a \emph{generalized Hermitian manifold}. In this case, the $\JJ$-compatibility of $\gcon^\phib$ with $*_\gamma$ is equivalent to the integrability of $\JJ_-$, i.e. it provides a \emph{generalized K\"ahler condition}.
\begin{theorem}\label{thm:genkahlercondition}
 On a generalized Hermitian manifold $\pair{M, \gamma; \GG^b, \JJ}$, let $\phi = \gamma + db$. Then $\pair{M, \gamma; \GG^b, \JJ}$ is a generalized K\"ahler manifold iff $\gcon^\phib$ is $\JJ$-compatible with $*_\gamma$.
\end{theorem}
\begin{proof}
 Starting with the $\JJ$-compatibility, then \eqref{eq:Jcompatibledorfman} and the integrability of $\JJ$ imply that $\ell_\pm$ and $\bar\ell_\pm$ are involutive with respect to $*_\gamma$, which further implies that $I_\pm$ are integrable. Let $x^{b+} \in \se{\ell_+}$, $y^{b-}\in \se{\ell_-}$, $X_+ = \pi\pair{x^{b+}}$ and 
 $Y_- = \pi\pair{y^{b-}}$, 
 then by \eqref{eq:phibconnectionexplicit} and Lemma \ref{lemma:phibJcompatiblealmosthermitian}
 \[x^{b+} *_\gamma \bar{y^{b-}} = \sair{\nabla^{-\phi}_{X_+} \bar Y_- + \pair{b-g} \nabla^{-\phi}_{X_+} \bar Y_-} - \sair{\nabla^{+\phi}_{\bar Y_-} X_+ + \pair{b+g} \nabla^{+\phi}_{\bar Y_-} X_+} \in \ell_+ \dsum \bar \ell_-\]
 Thus $\JJ_-$ is integrable. The opposite direction is left to the reader.
\end{proof}

\begin{remark}\label{remark:genkahlercondition}
 The generalized K\"ahler condition in Theorem \ref{thm:genkahlercondition} relates to the condition given in \cite{Gualtieri:Poisson:10} as follows. The condition $\nabla^{\pm \phi} I_\pm = 0$ on a generalized almost Hermitian manifold is equivalent to $\gcon^\bismut \JJ = 0$, which implies the equivalence of the integrability of $I_\pm$ to the type condition for the \emph{generalized torsion} as defined in \cite{Gualtieri:Poisson:10} -- the integrability of $\JJ$ then follows. In Theorem \ref{thm:genkahlercondition}, the integrability of $I_\pm$ follows from the integrability of $\JJ$ and \eqref{eq:Jcompatibledorfman}, obtaining the type of $\phi$ with respect to $I_\pm$ as a consequence.
\end{remark}

\begin{corollary}\label{coro:genkahlercondition}
 On a generalized Hermitian manifold $\pair{M, \gamma; \GG^b, \JJ}$, let $\phi = \gamma + db$. Then $\pair{M, \gamma; \GG^b, \JJ}$ is a generalized K\"ahler manifold iff $\gcon^\bismut \JJ = 0$.
\end{corollary}
\begin{proof}
 As stated in Remark \ref{remark:genkahlercondition}, $\gcon^\bismut \JJ = 0$ is equivalent to $\nabla^{\pm \phi} I_\pm = 0$, which by Lemma \ref{lemma:phibJcompatiblealmosthermitian} is equivalent to $\gcon^\phib$ being $\JJ$-compatible.
\end{proof}

Let $\gcon^\TT$ be any $\GG^b$-metric connection on a generalized Hermitian manifold $\pair{M, \gamma; \GG^b, \JJ}$. Its \emph{$\JJ$-Ricci form} $\rho_\JJ^\TT := \rho_\JJ\pair{\gcon^\TT} \in \Omega^2_\TT(M)$ is defined as
\[\rho_\JJ^\TT\pair{x, y} := \sum_i\sair{\grie^\TT\pair{x, y, \JJ e^{b+}_i, e^{b+}_i} + \grie^\TT\pair{x, y, \JJ e^{b-}_i, e^{b-}_i}}\]
The \emph{$\JJ$-scalar curvature} for $\gcon^\TT$ is
\[\gsca_\JJ^\TT := \sum_i\sair{\rho_\JJ^\TT\pair{\JJ e^{b+}_i, e^{b+}_i} + \rho_\JJ^\TT\pair{\JJ e^{b-}_i, e^{b-}_i}}\]

\subsection{Generalized K\"ahler manifolds}\label{subsec:genKahler}
Recall that the structure $\pair{M, \gamma; \GG^b, \JJ}$ defines a \emph{generalized K\"ahler manifold} if both $\JJ$ and $\JJ_- = \GG^b \JJ$ are integrable generalized almost complex structures with respect to ${\gamma}$ (\cite{Gualtieri:genKahler:14}). Let $\phi = \gamma + db$, Theorem \ref{thm:genkahlercondition} indicates that $\gcon^\phib$ is $\JJ$-compatible with $*_\gamma$. In particular, for $\theta \in \Omega^{p,0}_\JJ\pair{M}$, \eqref{eq:dTcomponentbar} gives
\begin{equation}\label{eq:dphibcomponents}
 \partial^\phib_\JJ \theta = d_{L} \theta \text{ and } \bar \partial^\phib_\JJ \bar\theta = d_{\bar L} \bar \theta
\end{equation}

\begin{example}[c.f. Example \ref{ex:ddbarfgeneral}]\label{ex:ddbarf}
 Set $\gcon^\TT = \gcon^\phib$ in \eqref{eq:vanishingcomponents}, then the components in the third identity can further be rewritten in terms of the $\GG^b$-eigendecomposition. For instance
 \begin{equation}\label{eq:ddbarf}
  \pair{\partial^\phib_\JJ \bar\partial^\phib_\JJ f}\pair{x_\pm, \bar y_\pm} = \pair{\partial_\pm \bar\partial_\pm f}\pair{X_\pm, \bar Y_\pm}
 \end{equation}
 where $\partial_\pm$ and $\bar\partial_\pm$ are the operators associated to the classical complex structures $I_\pm$, while $X_\pm, Y_\pm \in \se{T_{\pm; 1,0} M}$ are sections of the $I_\pm$-holomorphic tangent bundles respectively, and $x_\pm = X_\pm + \pair{b\pm g}X_\pm$ and so on.
\end{example}

Since $d_{\bar L}^2 = 0$, as a consequence of \eqref{eq:dphibcomponents}, the algebraic Bianchi identity \eqref{eq:firstbianchi} for $\grie^\phib$ implies that $\nabla^{+ \phi}$ (resp. $\nabla^{-\phi}$) induces a natural $I_-$-holomorphic (resp. $I_+$-holomorphic) structure on the eigenbundles of $I_+$ (resp. of $I_-$), providing an alternative proof of this well-known result in \cite{Gualtieri:Poisson:10}.
\begin{prop}\label{prop:LfirstBianchi}
 Let $(M, \gamma; \GG^b, \JJ)$ be a generalized K\"ahler manifold. Let $\phi = \gamma + db$, $x, y, z, w\in \se{\TT^{1,0}_\JJ M}$, and $\bar x$, etc. are their complex conjugates. Then
 \begin{equation}\label{eq:LfirstBianchi}
  \aair{\grie^{\phib}_{\bar x, \bar y} \bar z, w} + c.p. \text{ in } x, y, z = 0 \text{ and }
  \aair{\grie^{\phib}_{x, y} z, \bar w} + c.p. \text{ in } x, y, z = 0
 \end{equation}
 which give rise to
 \[R^{-\phi}_{\bar X_+, \bar Y_+} \bar Z_- = 0, R^{-\phi}_{\bar X_+, \bar Y_+} W_- = 0, R^{+\phi}_{\bar X_-, \bar Y_-} \bar Z_+ = 0 \text{ and } R^{+\phi}_{\bar X_-, \bar Y_-} W_+ = 0\]
 as well as their complex conjugates, where $X_\pm, Y_\pm, Z_\pm, W_\pm \in \se{T_{\pm; 1,0} M}$.
\end{prop}
\begin{proof}
 To see the first identity in \eqref{eq:LfirstBianchi}, note that $w \in \Omega^{0,1}_\JJ\pair{M}$, then 
 \[\bar\partial^\phib_\JJ \circ \bar\partial^\phib_\JJ w = d_{\bar L}^2 w = 0\]
 By \eqref{eq:jacobiatorcurvature} and \eqref{eq:curvaturemap}, for $x, y, z \in \se{\TT^{1,0}_\JJ M}$
 \begin{equation*}
  \begin{split}
   & 0 = \pair{\bar\partial^\phib_\JJ \circ \bar\partial^\phib_\JJ w} \pair{\bar x, \bar y, \bar z} = \pair{\dd^\phib \circ \dd^\phib w} \pair{\bar x, \bar y, \bar z} \\
   = & w\pair{\sair{\bar x\diamond_\phib \bar y\diamond_\phib \bar z}} = - 2 \aair{\grie^{\phib}_{\bar x, \bar y} \bar z, w} - c.p. \text{ in } x, y, z
  \end{split}
 \end{equation*}
 The second identity then follows by taking complex conjugation.
 
 To see the classical curvature identities, restrict $x, y, z, w$ in \eqref{eq:LfirstBianchi} further to the $\GG^b$-eigenbundles \eqref{eq:commoneigen}. For instance, consider $x, y \in \se{\ell_+}$ and $z, w \in \se{\ell_-}$. Set $X_+ = \pi\pair{x}$, $Z_- = \pi\pair{z}$ and so on, then apply \eqref{eq:firstbianchi} to the first identity in \eqref{eq:LfirstBianchi} gives
 \[- \aair{2 g\pair{R^{-\phi}_{\bar X_+, \bar Y_+} \bar Z_-}, w_-} = g\pair{R^{-\phi}_{\bar X_+, \bar Y_+} W_-, \bar Z_-} = 0\]
 Since $\pair{g, I_-}$ is Hermitian and $\nabla^{-\phi}$ preserves $I_-$, it gives the identities involving $R^{-\phi}$. The rest of the identities are obtained similarly.
\end{proof}

For the generalized Bismut connection $\gcon^\bismut$, with $\phi = \gamma + db$,  \eqref{eq:bismutriemann} implies that $\rho^\bismut_\JJ$ is computed by the sum of the respective classical Bismut-Ricci forms $\rho_\pm$ for $(g, I_\pm)$:
\[\rho_\JJ^\bismut\pair{x, y} = \rho_+\pair{X, Y} + \rho_-\pair{X, Y}\]
The corresponding $\JJ$-scalar curvature also decomposes:
\[\gsca_\JJ^\bismut = \sum_i\sair{\rho_\JJ^\bismut\pair{\JJ e^{b+}_i, e^{b+}_i} + \rho_\JJ^\bismut\pair{\JJ e^{b-}_i, e^{b-}_i}} = S_+ + 2S_{+-} + S_-\]
where $S_\pm$ are the respective classical Bismut scalar curvatures for $\pair{g, I_\pm}$ and $S_{+-}$ is the \emph{mixed Bismut scalar curvature}:
\[S_{+-} := \sum_{i,j} R^{+\phi}\pair{I_-X_i, X_i, I_+X_j, X_j} = \sum_{i,j} R^{-\phi}\pair{I_+X_i, X_i, I_-X_j, X_j}\]

\section{$\JJ$-holomorphic vector bundles}\label{sec:Jholobundles}

Let $\pair{V, h}$ be a Hermitian vector bundle on a generalized complex manifold $\pair{M, \gamma; \JJ}$. Recall that a \emph{$\JJ$-holomorphic structure} on $V$ is given by a flat $\TT_\JJ^{0,1}M$-connection:
\[\bar\partial_{\JJ} : \se V \to \Omega^{0,1}_\JJ\pair{V} := \se{\TT_\JJ^{1,0} M \tensor V} : \bar\partial_{\JJ} \pair{fv} = d_{\bar L} f\tensor v + f \bar\partial_{\JJ} v\]
for $f \in \se M$ and $v \in \se V$, such that 
\begin{equation}\label{eq:jholostructure}
 \bar\partial_\JJ \circ \bar\partial_\JJ = 0
\end{equation}
where the extension to $\Omega^{0,k}_\JJ\pair{V}$ is given by
\begin{equation}\label{eq:extensionpartialbar}
 \bar\partial_\JJ: \Omega^{0,k}_\JJ\pair{V} \to \Omega^{0,k+1}_\JJ\pair{V} : \bar\partial_\JJ\pair{\theta \tensor v} = d_{\bar L}\theta \tensor v + \pair{-1}^k \theta \wedge \bar\partial_\JJ v
\end{equation}
for $\theta \in \Omega^{0,k}_\JJ\pair{M}$ and $v \in \se V$.

If $\pair{M, \gamma; \GG^b, \JJ}$ is generalized K\"ahler, 
via the restriction to $C^b_\pm$, the $\JJ$-holomorphic structure $\bar\partial_\JJ$ induces on $V$ an $I_\pm$-holomorphic structure, which will be denoted $\bar\partial_\pm$ respectively:
\begin{equation}\label{eq:inducedIpmstructure}
 \bar\partial_{\pm, \bar X_\pm} v := \bar\partial_{\JJ, \bar x_\pm} v
\end{equation}
where $X_\pm \in T_{\pm; 1,0}M$, $x_\pm = X_\pm + \pair{b\pm g}X_\pm$ and $v \in \se{V}$.

\subsection{Connection on $\TT_\JJ^{1,0} M$}\label{subsec:genholotangentbundle}

Consider a generalized K\"ahler manifold $\pair{M, \gamma; \GG^b, \JJ}$, a natural $\bar L$-connections can be defined on $\TT_\JJ^{1,0} M$ using the diamond bracket $\diamond_\TT$ associated to any $\GG^b$-metric connection $\gcon^\TT$:
\begin{equation}\label{eq:pconnectiononL}
 \bar\partial^\TT_\diamond: \se{\TT_\JJ^{0,1} M} \tensor \se{\TT_\JJ^{1,0} M} \to \se{\TT_\JJ^{1,0} M}: \bar\partial^\TT_{\diamond, \bar x} y := \sair{\bar x \diamond_\TT y}_{1,0}
\end{equation}
where $x, y \in \se{\TT_\JJ^{1,0} M}$ and $\sair{\bullet}_{1,0}$ denotes taking the $(1,0)$-component with respect to $\JJ$.

It follows from Proposition \ref{prop:LfirstBianchi} that $\bar\partial_\diamond^\phib$ defined by $\gcon^\phib$ via \eqref{eq:pconnectiononL} induces an $I_\pm$-holomorphic structure on $\TT_\JJ^{1,0}M$ via its restriction to $\bar\ell_\pm \cong T_{\pm; 0,1} M$. 
Moreover, by \eqref{eq:jacobiatorcurvature} and \eqref{eq:Jcompatibledorfman} 
\begin{equation*}
 \begin{split}
  \pair{\bar\partial_{\diamond}^\phib \circ \bar\partial_{\diamond}^\phib}_{\bar x, \bar y} z = & 
  \sair{\bar x \diamond_\phib \pair{\bar y \diamond_\phib z} - \bar y \diamond_\phib \pair{\bar x \diamond_\phib z} - \pair{\bar x \diamond_\phib \bar y} \diamond_\phib z}_{1,0}\\
  = & \sair{\grie^\phib_{\bar x, \bar y} z + \grie^\phib_{\bar y, z} \bar y + \grie^\phib_{z, \bar y} \bar x}_{1,0}
 \end{split}
\end{equation*}
Thus, from \eqref{eq:firstbianchi}, $\bar\partial_\diamond^\phib$ defines a $\JJ$-holomorphic structure on $\TT_\JJ^{1,0} M$ iff 
\[\sair{g\pair{R^{-\phi}_{\bar X_+, Z_+} \bar Y_-}}_{1,0} = \sair{g\pair{R^{+\phi}_{\bar Y_-, Z_-} \bar X_+}}_{1,0} = 0\]
for all $X_\pm, Y_\pm, Z_\pm \in \se{T_{\pm; 1,0} M}$. Since $\nabla^{\pm\phi}$ preserves $I_\pm$, together with \eqref{eq:bismutdoubleswitch}, the above is equivalent to
\[R^{+\phi}\pair{\bar Y_-, \bar W_+, \bar X_+, Z_+} = R^{-\phi}\pair{\bar X_+, \bar W_-, \bar Y_-, Z_-} = 0\]
for all $X_\pm, Y_\pm, Z_\pm, W_\pm \in \se{T_{\pm; 1,0} M}$. The computations can be summarized as the following result.
\begin{theorem}\label{thm:naturalgenholo}
 Let $(M, \gamma; \GG^b, \JJ)$ be a generalized K\"ahler manifold, and $\phi = \gamma + db$. Let $\bar\partial_\diamond^\phib$ be the $\bar L$-connection on $\TT^{1,0}_\JJ M$ defined by
 \[\bar\partial^\phib_{\diamond, \bar x} y := \sair{\bar x \diamond_\phib y}_{1,0}\]
 for $x, y \in \se{\TT^{1,0}_\JJ M}$. It is a $\JJ$-holomorphic structure on $\TT^{1,0}_\JJ M$ iff
 \[R^{\pm \phi}_{\bar X_\mp, \bar Y_\pm} \bar Z_\pm = 0\]
 for all $X_\pm, Y_\pm, Z_\pm \in \se{T_{\pm;1,0}M}$. In particular, if $\nabla^{\pm \phi}$ are flat on $TM$, $\bar\partial^\phib_{\diamond}$ is a $\JJ$-holomorphic structure on $\TT^{1,0}_\JJ M$.
 \qed
\end{theorem}

\subsection{Chern curvature}\label{subsec:cherncurvature}
Analogous to the classical situation, on a $\JJ$-holomorphic Hermitian bundle $\pair{V, h}$, there exists a unique \emph{generalized Chern connection}:
\begin{equation}\label{eq:genchern}
 \gcon^{h,C} := \bar\partial_\JJ + \partial_\JJ
\end{equation}
where $\partial_\JJ : \se V \to \Omega^{1,0}_\JJ\pair{V} := \se{\TT_\JJ^{0,1} M \tensor V}$ is defined by
\[d_{\bar L} h\pair{v_1, v_2} = h\pair{\bar\partial_\JJ v_1, v_2} + h\pair{v_1, \partial_\JJ v_2}\]
for all $v_j \in \se V$. 
\begin{defn}\label{defn:cherncurvature}
 On a generalized complex manifold $\pair{M, \gamma; \JJ}$, let $\gcon^\TT$ be a $TM$-torsion free generalized connection on $\TT M$. For a $\JJ$-holomorphic Hermitian bundle $\pair{V, h, \bar\partial_\JJ}$, its \emph{$\gcon^\TT$-Chern curvature} $\gcurv^{\TT, C}$ is the $\gcon^\TT$-curvature of its generalized Chern connection \eqref{eq:genchern}.
\end{defn}

When $\gcon^\TT$ is $\JJ$-compatible with $*_\gamma$, \eqref{eq:dTcomponentbar} implies that the extension \eqref{eq:derivationextension} of $\gcon^\TT$ to $\Omega_\TT^*\pair{V}$ by $\gcon^{h,C}$ is compatible with the extension \eqref{eq:extensionpartialbar} of $\bar\partial_\JJ$ to $\Omega_\JJ^{0,*}\pair{V}$.
\begin{prop}\label{prop:cherncurvaturetype}
 Let $\gcon^\TT$ be a $\gamma$-$\JJ$-connection on a generalized complex manifold $\pair{M, \gamma; \JJ}$, the corresponding $\gcon^\TT$-Chern curvature $\gcurv^{\TT, C}$ is of type $(1,1)$ with respect to $\JJ$, i.e.:
\begin{equation}\label{eq:cherncurvature}
 \gcurv^{\TT, C} \in \Omega^{1,1}_\JJ\pair{\End\pair{V}} := \se{\TT_\JJ^{1,0} M \wedge \TT_\JJ^{0,1} M \tensor \End\pair{V}}
\end{equation}
 In particular, over a generalized K\"ahler manifold $\pair{M, \gamma; \GG^b, \JJ}$, the \emph{$\pair{\phib}$-Chern curvature} $\gcurv^{\phib, C}$ defined with $\gcon^\TT = \gcon^\phib$, is of type $\pair{1,1}$ with respect to $\JJ$.
 \qed
\end{prop}

\begin{example}[c.f. Example \ref{ex:hermitianline}]\label{ex:cherncurvaturelinebundle}
 Consider a $\JJ$-holomorphic Hermitian line bundle $\pair{V, h, \bar\partial_\JJ}$ and set $\gcon = \gcon^{h,C}$, the generalized Chern connection. Choose a unitary local section $s \in \se{V}$, there is a local section $u^{0,1}_\JJ \in \Omega^{0,1}_\JJ\pair{M}$ such that
 \[\bar\partial_{\JJ} s = u^{0,1}_\JJ \tensor s\]
 Let $\gcon^\TT$ be a $\gamma$-$\JJ$-connection, then $\gcurv^{\TT, C}\pair{\gcon} = \sqrt{-1}\dd^\TT u$ where $u \in \Omega_\TT^1\pair{M}$ and
 \[\sqrt{-1} u = u^{0,1}_\JJ - \bar{u^{0,1}_\JJ}\]
 By \eqref{eq:dTcomponentbar}, the flatness of $\bar\partial_\JJ$ implies that 
 \[\bar \partial^\TT_\JJ u^{0,1} = d_{\bar L} u^{0,1}_\JJ = 0\]
 It then gives 
 \[\gcurv^{\TT, C}\pair{\gcon} = \partial^\TT_\JJ u^{0,1} - \bar\partial^\TT_\JJ \bar{u^{0,1}} \in \Omega_\JJ^{1,1}\pair{M}\] 
\end{example}

Over a generalized K\"ahler manifold, the $\pair{\phib}$-Chern connection and curvature are related to the classical Chern connections and curvatures via the $\GG^b$-eigendecomposition.
\begin{lemma}\label{lemma:chernconnectiondecompose}
 Let $\pair{M, \gamma; \GG^b, \JJ}$ be a generalized K\"ahler manifold and $V \to M$ be a $\JJ$-holomorphic vector bundle. Let $\nabla^{h,C}_\pm$ be the $\GG^b$-eigendecomposition \eqref{eq:gconeigendecomposition} of $\gcon^{h,C}$ \eqref{eq:genchern}. Then $\nabla^{h,C}_\pm$ are the Chern connections for the induced $I_\pm$-holomorphic structures \eqref{eq:inducedIpmstructure} on $V$ respectively. Furthermore, the classical Chern curvatures are components of the $\GG^b$-eigendecomposition of the $\pair{\phib}$-Chern curvature.
\end{lemma}
\begin{proof}
 The statement about the connections follows from straightforward verification, while the statement about the curvature follows from Theorem \ref{thm:phibcurvaturedecompose}.
\end{proof}

\begin{example}\label{ex:phibcurvatureHermitianLine}
 Continue from Example \ref{ex:cherncurvaturelinebundle} and let $\pair{M, \gamma; \GG^b, \JJ}$ be generalized K\"ahler. The $I_\pm$-holomorphic structures induced by $\bar\partial_\JJ$ are given locally by $\alpha_\pm \in \Omega^{0,1}_\pm\pair{M}$:
 \[u^{0,1}_\JJ = \sair{g\inverse\pair{\alpha_+} + bg\inverse\pair{\alpha_+} + \alpha_+} - \sair{g\inverse\pair{\alpha_-} + bg\inverse\pair{\alpha_-} - \alpha_-}\]
 which satisfies
 \[\bar\partial_\pm\alpha_\pm = 0 \text{ and } \pair{\nabla^{-\phi}_{\bar X_+} \alpha_-}\pair{\bar Y_-} - \pair{\nabla^{+\phi}_{\bar Y_-}\alpha_+}\pair{\bar X_+} = 0\]
 where $X_\pm \in \se{T_{\pm;1,0} M}$ and so on. The generalized Chern connection is then
 \[\gcon^{h, C} s = \pair{u^{0,1}_\JJ - \bar{u^{0,1}_\JJ}} \tensor s\]
 while the classical Chern connections $\nabla^{b,C}_\pm$ are defined locally by $\nu_\pm = \alpha_\pm - \bar \alpha_\pm$ respectively. Then Example \ref{ex:hermitianlinephib} gives the $\GG^b$-eigendecomposition of the $\varphib$-Chern curvature.
\end{example}

\subsection{$\JJ$-Hermitian-Einstein equation}\label{subsec:JHermitianEinstein}
Let $\pair{M, \gamma; \GG^b, \JJ}$ be a generalized Hermitian manifold. The analogue to the contraction by the K\"ahler form is the \emph{$\JJ_-$-contraction}.

\begin{defn}\label{defn:Jminuscontraction}
 The \emph{$\JJ_-$-contraction $\Lambda_{\JJ_-} : \wedge^2 \TT_\CC M \to \RR$} is given by:
 \begin{equation}\label{eq:Jcontraction}
  \Lambda_{\JJ_-}\pair{x \wedge y} := \<\JJ_- x, y\> =  \GG^b\pair{\JJ x, y}
 \end{equation}
 where $x, y \in \se{\TT_\CC M}$.
\end{defn}
\noindent
In terms of the $\GG^b$-eigendecomposition, it corresponds to
\begin{equation}\label{eq:Jcontractiondecompose}
 \Lambda_{\JJ_-}\pair{x^{b\pm} \wedge y^{b\pm}} = \omega_\pm\pair{X, Y} ~~~ \text{ and } ~~~ \Lambda_{\JJ_-}\pair{x^{b\pm} \wedge y^{b\mp}} = 0
\end{equation}
where $X = \pi\pair{x^{b\pm}}$, $Y = \pi\pair{y^{b\pm}}$ and $\omega_\pm = g I_\pm$.

A version of the Hermitian-Einstein equation can thus be formulated in this context.
\begin{defn}\label{defn:hermitianEinstein}
 Let $\pair{M, \gamma; \GG^b, \JJ}$ be a generalized Hermitian manifold and $\gcon^\TT$ be a $\gamma$-$\JJ$-connection on $\TT M$. A Hermitian metric $h$ on a $\JJ$-holomorphic vector bundle $V \to M$ is \emph{$\gcon^\TT$-$\JJ$-Hermitian-Einstein} if the corresponding $\gcon^\TT$-Chern curvature satisfies the following \emph{$\gcon^\TT$-$\JJ$-Hermitian-Einstein equation}:
 \begin{equation}\label{eq:hermitianEinstein}
  \sqrt{-1}\Lambda_{\JJ_-}\pair{\gcurv^{\TT, C}(V)} = 2 c \Id_V
 \end{equation}
 for some $c \in \RR$. If $\pair{M, \gamma; \GG^b, \JJ}$ is generalized K\"ahler and $\gcon^\TT = \gcon^\phib$, the equation \eqref{eq:hermitianEinstein} will be simply called the \emph{$\JJ$-Hermitian-Einstein equation}, a solution of which is called a \emph{$\JJ$-Hermitian-Einstein metric}.
\end{defn}
 
When $\pair{M, \gamma; \GG^b, \JJ}$ is generalized K\"ahler and $\gcon^\TT = \gcon^\phib$, by \eqref{eq:phibcurvaturedecompose} and \eqref{eq:Jcontractiondecompose} the $\GG^b$-eigendecomposition of the left hand side of \eqref{eq:hermitianEinstein} is given by
\begin{equation}\label{eq:Jminuscontractioneigendecompose}
 \Lambda_{\JJ_-}\pair{\gcurv^{\phib, C}(V)} = \Lambda_+\pair{F^C_+(V)} + \Lambda_-\pair{F^C_-(V)}
\end{equation}
where $F^C_\pm$ are the curvatures of the classical Chern connections $\nabla^h_\pm$ on $V$. It follows that \eqref{eq:hermitianEinstein} is equivalent to an equation first proposed by Hitchin in \cite{Hitchin:GenHoloBundles:11}.
\begin{prop}\label{prop:hermitianEinsteinHitchin}
 Over a generalized K\"ahler manifold $\pair{M, \gamma; \GG^b, \JJ}$, the equation \eqref{eq:hermitianEinstein} is equivalent to
 \begin{equation}\label{eq:hermitianEinsteinHitchin}
  \dfrac{\sqrt{-1}}{2}\pair{F_+^C\pair{V} \wedge \omega_+^{m-1} + \pair{-1}^\varepsilon F_-^C\pair{V} \wedge \omega_-^{m-1}} = c \pair{m-1}! \Id_V d\vol_g
 \end{equation}
 where $\varepsilon = 0$ if $I_\pm$ induce the same orientation on $TM$ and $\varepsilon = 1$ otherwise.
\end{prop}
\begin{proof}
 Suppose that $d\vol_g = \dfrac{1}{m !} \omega_+^m$, then 
 \[m F_+^C\pair{V} \wedge \omega_+^{m-1} = \Lambda_+\pair{F_+^C\pair{V}} \omega_+^m = m! \Lambda_+\pair{F_+^C\pair{V}} d\vol_g\]
 Note that $\omega_+^m = \pair{-1}^\varepsilon\omega_-^m$, it follows from \eqref{eq:Jminuscontractioneigendecompose} that \eqref{eq:hermitianEinsteinHitchin} is equivalent to \eqref{eq:hermitianEinstein}.
\end{proof}

The $\JJ_-$-contraction naturally provides the definition of a degree.
\begin{defn}\label{defn:degreejholomorphic}
 On a generalized Hermitian manifold $\pair{M, \gamma; \GG^b, \JJ}$ with a $\gamma$-$\JJ$-connection $\gcon^\TT$, let $\pair{V, h,  \bar\partial_\JJ}$ be a $\JJ$-holomorphic Hermitian vector bundle.
 The \emph{$\gcon^\TT$-$\GG^b$-degree} of $V$ is given by:
 \begin{equation}\label{eq:Jdegree}
  \deg_{\GG^b}^\TT\pair{V, h} := \dfrac{\sqrt{-1}}{4\pi}\int_M \tr_V \sair{\Lambda_{\JJ_-}\pair{\gcurv^{\TT, C}(V)}} d\vol_g
 \end{equation}
 When $\pair{M, \gamma; \GG^b, \JJ}$ is generalized K\"ahler and $\gcon^\TT = \gcon^\phib$, the $\gcon^\phib$-$\GG^b$-degree is simply called the \emph{$\GG^b$-degree} and denoted $\deg_{\GG^b}\pair{V, h}$.
\end{defn}

Recall that the classical degrees of $\pair{V, h}$ with the induced $I_\pm$-holomorphic structure are
\begin{equation}\label{eq:classicaldegree}
 \deg_\pm\pair{V, h} := \dfrac{\sqrt{-1}}{2\pi} \int_M \tr_h\sair{\Lambda_\pm\pair{F^C_\pm(V)}} d\vol_g
\end{equation}
\begin{theorem}\label{thm:hitchindegreeformula}
 Let $\pair{M, \gamma; \GG^b, \JJ}$ be a generalized K\"ahler manifold and $\pair{V, h, \bar\partial_\JJ}$ be a $\JJ$-holomorphic vector bundle on $M$. Then
 \begin{equation}\label{eq:hitchindegree}
  \deg_{\GG^b}\pair{V, h} = \dfrac{1}{2}\sair{\deg_+\pair{V, h} + \deg_-\pair{V, h}}
 \end{equation}
 where $\deg_{\GG^b}$ and $\deg_\pm$ are as given in \eqref{eq:Jdegree} and \eqref{eq:classicaldegree} respectively.
\end{theorem}
\begin{proof}
 It follows from \eqref{eq:Jminuscontractioneigendecompose}.
\end{proof}

\begin{example}\label{ex:degreelinebundle}
Continue from Example \ref{ex:cherncurvaturelinebundle} and work now on a generalized Hermitian manifold $\pair{M, \gamma; \GG^b, \JJ}$, with a $\gamma$-$\JJ$-connection $\gcon^\TT$. For $f \in \se{M}$, let $h_1 = e^{2f}h$ be another Hermitian metric on the line bundle $V$. An $h_1$-unitary local section is then given by $s_1 = e^{-f} s \in \se V$, which leads to
 \[\bar\partial_{\JJ} s_1 = \pair{u^{0,1}_\JJ - d_{\bar L} f} \tensor s_1\]
 Let $\gcon_1$ be the corresponding generalized Chern connection, whose $\gcon^\TT$-Chern curvature is given by
 \[\gcurv^{\TT, C}\pair{\gcon_1} = \dd^\TT \pair{\sqrt{-1}u + d_L f - d_{\bar L} f} = \gcurv^{\TT, C}\pair{\gcon} - 2\partial_\JJ^\TT d_{\bar L} f\]
 where the last step is due to \eqref{eq:vanishingcomponents}. Hence
 \begin{equation}\label{eq:degreediff}
  \deg_{\GG^b}^\TT\pair{V, h_1} - \deg_{\GG^b}^\TT\pair{V, h} = - \dfrac{\sqrt{-1}}{2\pi} \int_M \Lambda_{\JJ_-}\pair{\partial_\JJ^\TT d_{\bar L}f} d\vol_g
 \end{equation}
 It follows that $\deg_{\GG^b}^\TT$ is independent of the Hermitian metric on $V$ iff the right hand side of \eqref{eq:degreediff} vanishes for all $f \in \se{M}$.
\end{example}

 The integrand in \eqref{eq:degreediff} gives rise to the second order operator for $f \in \se{M}$:
 \[P_\JJ \pair{f}: = -\sqrt{-1} \Lambda_{\JJ_-}\pair{\partial_\JJ^\TT d_{\bar L}f}\]
 Similar to the classical case (e.g. \cite{LubkeTeleman:KobayashiHitchin:95}), $P_\JJ$ is elliptic since its symbol is given by
 \[\sigma\pair{P_\JJ} = 4 \sqrt{-1}\Lambda_{\JJ_-}\pair{\xi^{0,1}_\JJ \wedge \xi^{1,0}_\JJ} = 4 \GG\pair{\xi^{0,1}_\JJ, \xi^{1,0}_\JJ} = \norm{\xi}_g^2 \]
 where $\xi^{0,1}_\JJ$ is the projection of $\xi \in T^*M$ to $\TT_\JJ^{1,0} M$ and $\xi^{1,0}_\JJ$ is the complex conjugate.
\begin{defn}\label{defn:genGauduchon}
 For a generalized Hermitian manifold $\pair{M, \gamma; \GG^b, \JJ}$, the metric $\GG^b$ is \emph{$\gcon^\TT$-$\JJ$-Gauduchon} if the right hand side of \eqref{eq:degreediff} vanishes for all $f \in \se{M}$. When the structure is generalized K\"ahler, a $\gcon^\phib$-$\JJ$-Gauduchon the metric $\GG^b$ is simply said to be \emph{$\JJ$-Gauduchon}.
\end{defn}
\noindent
As in the classical situation, if $\pair{M, \gamma; \GG^b, \JJ}$ is $\gcon^\TT$-$\JJ$-Gauduchon and $\pair{V, h}$ solves \eqref{eq:hermitianEinstein}, then the constant $c$ is given by
 \[c = \dfrac{2\pi \deg_{\GG^b}^\TT\pair{V}}{\rank\pair{V} \Vol_g\pair{M}}\]
It is then natural to extend the notions of slope and stability to $\JJ$-holomorphic vector bundles over a $\gcon^\TT$-$\JJ$-Gauduchon generalized Hermitian manifold. The notion of coherent subsheaf in this context can be adopted from Definition 3.4 of \cite{HuMoraruSeyyidali:KobayashiHitchin:16}.
\begin{defn}\label{defn:stability}
 Let $\GG^b$ be a $\gcon^\TT$-$\JJ$-Gauduchon metric. The \emph{$\gcon^\TT$-$\GG^b$-slope} of a $\JJ$-holomorphic vector bundle $\pair{V, \bar\partial_\JJ}$ over $M$ is
 \begin{equation}\label{eq:gslope}
  \mu^\TT_{\GG^b}\pair{V} := \dfrac{\deg^\TT_{\GG^b}\pair{V}}{\rank\pair{V}}
 \end{equation}
 The bundle $V$ is \emph{$\gcon^\TT$-$\GG^b$-semistable} if for any coherent $\JJ$-holomorphic subsheaf $W$ of $V$:
 \begin{equation}\label{eq:gsemistable}
  \mu^\TT_{\GG^b}\pair{W} \varleq \mu^\TT_{\GG^b}\pair{V}
 \end{equation}
 $V$ is said to be \emph{$\gcon^\TT$-$\GG^b$-stable} if strict inequality holds in \eqref{eq:gsemistable}. Over a $\JJ$-Gauduchon generalized K\"ahler manifold, the corresponding notions are simply refered to as the \emph{$\GG^b$-slope}, \emph{$\GG^b$-semistable} and \emph{$\GG^b$-semistable} respectively.
\end{defn}

Recall that over a Hermitian manifold $\pair{M, g, I}$, the degree of any holomorphic vector bundle $V$ is independent of a Hermitian metric on $V$ iff $g$ is Gauduchon, i.e.
\begin{equation}\label{eq:classicalGauduchon}
 \partial\bar\partial\pair{\omega^{m-1}} = 0
\end{equation}
where $m$ is the complex dimension of $M$. On a generalized K\"ahler manifold, the $\JJ$-Gauduchon condition can be expressed in a similar fashion.

\begin{prop}\label{prop:genGauduchon}
 A generalized K\"ahler manifold $\pair{M, \gamma; \GG^b, \JJ}$ is $\JJ$-Gauduchon iff
 \begin{equation}\label{eq:genGauduchon}
  \partial_+\bar\partial_+\pair{\omega_+^{m-1}} + \pair{-1}^\varepsilon \partial_-\bar\partial_-\pair{\omega_-^{m-1}} = 0
 \end{equation}
 where $\varepsilon = 0$ if $I_\pm$ induce the same orientation on $TM$ and $\varepsilon = 1$ otherwise.
\end{prop}
\begin{proof}
 Note that the degree of a Hermitian vector bundle coincides with that of its determinant line bundle, thus it's sufficient to consider line bundles.
 By \eqref{eq:ddbarf} and \eqref{eq:Jcontractiondecompose}
 \begin{equation*}
  \int_M \Lambda_{\JJ_-}\pair{\partial_\JJ^\phib d_{\bar L}f} d\vol_g =  \int_M\sair{\Lambda_+\pair{\partial_+\bar\partial_+f} + \Lambda_-\pair{\partial_-\bar\partial_-f}} d\vol_g
 \end{equation*}
 The statement then follows from integration by parts.
\end{proof}

\begin{remark}\label{remark:genGauduchon}
 Evidently, \eqref{eq:genGauduchon} holds for $m < 3$, since the $3$-form below is closed:
 \[\phi = \gamma + db = \mp d^c_\pm \omega_\pm\]
 For $m \vargeq 3$, \eqref{eq:genGauduchon} is equivalent to
 \[\phi \wedge d\pair{\omega_+^{m-2} - \pair{-1}^\varepsilon \omega_-^{m-2}} = 0\]
 In particular, $\pair{M, \gamma; \GG^b, \JJ}$ is $\JJ$-Gauduchon if the difference $\omega_+^{m-2} - \pair{-1}^\varepsilon\omega_-^{m-2}$ is closed. When $m = 3$, \eqref{eq:genGauduchon} can also be rewritten as
 \[\phi_+^{\pair{2,1}} \wedge \phi_+^{\pair{{1,2}}} + \pair{-1}^\varepsilon\phi_-^{\pair{2,1}} \wedge \phi_-^{\pair{{1,2}}} = 0\]
 where, for instance, $\phi_\pm^{\pair{2,1}}$ denote the $\pair{2,1}$-component of $\phi$ with respect to $I_\pm$ respectively. It is clear that the generalized K\"ahler manifold is $\JJ$-Gauduchon if $g$ is Gauduchon with respect to both $I_\pm$.
\end{remark}

\section{Geometric Lax flows}\label{sec:laxflows}

A \emph{Lax pair} \cite{Lax:laxpair:68} consists of two families of operators $\bair{\pair{L_t, P_t}: t \in I \subseteq \RR}$ such that
\begin{equation}\label{eq:laxpair}
 \dfrac{d}{dt}P_t = \sair{L_t, P_t}
\end{equation}
where $\bair{L_t}$ is the \emph{Lax operator} and it is assumed that $0 \in I$. The equation \eqref{eq:laxpair} is also said to be in \emph{the Lax form}. Suppose that $\bair{\Psi_t}$ is generated by $\bair{L_t}$, i.e. it solves the equation
\begin{equation}\label{eq:laxgenerated}
 \dfrac{d}{dt}\Psi_t = L_t\Psi_t, \text{ with } \Psi_0 = \Id
\end{equation}
then $\bair{P_t}$ can be obtained from pushing-forward of an initial operator $P_0$ by $\bair{\Psi_t}$:
\[P_t = \Psi_t P_0 \Psi_t\inverse\]
In particular, $\bair{P_t}$ is an isospectral family.

Let $\bair{A_t}$ be a smooth family of operators on the same space, then the $t$-differential of $\bair{\Psi_t\inverse A_t\Psi_t}$ describes the extend to which $\bair{\pair{L_t, A_t}}$ fails to be a Lax pair as well.

\begin{defn}\label{defn:laxdiff}
 Suppose that $\bair{\pair{L_t, P_t}}$ is a Lax pair and $\bair{A_t}$ be a smooth family of operators on the same space. The \emph{$L_t$-differential of $A_t$ (along the Lax flow)} is:
 \begin{equation}\label{eq:laxriccidiff}
  \delta_{L} A_t := \frac{d}{dt} A_t - \sair{L_t, A_t}
 \end{equation}
\end{defn}
It is straightforward to verify that commutativity with $P_t$ is preserved by $\delta_L$.
\begin{lemma}\label{lemma:commutinglaxdifferential}
 If $A_t$ and $P_t$ commute for all $t$, then $\delta_L A_t$ also commute with $P_t$ for all $t$.
\end{lemma}
\begin{proof}
 Notice that
 \[\dfrac{d}{dt}\pair{\Psi_t\inverse A_t\Psi_t} = \Psi_t\inverse \pair{\delta_L A_t} \Psi_t\]
 from which the statement follows.
\end{proof}

When a pair of geometric quantities forms a Lax pair, the corresponding equation \eqref{eq:laxpair} is said to generate a \emph{geometric Lax flow}. Two main classes of examples will be described, where the operators $\bair{P_t}$ are either generalized metrics or generalized almost complex structures. Such Lax pairs impose certain necessary conditions on the Lax operator.

\begin{defn}\label{defn:laxcompatible}
 Let $\GG^b$ be a generalized metric and $\JJ$ a generalized (almost) complex structure. An operator $\LL$ on $\se{\TT M}$ is \emph{Lax compatible with $\GG^b$} if
 \begin{equation}\label{eq:gmetricdeformation}
  \aair{\LL x^{b+}, y^{b-}} + \aair{x^{b+}, \LL y^{b-}} = 0
 \end{equation}
 for all $x^{b+} \in \se{C_+^{b}}$ and $y^{b-} \in \se{C_-^{b}}$.
 The operator $\LL$ is \emph{Lax compatible with $\JJ$} if
  \begin{equation}\label{eq:gcplexdeformation}
   \aair{\LL x, y} + \aair{x, \LL y} = \aair{\LL \bar x, \bar y} + \aair{\bar x, \LL \bar y} = 0
  \end{equation}
  for all $x, y \in \se{\TT^{1,0}_{\JJ} M}$.
\end{defn}

\begin{lemma}\label{lemma:laxpairconditions}
 Suppose that $\bair{\pair{\LL_t, \PP_t}}$ is a Lax pair of operators on $\se{\TT M}$. If $\PP_t$ is a smooth family of generalized metrics or almost complex structures, then $\LL_t$ is Lax compatible with $\PP_t$ for all $t$.
\end{lemma}
\begin{proof}
 It follow from the orthogonality of $\PP_t$ with respect to the pairing $\aair{,}$ and that $\PP_t^2$ are constant operators.
\end{proof}

The following simplified criteria are useful in practice.
\begin{corollary}\label{coro:laxoperatorskew}
 Let $\LL \in \se{\End\pair{\TT M}}$ and $x, y \in \TT M$. Then 
 \begin{enumerate}
  \item
  \eqref{eq:gmetricdeformation} holds for $\LL$ if the bilinear forms $\GG^{b}\pair{\LL x, y}$ are symmetric, i.e.
  \[\GG^{b}\pair{\LL x, y} = \GG^{b}\pair{\LL y, x}\]
  \item
  \eqref{eq:gmetricdeformation} and \eqref{eq:gcplexdeformation} hold for $\LL$ if the bilinear forms $\aair{\LL x, y}$ are skew-symmetric, i.e.
  \[\aair{\LL x, y} + \aair{\LL y, x} = 0\]
 \end{enumerate}
\end{corollary}
\begin{proof}
 Left for the reader.
\end{proof}

\subsection{$\theta \in \Omega_\TT^2\pair{M}$ as the Lax operator}\label{subsec:geometriclaxmorphism}

Let $\bair{P_t \in \se{T^*M^{\tensor 2}}}$ be a smooth family of $2$-tensors on $M$, whose symmetric and skew-symmetric parts are respectively $P^s_t$ and $P^a_t$, i.e.
\[P_t = P_t^s + P_t^a \text{ with } P_t^s\pair{X, Y} = P_t^s\pair{Y, X}, P_t^a\pair{X, Y} = -P_t^a\pair{Y, X}\]
for $X, Y \in \se{TM}$. Then $\bair{P_t}$ defines an initial value problem for a family of generalized metrics $\GG^{b_t}_t$ as follows:
\begin{equation}\label{eq:2tensorflow}
 \begin{cases}
  \dfrac{d}{dt} g_t & = - P_t^s\\
  \dfrac{d}{dt} b_t & = P_t^a
 \end{cases}
\end{equation}

The system \eqref{eq:2tensorflow} can be reformulated into a Lax flow for $\GG^{b_t}_t$.
\begin{lemma}\label{lemma:2tensorflow}
 Let $\bair{P_t \in \se{T^*M^{\tensor 2}}}$ and $\bair{\theta_{t} \in \Omega_\TT^2\pair{M}}$ be a smooth family of $\TT M$-forms such that
 \begin{equation}\label{eq:2tensortoform}
  \theta_{t}\pair{x^{b_t - }, y^{b_t+}} = P_t(X, Y)
 \end{equation}
 where $x^{b_t \pm} \in \se{C_\pm^{b_t}}$ with $X = \pi\pair{x^{b_t \pm}}$ and so on. Let $\theta_{t}: \TT M \to \TT M$ be given by
 \[2\aair{\theta_t\pair{x}, y} := \theta_t\pair{x, y}\]
 Then \eqref{eq:2tensorflow} is equivalent to the \emph{$2$-tensor Lax flow}
 \begin{equation}\label{eq:2formlax}
  \dfrac{d}{dt} \GG^{b_t}_t = \sair{\theta_t, \GG^{b_t}_t}
 \end{equation}
\end{lemma}
\begin{proof}
 Note that $\theta_t$ satisfies Corollary \ref{coro:laxoperatorskew} $(2)$. Since $\pair{\GG^b}^2 = \1$, the differential of a smooth family of generalized metrics is skew with respect to the generalized metric at each time. Obviously the left hand side of \eqref{eq:2formlax} is skew with respect to $\GG^{b_t}_t$. Thus only the mixed $\GG^{b_t}_t$-eigencomponents of the left hand side are non-trivial, one of which goes as follows:
 \begin{equation*}
  \begin{split}
   & 0 = \dfrac{d}{dt}\<\GG_t^{b_t} x^{b_t-}, y^{b_t+}\> \\
   = & \<\dfrac{d}{dt} \GG_t^{b_t} x^{b_t-}, y^{b_t+}\> + \<\GG_t^{b_t} \dfrac{d}{dt}\pair{b_t - g_t} X, y^{b_t+}\> + \<\GG_t^{b_t} x^{b_t-}, \dfrac{d}{dt}\pair{b_t + g_t} Y\> \\
   = & \<\dfrac{d}{dt} \GG_t^{b_t} x^{b_t-}, y^{b_t+}\> + \dfrac{d}{dt}\pair{b_t - g_t} \pair{X, Y}
  \end{split}
 \end{equation*}
 The right hand side is given by
 \begin{equation*}
  \begin{split}
   & \aair{\sair{\theta_t, \GG^{b_t}_t} x^{b_t-}, y^{b_t+}} = \aair{\pair{\theta_t\GG^{b_t}_t - \GG^{b_t}_t\theta_t} x^{b_t-}, y^{b_t+}} = -2\aair{\theta_t\pair{x^{b_t-}}, y^{b_t+}} = -P_t\pair{X, Y}
  \end{split}
 \end{equation*}
 Thus \eqref{eq:2formlax} gives rise to
 \[\dfrac{d}{dt}\pair{b_t - g_t} \pair{X, Y} = P_t\pair{X, Y}\]
 from which \eqref{eq:2tensorflow} follows. The other direction is left for the reader.
\end{proof}

A smooth conformal family of metrics $\bair{g_t}$ can be seen as a solution to \eqref{eq:2tensorflow}, by setting $P_t = P_t^s = f_t g_0$, where $\bair{f_t \in \se{M}}$ is a smooth family of functions. In this case, $b_t \equiv b_0$ is constant throughout the flow. The corresponding $\TT M$-forms as given in Lemma \ref{lemma:2tensorflow} can be chosen to be dependent only on the generalized metric $\GG^b$.

\begin{defn}\label{defn:conformalform}
 Let $\GG^b$ be a generalized metric defined by $\pair{g, b}$ where $g$ is a Riemannian metric. A $\TT M$-form $\theta \in \Omega_\TT^2\pair{M}$ is \emph{$\GG^b$-conformal} if there are $r, s \in \se{M}$ such that
 \[\theta\pair{x^{b-}, y^{b+}} = r g\pair{X, Y} + s b\pair{X, Y}\text{ and } \theta\pair{x^{b\pm}, y^{b\pm}} = 0\]
 where $x^{b \pm} \in \se{C_\pm^{b}}$ with $X = \pi\pair{x^{b \pm}}$ and so on. $r, s$ are called the \emph{conformal weights}. 
\end{defn}

\begin{remark}\label{remark:conformalform}
 When the conformal weights coincide, i.e. $r = s$, a family of $\GG^b$-conformal forms generates conformal deformations of the generalized metric $\GG^b$. Otherwise, the metric $g$ and $b$ are deformed by different factors. In particular, when the conformal weight $s = 0$, it corresponds to classical conformal deformation of the metric $g$.
\end{remark}

\subsection{$\varphib$-curvature Lax flow}\label{subsec:phibcurvatureflow}
A special case of \eqref{eq:2formlax} is when $\theta_t$ are $\dd^\phib$-exact, i.e. $\theta_t = \dd^{\phibt} u_t$ for a smooth family of sections $\bair{u_t \in \Omega^1_\TT\pair{M}}$ and $3$-forms $\bair{\phi_t \in \Omega^3\pair{M}}$. Suppose that $u_t = Z_t + \zeta_t$ then
\begin{equation*}
 \begin{split}
  & \pair{\dd^{\phibt} u_t} \pair{x^{b_t -}, y^{b_t +}} = X u_t\pair{y^{b_t +}} - Y u_t\pair{x^{b_t -}} - u_t\pair{x^{b_t -} \diamond_{\phibt} y^{b_t +}}\\
  = & \pair{\LLC_{Z_t} g_t} \pair{X, Y} + \sair{d\pair{\zeta_t - \iota_{Z_t} b_t} - \iota_{Z_t}\phi_t}\pair{X, Y}
 \end{split}
\end{equation*}
It follows that in this case, the Lax flow \eqref{eq:2formlax} is equivalent to
\begin{equation}\label{eq:generalizedpushforward}
 \begin{cases}
  \dfrac{d}{dt} g_t & = - \LLC_{Z_t}g_t\\
  \dfrac{d}{dt} b_t & = -\LLC_{Z_t} b_t + d\zeta_t - \iota_{Z_t}\pair{\phi_t - db_t} 
 \end{cases}
\end{equation}
When $\gamma = \phi_t - db_t$ is a fixed closed $3$-form, \eqref{eq:generalizedpushforward} describes the push-forward of an initial generalized metric by the family of generalized diffeomorphism of $\TT M$ generated by $u_t$.
\begin{prop}\label{prop:genmorphismaslaxflow}
 Let $\gamma \in \Omega^3\pair{M}$ be a closed $3$-form and $u_t = Z_t + \zeta_t \in \se{\TT M}$ be a smooth family of sections. Let $\pair{\lambda_t, \beta_t}$ be the family of generalized diffeomorphisms generated by $\bair{u_t}$ under $*_\gamma$ (recalled below). Set $\phi_t = \gamma + db_t$, then \eqref{eq:generalizedpushforward} coincides with the infinitesimal action by $\pair{\lambda_t, \beta_t}$ on a generalized metric $\GG^b$ via push-forward.
\end{prop}
\begin{proof}
 It's more straightforward to work with pull-back action on generalized metrics. Recall \cite{HuUribe:gensymmetry:07} that in $\pair{\lambda_t, \beta_t}$, $\lambda_t$ is the $1$-parameter family of diffeomorphisms generated by $Z_t$, and 
 \[\beta_t := \int_0^t \lambda_s^* \pair{d\zeta_s - \iota_{Z_s} \gamma} ds\]
 The push-forward of $x = X + \xi \in \TT M$ by $\pair{\lambda_t, \beta_t}$ is given by
 \[\pair{\lambda_t, \beta_t}_* x = \lambda_{t*}\pair{x + \iota_X \beta_t} = \lambda_{t*} X + \pair{\lambda_t\inverse}^*\pair{\xi + \iota_X \beta_t}\]
 with the corresponding infinitesimal action $x \mapsto  - u_t *_\gamma x$. 
 
 The pull-back of $\GG^b$ by $\pair{\lambda_t, \beta_t}$ gives the family of generalized metrics
 \[\GG^{b_t}_t\pair{x, y} := \GG^b\pair{\pair{\lambda_t, \beta_t}_* x, \pair{\lambda_t, \beta_t}_* y}\]
 where $x, y \in \se{\TT M}$ are independent of $t$. Analogous to the computations in the classical case, differentiate the above with respect to $t$ gives by the left hand side
 \begin{equation*}
  \begin{split}
   & \pair{\dfrac{d}{dt}\GG^{b_t}_t}\pair{x, y} = \dfrac{d}{dt}\sair{g_t\pair{X, Y} + g_t\inverse\pair{\xi - \iota_X b_t, \eta - \iota_Y b_t}}\\
   = & \half \sair{\pair{\dfrac{d}{dt} g_t}\pair{X, Y} - \pair{\dfrac{d}{dt} g_t}\pair{g_t\inverse\xi'_t, g_t\inverse\eta'_t} - \pair{\dfrac{d}{dt} b_t}\pair{Y, g_t\inverse \xi'_t} - \pair{\dfrac{d}{dt} b_t}\pair{X, g_t\inverse \eta'_t}}
  \end{split}
 \end{equation*}
 where $\xi'_t = \xi - \iota_X b_t$ and $\eta'_t = \eta - \iota_Y b_t$; while by the right hand side
 \begin{equation*}
  \begin{split}
   & \pair{\dfrac{d}{dt}\GG^{b_t}_t}\pair{x, y} = Z_t \GG^{b_t}_t\pair{x, y} + \GG^{b_t}_t\pair{-u_t *_\gamma x, y} +  \GG^{b_t}_t\pair{x, -u_t *_\gamma y} \\
   = & \half \sair{\pair{\LLC_{Z_t} g_t}\pair{X, Y} - \pair{\LLC_{Z_t} g_t}\pair{g_t\inverse \xi'_t, g_t\inverse \eta'_t}}\\
   & + \half \sair{\pair{-\LLC_{Z_t} b + d\zeta_t - \iota_{Z_t} \gamma}\pair{Y, g_t\inverse\xi'_t} + \pair{-\LLC_{Z_t} b + d\zeta_t - \iota_{Z_t} \gamma}\pair{X, g_t\inverse\eta'_t}}
  \end{split}
 \end{equation*}
 
 Compare these two ways of computing $\pair{\dfrac{d}{dt}\GG^{b_t}_t}\pair{x, y}$, the pull-back action on $\GG^b$ gives
 \[
 \begin{cases}
  \dfrac{d}{dt} g_t & = \LLC_{Z_t} g_t\\
  \dfrac{d}{dt} b_t & = \LLC_{Z_t} b_t - \pair{d\xi_t - \iota_{Z_t}\gamma}
 \end{cases}
 \]
 The push-forward action reverses the signs on the right hand side, which is \eqref{eq:generalizedpushforward}.
\end{proof}

A $\dd^\phib$-exact $\TT M$-form in $\Omega_\TT^2\pair{M}$ can be seen as the $\varphib$-curvature of a unitary generalized connection on the trivial Hermitian line bundle. In general, let $\pair{V, h}$ be a Hermitian line bundle with a family of unitary generalized connections $\bair{\gcon_t}$. By Example \ref{ex:hermitianlinephib}, the Lax flow \eqref{eq:2formlax} defined by $\bair{\theta_t = \sqrt{-1}\gcurv^{\phibt}\pair{\gcon_t}}$ where $\phi_t = \gamma + db_t$ is equivalent to the system
\begin{equation}\label{eq:generalizedpushforwardnonexact}
 \begin{cases}
  \dfrac{d}{dt} g_t & = - \LLC_{g_t\inverse \psi_t} g_t\\
  \dfrac{d}{dt} b_t & = \sqrt{-1}F^{b_t}_0 - \iota_{g_t\inverse \psi_t}\phi_t
 \end{cases}
\end{equation}
which reduces to \eqref{eq:generalizedpushforward} when $V$ is trivial, and admits similar interpretation in terms of (not necessarily exact) generalized diffeomorphisms.
\begin{theorem}\label{thm:linebundlemorphism}
 Given a family of unitary generalized connections $\bair{\gcon_t}$ on a line bundle $V$ with Hermitian metrics $\bair{h_t}$, the Lax flow \eqref{eq:2formlax} defined by $\bair{\theta_t = \gcurv^{\phibt}\pair{\gcon_t}}$, where $\phi_t = \gamma + db_t$, corresponds to the push-forward of an initial generalized metric $\GG^{b}$ by a family of generalized diffeomorphisms, which may not be \emph{exact}, i.e. is not generated by global sections of ${\TT M}$. 
\end{theorem}
\begin{proof}
 In the notations of Example \ref{ex:hermitianlinephib}, the local section defining $\gcon_t$ is given by $\sqrt{-1}u_t = \sqrt{-1}\pair{Z_t + \zeta_t}$, where $Z_t$ is a global vector field and $\zeta_t$ is a locally defined $1$-form, with
 \[Z_t = g_t\inverse \psi_t = \half g_t\inverse \pair{\nu_{t,+} - \nu_{t,-}} \text{ and } \zeta_t = \iota_{Z_t} b_t + \half \pair{\nu_{t,+} + \nu_{t,-}}\]
 Then $\nabla^{b_t}_0$ is defined by the local section $\half \pair{\nu_{t,+} + \nu_{t,-}}$, which implies that
 \[\sqrt{-1}F^{b_t}_0 = d \pair{\zeta_t - \iota_{Z_t} b_t} = - \LLC_{Z_t} b_t + d\zeta_t + \iota_{Z_t} db_t\]
 Locally, the second equation in \eqref{eq:generalizedpushforwardnonexact} is thus
 \[ \dfrac{d}{dt} b_t = - \LLC_{Z_t} b_t + d\zeta_t - \iota_{Z_t}\gamma\]
 noting that $\phi_t = \gamma + db_t$. By Proposition \ref{prop:genmorphismaslaxflow}, \eqref{eq:generalizedpushforwardnonexact} corresponds to pushing-forward by the local exact generalized diffeomorphisms generated by $\bair{\pair{Z_t, \zeta_t}}$. Globally, \eqref{eq:generalizedpushforwardnonexact} corresponds to pushing-forward by possibly non-exact generalized diffeomorphisms.
\end{proof}

\begin{example}\label{ex:hermitianlineflow}
 Let $\pair{M, \gamma; \JJ}$ be a generalized complex manifold and $\bair{\theta_t \in \Omega_\TT^2\pair{M}}$ be a smooth family of $\TT M$-forms. The Lax flow with initial value $\JJ_0 = \JJ$
 \[\dfrac{d}{dt}\JJ_t = \sair{\theta_t, \JJ_t}\]
 consists of generalized almost complex structures. The flow above preserves $\JJ$ iff $\sair{\theta_t, \JJ} = 0$ for all $t$, which is equivalent to $\theta_t \in \Omega_\JJ^{1,1}\pair{M}$. Thus, starting with a generalized Hermitian manifold $\pair{M, \gamma; \GG^b, \JJ}$, the Lax flow \eqref{eq:2formlax} defined by $\bair{\theta_t \in \Omega_\JJ^{1,1}\pair{M}}$ produces a family of generalized Hermitian structures with the same $\JJ$.
 
 Suppose furthermore that $\pair{M, \gamma; \GG^b, \JJ}$ admits a $\gamma$-$\JJ$-connection $\gcon^\TT$. Fix a $\JJ$-holomorphic line bundle $V$ and let $\bair{\theta_t}$ be the $\gcon^\TT$-Chern curvatures for a smooth family of Hermitian metrics $\bair{h_t}$ on $V$. By Proposition \ref{prop:cherncurvaturetype}, $\bair{\theta_t}$ consists of $\pair{1,1}$-forms with respect to $\JJ$. By Theorem \ref{thm:linebundlemorphism}, the corresponding Lax flow \eqref{eq:2formlax} corresponds to the push-forward of $\GG^b$ by a family of generalized diffeomorphisms $\bair{\pair{\lambda_t, \beta_t}}$. Alternatively, it can be seen as a family of generalized Hermitian structures $\pair{M, \gamma; \GG^b, \pair{\lambda_t, \beta_t}^* \JJ}$ with $\GG^b$ fixed.
\end{example}

\subsection{Ricci Lax flow} \label{subsec:ricciflow}

When the Ricci curvature of a generalized connection $\gcon^\TT$ satisfies \eqref{eq:gmetricdeformation}, it can serve as the Lax operator in a Lax pair involving the generalized metric.
\begin{defn}\label{defn:laxricciflow}
 A smooth family of pairs $\bair{\pair{\GG_t^{b_t}, \gcon^{\TT_t}}}$ of generalized metrics $\GG_t^{b_t}$ and $\GG_t^{b_t}$-metric connections is a solution to the \emph{$\gcon^\TT$-Ricci Lax flow} if $\bair{\Ric^{\TT_t}}$ satisfies \eqref{eq:gmetricdeformation} and the pair $\bair{\pair{\Ric^{\TT_t}, \GG_t^{b_t}}}$ form a Lax pair, i.e.
\begin{equation}\label{eq:laxricci}
 \frac{d}{dt} \GG^{b_t}_t = \sair{\Ric^{\TT_t}, \GG^{b_t}_t}
\end{equation}
When $\gcon^{\TT_t}$ is prescribed to depend on $\GG_t^{b_t}$, \eqref{eq:laxricci} becomes an equation for $\GG_t^{b_t}$ only, in which case the family of generalized metrics $\GG_t^{b_t}$ is said to be \emph{a solution} to \eqref{eq:laxricci}. If $\gcon^{\TT_t} = \gcon^{\phibt}$, then the flow \eqref{eq:laxricci} is simply called the \emph{Ricci Lax flow}.
\end{defn}

Since $\Rc^\phib$ is symmetric, it satisfies the condition in Corollary \ref{coro:laxoperatorskew}. 
Theorem \ref{thm:genricciflow} below shows that the Ricci Lax flow is equivalent to
\begin{equation}\label{eq:riemannlaxriccidecomposition}
 \dfrac{d}{dt}\pair{g_t \mp b_t} = -2 Rc_t^{\pm \phi_t} = -2 Rc_t + \dfrac{1}{2}\phi_t^2 \pm d^*\phi_t
\end{equation}
Let $\gamma \in \Omega^3\pair{M}$ such that $d\gamma = 0$. The following \emph{generalized Ricci flow} \cite{GarciaFernandezStreets:generalizedRicciflow:20, StreetsTian:GKpluriclosedflow:12} is a system for Riemannian metrics $g_t$ and $b_t \in \Omega^2\pair{M}$:
\begin{equation}\label{eq:genricciflow}
 \begin{cases}
  \dfrac{d}{dt} g_t & = -2 Rc_{t} + \dfrac{1}{2}\pair{\gamma + db_t}^2\\
  \dfrac{d}{dt} b_t & = - d^*\pair{\gamma + db_t}
 \end{cases}
\end{equation}
then \eqref{eq:riemannlaxriccidecomposition} and \eqref{eq:genricciflow} coincide when $\phi_t = \gamma + db_t$.

\begin{theorem}\label{thm:genricciflow}
 Under the further constraints that $\gcon^\TT = \gcon^\phib$ and is metric compatible with the Dorfman bracket $*_\gamma$ \eqref{eq:diamondmetriccompatible}, i.e. $\phi_t - db_t = \gamma \in \Omega^3(M)$ for all $t$, 
 the Ricci Lax flow \eqref{eq:laxricci} is equivalent to the generalized Ricci flow \eqref{eq:genricciflow}.
\end{theorem}
\begin{proof}
 It's obvious that \eqref{eq:genricciflow} follows from \eqref{eq:riemannlaxriccidecomposition}, by matching the symmetric and skew-symmetric terms on both sides. To obtain \eqref{eq:riemannlaxriccidecomposition} from \eqref{eq:laxricci}, fix $X, Y \in TM$ and consider 
 \[x^{b_t -} = X + \pair{b_t - g_t}X \text{ and } y^{b_t +} = Y + \pair{b_t + g_t}Y \]
 The left hand side of \eqref{eq:riemannlaxriccidecomposition} is computed in Lemma \ref{lemma:2tensorflow}, while the right hand side becomes:
 \begin{equation*}
  \begin{split}
   & \< \sair{\Ric^{\phibt}, \GG^{b_t}_t} x^{b_t-}, y^{b_t+}\> = \< \pair{\Ric^{\phibt}\GG^{b_t}_t - \GG^{b_t}_t\Ric^{\phibt}} x^{b_t-}, y^{b_t+}\>\\
   = & -2\Rc^{\phibt}\pair{x^{b_t-}, y^{b_t+}} = -2Rc^{+\phi_t}\pair{X, Y}
  \end{split}
 \end{equation*}
 This gives half of \eqref{eq:riemannlaxriccidecomposition}. The other half is equivalent.
\end{proof}

\begin{example}\label{ex:genriccieigen}
 Consider two $3$-dimensional Lie groups: the $3$-torus $T = \pair{U\pair{1}}^3$ and $G = SU\pair{2}$, with their invariant metrics $g$ and invariant volume forms $\phi$. For both, set $b = 0$ and consider the corresponding $\pair{\phi, 0}$-Ricci curvature $\Ric^{\phi, 0}$.
 
 For $T$, the invariant metric is flat, thus $Ric = 0$ while $Ric^{\pm\phi} \neq 0$ by \eqref{eq:torsionricci}.
 In the $\GG^b$-eigendecomposition, $\Rc^{\phi, 0}$ is of the form
 \[\Rc^{\phi,0} = \begin{pmatrix}0 & R\\R^T & 0\end{pmatrix}\]
 where $R = \1_3$, the $3\times 3$ identity matrix.

 For $G$, the invariant metric is the standard round metric on $S^3$ which is not flat, while the connections with torsion $\pm\phi$ are flat, hence $Ric^{\pm\phi} = 0$. In this case, $\Rc^{\phi, 0}$ is of the form
 \[\Rc^{\phi, 0} = \begin{pmatrix}R & 0\\0 & R\end{pmatrix}\]
 where $R$ here is the classical Ricci tensor for the round metric (c.f. \S \ref{subsec:compactLiegroups}).
\end{example}

Definition \ref{defn:conformalform} and Theorem \ref{thm:genricciflow} lead to the natural generalization of Ricci soliton.
\begin{defn}\label{defn:laxriccisoliton}
 Let $\gamma \in \Omega^3\pair{M}$ be a closed $3$-form. A smooth family of generalized metrics $\bair{\GG^{b_t}_t}$ is a \emph{$\gamma$-Ricci Lax soliton} if there exists a smooth family of sections $\bair{u_t \in \se{\TT M}}$ and $\GG^{b_t}_t$-conformal forms $\bair{\theta_t \in \Omega_\TT^2\pair{M}}$ with constant conformal weights $r_t$ and $s_t$, and
 \begin{equation}\label{eq:laxriccisoliton}
  \sair{\Ric^\phibt - \dd^\phibt u_t - \theta_t, \GG^{b_t}_t} = 0
 \end{equation}
 where $\phi_t = \gamma + db_t$. The family is a \emph{gradient $\gamma$-Ricci Lax soliton} if furthermore there is $f \in \se{M}$ such that $\bair{u_t = - \half \GG^{b_t}_t \pair{d f}}$.
\end{defn}

Notice that $\GG_t^{b_t}\pair{df} = \nabla^t f + \iota_{\nabla^t f} b_t$, where $\nabla^t f$ is the gradient of $f$ with respect to $g_t$. When $s_t = 0$, the gradient $\gamma$-Ricci Lax soliton equation is then equivalent to the system
\begin{equation}\label{eq:gradientlaxriccisoliton}
 \begin{cases}
  Rc_t - \dfrac{1}{4}\phi_t^2 + \pair{\nabla^t}^2 f & = r_t g_t \\
  d^*\phi_t + \iota_{\nabla^t f} \phi_t & = 0
 \end{cases}
\end{equation}
 where $\pair{\nabla^t}^2 f$ is the Hessian of $f$ with respect to the Levi-Civita connection of $g_t$. When $r_t \equiv r$ is a constant function independent of $t$, the system \eqref{eq:gradientlaxriccisoliton} is exactly the \emph{generalized Ricci soliton} equation (Apostolov-Streets-Ustinovskiy \cite{ApostolovStreetsUstinovskiy:KahlerRicciSoliton:21} and references therein).

\subsection{Bismut-Ricci Lax flow}\label{subsec:bismutlaxricci}
Even though $\Ric^\bismut$ is neither symmetric nor skew-symmetric, it satisfies \eqref{eq:gmetricdeformation}, and thus can be used to define a Lax flow for generalized metrics. 
In \eqref{eq:laxricci}, taking $\Ric^\bismut$ as the Lax operator leads to the \emph{Bismut-Ricci Lax flow}:
\begin{equation}\label{eq:bismutlaxricci}
 \dfrac{d}{dt} \GG_t^{b_t} = \sair{\Ric_t^{\bismutt}, \GG_t^{b_t}}
\end{equation}
Since $\Rc^\bismut$ has the same mixed components as $\Rc^\phib$, i.e.
\[\sair{\Ric^{\phib}, \GG^{b}} = \sair{\Ric^{\bismut}, \GG^{b}}\]
by Theorem \ref{thm:genricciflow}, \eqref{eq:bismutlaxricci} is equivalent to \eqref{eq:riemannlaxriccidecomposition}, so they generate the same flow for the generalized metric \cite{Streets:genTdualityRenorm:13}. The $\Ric^\bismut$-differential of $\gcon^{\bismutt}$ takes a particularly simple form.
\begin{prop}\label{prop:riccidifferentialbismut}
 Fix $X, Y, Z \in \se{TM}$ and $x \in \se{\TT M}$, with $\pi\pair{x} = X$. Let $y^{b_t\pm}, z^{b_t\pm} \in C^{b_t}_\pm$ so that $\pi\pair{y^{b_t\pm}} = Y$ and $\pi\pair{z^{b_t\pm}} = Z$. Then
\begin{equation}\label{eq:riccidifferentialbismut}
 \GG^{b_t}_t\pair{\pair{\delta_{\Ric^\bismut}\gcon^{\bismutt}_{x}} y^{b_t\pm}, z^{b_t\pm}} = - g_t\pair{Y, \dfrac{d}{dt} \nabla^{\pm{\phi_t}}_X Z}
\end{equation}
\end{prop}
\begin{proof}
 Only the case for $C^{b_t}_+$ is shown here and the case for $C^{b_t}_-$ is similar.
 \begin{equation*}
  \begin{split}
   & \GG^{b_t}_t\pair{\dfrac{d}{dt}\pair{\gcon^{\bismutt}_{x}} y^{b_t+}, z^{b_t+}}\\
   = & \dfrac{d}{dt}\sair{\GG^{b_t}_t\pair{\gcon^{\bismutt}_{x} y^{b_t+}, z^{b_t+}}} - \GG^{b_t}_t\pair{\gcon^{\bismutt}_{x} \sair{\dfrac{d}{dt}\pair{b_t + g_t}Y}, z^{b_t+}} - \GG^{b_t}_t\pair{\gcon^{\bismutt}_{x} y^{b_t+}, \dfrac{d}{dt}\pair{b_t + g_t}Z}\\
   = & g_t\pair{\dfrac{d}{dt} \nabla^{+\phi_t}_X Y, Z} + \sair{X Rc_t^{-\phi}\pair{Y, Z} - Rc_t^{-\phi}\pair{Y, \nabla^{+\phi_t}_X Z} - Rc_t^{-\phi} \pair{\nabla^{+\phi_t}_X Y, Z}}
  \end{split}
 \end{equation*}
 where the last equality follows from \eqref{eq:riemannlaxriccidecomposition} together with the fact that $\gcon^{\bismutt}$ preserves $\GG_t^{b_t}$. Next, \eqref{eq:bismutriccieigen} implies that
 \begin{equation*}
  \begin{split}
   & \GG^{b_t}_t\pair{ \sair{\Ric^\bismutt_t, \gcon^{\bismutt}_{x}} y^{b_t +}, z^{b_t +} } = \GG^{b_t}_t\pair{ \Ric^\bismutt_t \gcon^{\bismutt}_{x} y^{b_t +} - \gcon^{\bismutt}_{x} \Ric^{\phi, \mc B}_t y^{b_t +}, z^{b_t +} }\\
   = & \Rc^{\bismutt}_t\pair{\gcon^{\bismutt}_{x} y^{b_t +}, z^{b_t +}} - X \Rc^{\bismutt}_t\pair{ y^{b_t +}, z^{b_t +}} + \Rc^{\bismutt}_t\pair{y^{b_t +}, \gcon^{\bismutt}_{x} z^{b_t +}}\\
   = &  - X Rc_t^{+\phi_t}\pair{Y, Z} + Rc_t^{+\phi_t}\pair{Y, \nabla^{+\phi_t}_X Z} + Rc_t^{+\phi_t}\pair{\nabla^{+\phi_t}_X Y, Z}
  \end{split}
 \end{equation*}
 Combining the results above and applying \eqref{eq:riemannlaxriccidecomposition} again lead to
 \begin{equation*}
  \begin{split}
   & \GG^{b_t}_t\pair{\pair{\delta_{\Ric^\bismut}\gcon^{\bismutt}_{x}} y^{b_t\pm}, z^{b_t\pm}}\\
   = & g_t\pair{\dfrac{d}{dt} \nabla^{+\phi_t}_X Y, Z} - \sair{X \pair{\dfrac{d}{dt} g_t}\pair{Y, Z} - \pair{\dfrac{d}{dt} g_t}\pair{Y, \nabla^{+\phi_t}_X Z} - \pair{\dfrac{d}{dt} g_t} \pair{\nabla^{+\phi_t}_X Y, Z}} \\
  = & - g_t\pair{Y, \dfrac{d}{dt} \nabla^{+\phi_t}_X Z}
  \end{split}
 \end{equation*}
 where the last equality follows from the fact that $\nabla^{+\phi_t}$ preserves $g_t$ for all $t$.
\end{proof}

Consider a smooth family of generalized almost Hermitian structures $\pair{M, \gamma; \GG_t^{b_t}, \JJ^t}$, where $\GG_t^{b_t}$ is a solution to the Bismut-Ricci Lax flow \eqref{eq:bismutlaxricci}. By Lemma \ref{lemma:commutinglaxdifferential}, the $\Ric^{\bismut}$-differential $\delta_{\Ric^\bismut}\JJ^t$ preserves the $\GG_t^{b_t}$-eigenbundles. Following computations similar to those in Proposition \ref{prop:riccidifferentialbismut} gives:
\begin{equation}\label{eq:riccidifferentialJ}
 \GG_t^{b_t}\pair{\pair{\delta_{\Ric^{\phi, \mc B}}\JJ^t} y^{b\pm}, z^{b\pm}} = - g_t\pair{Y, \pair{\dfrac{d}{dt} I_\pm^t} Z}
\end{equation}
Set $\JJ_-^t := \GG_t^{b_t} \JJ^t$ then
\[\delta_{\Ric^{\phi, \mc B}}\JJ_-^t = \pair{\delta_{\Ric^{\phi, \mc B}}\GG_t^{b_t} } \JJ^t + \GG_t^{b_t}\pair{\delta_{\Ric^{\phi, \mc B}} \JJ^t } = \GG_t^{b_t}\pair{\delta_{\Ric^{\phi, \mc B}} \JJ^t }\]
It implies that $\JJ^t$ solves the following Lax flow \eqref{eq:bismutlaxricciJ}, in which case as does $\JJ_-^t$, iff $I_\pm^t$ are constant almost complex structures:
\begin{equation}\label{eq:bismutlaxricciJ}
 \dfrac{d}{dt} \JJ^t = \sair{\Ric^\bismutt, \JJ^t}
\end{equation}
In terms of the $2$-forms $\omega_\pm$, \eqref{eq:bismutlaxricciJ} is equivalent to the following simultaneous almost Hermitian Ricci flows:
\begin{equation}\label{eq:bismutlaxricciomega}
 \pair{\dfrac{d}{dt}\omega_\pm^t}\pair{X, Y} 
 = 2 Rc\pair{X, I_\pm^t Y} - \half \phi_t^2\pair{X, I_\pm^t Y}
\end{equation}
where $X, Y \in \se{TM}$. It has the classical K\"ahler-Ricci flow as a special case.

\begin{example}\label{ex:kahlerricci}
Take a family of classical K\"ahler structures $\pair{g_t, I, \omega_t}$, set
\[\JJ^t = \begin{bmatrix}I & 0\\ 0 & -I^*\end{bmatrix}, \JJ_-^t = \begin{bmatrix}0 & -\omega_t\inverse\\ \omega_t & 0\end{bmatrix}\]
The generalized Bismut connections, as well as the $\pair{\phi, b}$-connections, are simply the lift of the Levi-Civita connections for $g_t$. Let $X, Y \in TM$, since $\phi = 0$, 
\eqref{eq:bismutlaxricciomega} becomes
\[\pair{\dfrac{d}{dt}\omega_t}\pair{X, Y} = 2Rc_t\pair{X, IY} = - 2 \rho_t\pair{X, Y}\]
where $\rho_t\pair{X, Y} = Rc_t\pair{IX, Y}$ is the Ricci form. Thus \eqref{eq:bismutlaxricciJ} recovers exactly the equation for the classical K\"ahler-Ricci flow.
\end{example}

\end{document}